\documentclass{article}
\usepackage{amsfonts,amsthm}
\usepackage{authblk}
\usepackage{newpxtext,newpxmath}
\usepackage{pdflscape}
\usepackage{multirow}
\usepackage{multicol}
\usepackage{tikz}
\usetikzlibrary{arrows}
\usepackage{geometry}
\usepackage{relsize}
\usepackage{xfrac}
\usepackage{arydshln}
\makeatletter
\newcommand{\oset}[3][0ex]{%
  \mathrel{\mathop{#3}\limits^{
    \vbox to#1{\kern-2\ex@
    \hbox{$\scriptstyle#2$}\vss}}}}
\makeatother
\usepackage{mathrsfs}

\theoremstyle{plain}% default
\newtheorem{theorem}{Theorem}[section]

\newtheorem{proposition}[theorem]{Proposition}
\newtheorem{corollary}[theorem]{Corollary}

\theoremstyle{definition}
\newtheorem{definition}{Definition}[section]

\newtheorem{example}{Example}[section]

\theoremstyle{remark}
\newtheorem{remark}{Remark}[section]

\title{On the variety of Lie algebras endowed with complex structures: degenerations and deformations.}

\makeatletter
\def\blfootnote{\xdef\@thefnmark{}\@footnotetext}
\makeatother

\makeatletter
\renewcommand\@date
{{%

\vspace{-1.5cm}%
\large\centering

  \begin{tabular}{@{}c@{}}
    Edison Alberto \textsc{Fern\'andez-Culma} \textsuperscript{1,3,$\ddagger$}\blfootnote{${}^{\ddagger}$ Partially supported by CONICET and SECyT-UNC: 33620180100523CB} \\
    \texttt{\small efernandez@famaf.unc.edu.ar}
  \end{tabular}%

  and

  \begin{tabular}{@{}c@{}}
   Nadina  \textsc{Rojas} \textsuperscript{1,2,3,$\ddagger$} \\
     \texttt{ \small nadina.rojas@unc.edu.ar}
  \end{tabular}

  \textsuperscript{1}{\small Centro de Investigaci\'on y Estudios en Matem\'atica de C\'ordoba, CONICET}\par
  \textsuperscript{2}{\small Facultad de Ciencias Exactas, F\'isicas y Naturales, UNC}\par
  \textsuperscript{3}{\small Facultad de Matem\'atica, Astronom\'ia, F\'isica y Computaci\'on, UNC}\par
  \textsuperscript{ }{\small C\'ordoba, Argentina}\par
}}
\makeatother

\author{}

\begin{document}

\maketitle

\blfootnote{2020 \textit{Mathematics Subject Classification:} 17B99. \textit{Key words and phrases:} Lie algebras, Orbit closure problem, Complex structures, Left-invariant Hermitian structures.}

\begin{abstract}
We study the space of Lie algebras equipped with left-invariant complex structures, \( \mathcal{L}_{ J_{\tiny{\mbox{cn}}} }(\mathbb{R}^{2n}) \), with particular attention to their degenerations and deformations. To this end, we identify certain invariants that remain well-behaved under degenerations while preserving the complex structure. These concepts are then applied to the four-dimensional case. Additionally, we explore applications to the study of left-invariant Hermitian structures on Lie groups, and we discuss some aspects of the deformation theory within \( \mathcal{L}_{ J_{\tiny{\mbox{cn}}} }(\mathbb{R}^{2n}) \).
\end{abstract}

\section{Introduction}{$ $}

Invariant complex structures on Lie groups offer a powerful framework for studying their geometric and algebraic properties, playing a key role in both differential geometry and mathematical physics. The interaction between algebra and geometry in Lie groups makes the study of these structures particularly challenging, but also provides valuable insights that serve as a stepping stone to understanding more intricate phenomena in complex geometry and the classification of special geometric structures. Research on complex structures in Lie algebras
remains a dynamic field, with significant advancements in their classification and in understanding their connections to complex and symplectic geometry (see, for example, \cite{CaoZheng,FinoMainenti,LatorreUgarteVillacampa,ABDGH} and the references therein).

On the other hand, the study of degenerations of algebraic structures---particularly associative and Lie algebras---has its origins in applications within mathematical physics, geometry, and algebra. Over time, this research has proven essential across various fields. Investigating the degeneration of these structures often leads to the discovery of new invariants and classification results. For instance, the role of closed orbits in reductive representations is a key example, and in the realm of geometry, John Milnor's groundbreaking work on left-invariant metrics on Lie groups provides an important context. In his work, degenerations of Lie algebras are implicitly used to demonstrate the existence of Riemannian metrics with specific curvature properties (see, for example, \cite[Theorem 2.5]{Milnor}). In recent years, this area of research has garnered increasing attention in algebra (see, e.g., \cite[Section 2]{Kaygorodov}). In \cite{FCRojitas}, we explore the degenerations of Lie algebras endowed with symplectic structures, offering new insights into their algebraic and geometric properties, as well as their significance for understanding the almost-K\"{a}hler geometry of Lie groups.

In this work, we investigate the variety of Lie algebras equipped with a complex structure, denoted by \( \mathcal{L}_{ J_{\tiny{\mbox{cn}}} }(\mathbb{R}^{2n}) \), with a particular focus on their deformations and degenerations. The deformation theory of algebraic structures, which traditionally examines how an algebra can be perturbed into a nearby structure, naturally extends to the setting where the integrability of the complex structure is preserved. On the other hand, degenerations provide a complementary perspective by describing how families of Lie algebras with complex structures collapse into limiting objects, shedding light on fundamental algebraic and geometric invariants.

The paper is structured as follows. In Section 2, we introduce the necessary background on the variety of Lie algebras endowed with complex structures, the notion of degeneration, and key invariants that remain well behaved under this process. Section 3 is dedicated to the classification of degenerations in dimension 4, where we illustrate the classification through concrete examples. Finally, in Section 5, we explore applications to the Hermitian geometry of Lie groups and discuss relevant observations on deformations within \( \mathcal{L}_{ J_{\tiny{\mbox{cn}}} }(\mathbb{R}^{2n}) \).

\section{Preliminaries}

\begin{definition}[Complex structure on Lie algebras]
Let $\mathfrak{g}$ be a finite-dimensional real Lie algebra. A \emph{complex structure} on $\mathfrak{g}$ is a linear transformation $J: \mathfrak{g} \to \mathfrak{g}$ satisfying
$J^2 = -\operatorname{Id}_{\mathfrak{g}},$ with the property that the Nijenhuis tensor
\begin{eqnarray}\label{nijenhuis}
N_J(X,Y) = [X,Y] - [JX,JY] + J [JX,Y] + J [X,JY]
\end{eqnarray}
vanishes for all \(X,Y \in \mathfrak{g}\).
\end{definition}

\begin{definition}[Bi-invariant and Abelian Complex Structures]
Let $\mathfrak{g}$ be a finite-dimensional real Lie algebra and $J:\mathfrak{g}\to\mathfrak{g}$ a linear map satisfying $J^2 = -\operatorname{Id}_{\mathfrak{g}}$.

For a fixed integer $k\in\mathbb{Z}$, suppose that
\begin{eqnarray}\label{definicionAbBi}
[X,Y] - J^k\bigl[JX,J^{k+1}Y\bigr] &=& 0
\end{eqnarray}
for all $X,Y\in\mathfrak{g}$. Then, $J$ is called an \emph{abelian complex structure} if $k$ is even, and a \emph{bi-invariant complex structure} if $k$ is odd.

In particular, when $k$ is even, equation \eqref{definicionAbBi} reduces to
$$
[X,Y] - [JX,JY] = 0.
$$
On the other hand, if $k$ is odd, it is equivalent to
$$
[JX,Y] = J[X,Y].
$$
\end{definition}

A direct computation shows that condition \eqref{definicionAbBi} necessarily implies \eqref{nijenhuis}, ensuring that any bi-invariant or abelian complex structure is indeed a complex structure.

\begin{definition}[Equivalence]
Let \( (\mathfrak{g}_1, J_1) \) and \( (\mathfrak{g}_2, J_2) \) be two Lie algebras with complex structures. The pairs \( (\mathfrak{g}_1, J_1) \) and \( (\mathfrak{g}_2, J_2) \) are called \textit{holomorphically isomorphic} if there exists a Lie algebra isomorphism \( g: \mathfrak{g}_1 \to \mathfrak{g}_2 \) such that \( g J_1 = J_2 g \). In particular, when \( (\mathfrak{g}_1, J_1) \)= \( (\mathfrak{g}_2, J_2)\), the isomorphism \( g \) is called a \textit{holomorphic automorphism}.
\end{definition}

\subsection{The variety of Lie algebras endowed with complex structures}

To aid in the following discussion, we introduce some useful notation. Let \( V \) and \( W \) be finite-dimensional real vector spaces. Define \( L^{k}(V; W) \) as the space of all multilinear maps from \( V^k \), the \( k \)-fold Cartesian product of \( V \), to \( W \). Moreover, let \( C^{k}(V;W) \subseteq L^{k}(V;W) \) be the subspace consisting of \textit{alternating} multilinear maps, i.e., maps \( \mu \in L^{k}(V;W) \) such that for any permutation \( \sigma \in S_k \), the following identity holds:

\[
\mu(X_{\sigma(1)}, X_{\sigma(2)}, \ldots , X_{\sigma(k)}) = \operatorname{sign}(\sigma) \, \mu(X_1, X_2, \ldots, X_k).
\]

Next, we describe the action of the general linear group of \( V \), denoted \( \operatorname{GL}(V) \), on \( L^{k}(V;V) \) and \( L^{k}(V; \mathbb{R}) \) through \textit{change of basis}. Specifically, for any \( g \in \operatorname{GL}(V) \), \( \mu \in L^{k}(V;V) \), and \( \alpha \in L^{k}(V;\mathbb{R}) \), the actions are defined as follows:

\[
g\cdot\mu(\cdot, \ldots, \cdot) := g\mu(g^{-1}\cdot, \ldots, g^{-1}\cdot),
\]
and
\[
g\cdot\alpha(\cdot, \ldots, \cdot) := \alpha(g^{-1}\cdot, \ldots, g^{-1}\cdot).
\]

In the study of complex structures on vector spaces, recall that if \( J: V \to V \) is a linear map satisfying $J^2 = -\operatorname{Id}_V$, then by the \textbf{real Jordan form theorem} there exists an ordered basis $
\{e_1, \ldots, e_n, e_{n+1}, \ldots, e_{2n}\}$ of \( V \) such that
\[
J e_i = e_{n+i} \quad \text{for all } i = 1, \ldots, n.
\]
Therefore, from now on we assume that \( V = \mathbb{R}^{2n} \) and introduce the \textit{canonical complex structure} \( J_{\tiny{\mbox{cn}}} \) on \( \mathbb{R}^{2n} \), defined by
\[
J_{\tiny{\mbox{cn}}} e_i := e_{n+i}, \quad J_{\tiny{\mbox{cn}}} e_{n+i} := -e_i, \quad \text{for } i = 1, \ldots, n.
\]
Now, let us denote by
\[
\operatorname{GL}(\mathbb{R}^{2n}, J_{\tiny{\mbox{cn}}} ) := \{ g \in \operatorname{GL}(\mathbb{R}^{2n}) : g J_{\tiny{\mbox{cn}}} g^{-1} = J_{\tiny{\mbox{cn}}} \}
\]
the Lie subgroup of \( \operatorname{GL}(\mathbb{R}^{2n}) \) that preserves \( J_{\tiny{\mbox{cn}}} \).

The Lie algebra of \( \operatorname{GL}(\mathbb{R}^{2n}, J_{\tiny{\mbox{cn}}} ) \), which we denote by \( \mathfrak{gl}(\mathbb{R}^{2n}, J_{\tiny{\mbox{cn}}}) \), is the Lie subalgebra of \( \mathfrak{gl}(\mathbb{R}^{2n}) \) given by
\[
C( J_{\tiny{\mbox{cn}}} ) = \{ T \in \mathfrak{gl}(\mathbb{R}^{2n}) : T J_{\tiny{\mbox{cn}}} - J_{\tiny{\mbox{cn}}} T = 0 \}.
\]
Observe that \( C( J_{\tiny{\mbox{cn}}} ) \) is, in fact, an associative subalgebra of \( L(\mathbb{R}^{2n}; \mathbb{R}^{2n}) \) (with the usual composition of functions) and is naturally isomorphic to \( L^{1}(\mathbb{C}^n; \mathbb{C}^n) \). The matrix representation of any linear transformation in \( C( J_{\tiny{\mbox{cn}}} ) \), with respect to the canonical basis of \( \mathbb{R}^{2n} \), is a block matrix of the form
\[
M = \begin{pmatrix} A & -B \\ B & A \end{pmatrix},
\]
with \( A, B \in \operatorname{M}(n, \mathbb{R}) \). An isomorphism between these associative algebras is given by
\begin{eqnarray}\label{isomorphism}
\mathcal{I}: C(  J_{\tiny{\mbox{cn}}} ) \to L(\mathbb{C}^n; \mathbb{C}^n), \quad M \mapsto A + \sqrt{-1}B,
\end{eqnarray}
where we identify elements of \( L(\mathbb{C}^n; \mathbb{C}^n) \) with their matrices with respect to the canonical basis of \( \mathbb{C}^n \). Since the map $\mathcal{I}$ satisfies that
$\operatorname{det}(M) = | \operatorname{det} (\mathcal{I}(M)) |^2$, it follows that $\mathcal{I}$ restricts to a Lie group isomorphism between $ \operatorname{GL}(\mathbb{R}^{2n}, J_{\tiny{\mbox{cn}}} ) $ and
$ \operatorname{GL}(\mathbb{C}^{n})  $.

In what follows, we denote by $C^{2}_{ J_{\tiny{\mbox{cn}}}  }(\mathbb{R}^{2n};\mathbb{R}^{2n})$ the vector space $\left\{\mu \in C^{2}(\mathbb{R}^{2n};\mathbb{R}^{2n})  : N_{J_{\tiny{\mbox{cn}}}}(\mu) =0 \right\}$, where
$$
N_{ J_{\tiny{\mbox{cn}}} }(\mu)(X,Y) = \mu(X,Y) - \mu(J_{\tiny{\mbox{cn}}} X,J_{\tiny{\mbox{cn}}} Y) + J_{\tiny{\mbox{cn}}}  \mu(J_{\tiny{\mbox{cn}}} X,Y) + J_{\tiny{\mbox{cn}}} \mu(X,J_{\tiny{\mbox{cn}}} Y)
$$
From this space, we single out those bilinear maps that endow \(\mathbb{R}^{2n}\) with a Lie algebra structure. More precisely, we define
\[
\mathcal{L}_{ J_{\tiny{\mbox{cn}}}  }(\mathbb{R}^{2n}) \coloneqq \left\{ \mu \in C^{2}_{ J_{\tiny{\mbox{cn}}}  }(\mathbb{R}^{2n}; \mathbb{R}^{2n}) : \operatorname{Jac}(\mu)=0 \right\},
\]
or equivalently,
\[
\mathcal{L}_{ J_{\tiny{\mbox{cn}}}  }(\mathbb{R}^{2n}) = \left\{ \mu \in C^{2}(\mathbb{R}^{2n}; \mathbb{R}^{2n}) : N_{ J_{\tiny{\mbox{cn}}}  }(\mu)=0 \text{ and } \operatorname{Jac}(\mu)=0 \right\}.
\]

Here, \(\operatorname{Jac}(\mu)\) denotes the Jacobiator---i.e., the alternating trilinear map defined by
\[
\operatorname{Jac}(\mu)(X_1,X_2,X_3) \coloneqq \sum_{\sigma \in S_3} \operatorname{sign}(\sigma)\,\mu\big(\mu(X_{\sigma(1)},X_{\sigma(2)}),X_{\sigma(3)}\big),
\]
so that the condition \(\operatorname{Jac}(\mu)=0\) is equivalent to the Jacobi identity.

We refer to \(\mathcal{L}_{  J_{\tiny{\mbox{cn}}}  }(\mathbb{R}^{2n})\) as the algebraic variety of Lie algebra structures on \(\mathbb{R}^{2n}\) endowed with a complex structure. This subset is invariant under the natural action of \( \operatorname{GL}(\mathbb{R}^{2n},  J_{\tiny{\mbox{cn}}} ) \), and the quotient $\mathcal{L}_{  J_{\tiny{\mbox{cn}}}  }(\mathbb{R}^{2n}) \big/ \operatorname{GL}(\mathbb{R}^{2n},  J_{\tiny{\mbox{cn}}} )$
parametrizes, up to holomorphic isomorphism, the family of all \(2n\)-dimensional Lie algebras endowed with a complex structures.

Note that the Equation (\ref{definicionAbBi}) defines two closed subvarieties of $\mathcal{L}_{ J_{\tiny{\mbox{cn}}}  }(\mathbb{R}^{2n})$, which are also $ \operatorname{GL}(\mathbb{R}^{2n},  J_{\tiny{\mbox{cn}}} )$-invariants.

We consider the Hausdorff vector topology  on \( C^{2}_{  J_{\tiny{\mbox{cn}}}  }(\mathbb{R}^{2n}; \mathbb{R}^{2n}) \) and its induced topology on \( \mathcal{L}_{  J_{\tiny{\mbox{cn}}}  }(\mathbb{R}^{2n}) \), providing the framework to define degeneration. The closure of any subset \( S \) of \( C^{2}_{  J_{\tiny{\mbox{cn}}}  }(\mathbb{R}^{2n}; \mathbb{R}^{2n}) \) in this topology is denoted by \( \overline{S} \).

\begin{definition}[Degeneration]
Let $\mu$ and $\lambda$ be Lie algebra structures in $\mathcal{L}_{J_{\tiny{\mbox{cn}}}}(\mathbb{R}^{2n})$. We say that the Lie algebra with complex structure \((\mathbb{R}^{2n}, \mu, J_{\tiny{\mbox{cn}}})\) \textit{degenerates to} \((\mathbb{R}^{2n}, \lambda, J_{\tiny{\mbox{cn}}})\) with respect to $\operatorname{GL}(\mathbb{R}^{2n}, J_{\tiny{\mbox{cn}}})$ if $\lambda$ belongs to the closure $\overline{\operatorname{GL}(\mathbb{R}^{2n}, J_{\tiny{\mbox{cn}}}) \cdot \mu}$. Furthermore, if $\lambda \notin \operatorname{GL}(\mathbb{R}^{2n}, J_{\tiny{\mbox{cn}}} ) \cdot \mu$, we say that the degeneration is \textit{proper}.
\end{definition}

Thus, the central problem we address is the following: given two Lie brackets $\mu$ and $\lambda$ in $\mathcal{L}_{ J_{\tiny{\mbox{cn}}} }(\mathbb{R}^{2n})$, we seek an elementary and verifiable criterion to determine whether $\lambda$ belongs to the orbit closure $\overline{ \operatorname{GL}(\mathbb{R}^{2n}, J_{\tiny{\mbox{cn}}}) \cdot \mu }$. In the classical study of degenerations of (Lie) algebras, where one analyzes the closure of $\operatorname{GL}(m, \mathbb{R})$-orbits, it is common practice to employ well-established algebraic invariants to address this question (see, for example, \cite[Theorem 1]{nesterenko}). In our setting, however, it is necessary to introduce invariants that remain stable under degeneration with respect to the action of $\operatorname{GL}(\mathbb{R}^{2n}, J_{\tiny{\mbox{cn}}})$. A key invariant of this type arises from the well-known \textit{closed orbit lemma} of Armand Borel (see \cite[Proposition 15.4]{Borel1}) and its adaptation to the real case by Michael Jablonski.

\begin{proposition}[{\cite[Proposition 3.2-(b)]{Jablonski1}}]\label{orden}
Let $G$ be a real reductive algebraic group acting linearly and rationally on a real vector space $V$. Then, the boundary of a $G$-orbit, given by $\partial(G \cdot v) = \overline{G\cdot v} \setminus G\cdot v$, consists of $G$-orbits of strictly smaller dimension.
\end{proposition}

If we denote by $G_p := \{ g \in G \mid g \cdot p = p \}$ the \textit{isotropy group} (or \textit{stabilizer}) of $p$, whose Lie algebra $\mathfrak{g}_p$ consists of elements of $\mathfrak{g} = \operatorname{Lie}(G)$ that act trivially at $p$, then a standard result states that
\[
\dim(G\cdot p) = \dim(\mathfrak{g}) - \dim(\mathfrak{g}_p).
\]
In our setting, when $G = \operatorname{GL}(\mathbb{R}^{2n}, J_{\tiny{\mbox{cn}}})$ and $\mu \in \mathcal{L}_{ J_{\tiny{\mbox{cn}}} }(\mathbb{R}^{2n})$, the stabilizer $G_{\mu}$ corresponds to the group of holomorphic automorphisms of $(\mathbb{R}^{2n}, \mu, J_{\tiny{\mbox{cn}}})$. Regarding $\operatorname{Lie}(G_{\mu})$, we introduce the following definition:

\begin{definition}
Let $(\mathfrak{g}, J)$ be a Lie algebra endowed with a complex structure. A \textit{holomorphic derivation} of $(\mathfrak{g}, J)$ is a linear map $D: \mathfrak{g} \to \mathfrak{g}$ satisfying:
\begin{itemize}
    \item (\textit{Derivation property}) $D [X,Y] = [D X, Y] + [X,D Y]$ for all $X,Y \in \mathfrak{g}$,
    \item (\textit{Commutativity with} $J$) $D J = J D$.
\end{itemize}
We denote by $\operatorname{Der}_{J}(\mathfrak{g})$ the real vector space of all holomorphic derivations of $(\mathfrak{g}, J)$.
\end{definition}

Using this definition and the orbit structure described above, an immediate consequence of Proposition \ref{orden} is that the dimension of the algebra of holomorphic derivations of a Lie algebra with a complex structure acts as an \textit{obstruction} to its degeneration.

\begin{corollary}\label{dercomp}
If a Lie algebra with complex structure $(\mathbb{R}^{2n}, \mu, {J_{\tiny{\mbox{cn}}}} )$ degenerates properly to $(\mathbb{R}^{2n}, \lambda, {J_{\tiny{\mbox{cn}}}} )$ with respect to $J_{\tiny{\mbox{cn}}}$, then
\[
\dim \operatorname{Der}_{J_{\tiny{\mbox{cn}}}}(\mathbb{R}^{2n}, \mu ) < \dim \operatorname{Der}_{J_{\tiny{\mbox{cn}}}}(\mathbb{R}^{2n}, \lambda ).
\]
\end{corollary}

Another key consequence of Proposition \ref{orden} is that the notion of degeneration for Lie algebras with complex structures induces a \textit{partial order} on the orbit space $
\mathcal{L}_{J_{\tiny{\mbox{cn}}}}(\mathbb{R}^{2n}) / \operatorname{GL}(\mathbb{R}^{2n}, J_{\tiny{\mbox{cn}}}) $ given by $
\operatorname{GL}(\mathbb{R}^{2n}, J_{\tiny{\mbox{cn}}}) \cdot \mu \geq \operatorname{GL}(\mathbb{R}^{2n}, J_{\tiny{\mbox{cn}}}) \cdot \lambda $ if $(\mathbb{R}^{2n}, \mu, {J_{\tiny{\mbox{cn}}}} )$ degenerates to $(\mathbb{R}^{2n}, \lambda, {J_{\tiny{\mbox{cn}}}} )$ with respect to $\operatorname{GL}(\mathbb{R}^{2n}, J_{\tiny{\mbox{cn}}})$.

The following observation provides an additional tool for studying orbit closures, and can be interpreted as stating that degenerations are determined by the orbit closures of Borel subgroups.

\begin{proposition}[{\cite[Proposition 1.7]{GrunewaldOHalloran}}]\label{Borel}
Let $\operatorname{G} \subseteq \operatorname{GL}(n,\mathbb{C})$ be a reductive complex algebraic group acting rationally on an algebraic set $\mathcal{Z}$, and let $\operatorname{K} \operatorname{B}$ be an \textit{Iwasawa decomposition} of $\operatorname{G}$. Then, for all $p \in \mathcal{Z}$,
\[
\overline{\operatorname{G}\cdot p} = \operatorname{K} \cdot \overline{\operatorname{B}\cdot p}.
\]
\end{proposition}

The proof of this result can be easily extended to the context of continuous representations of real reductive Lie groups.

Using the Lie group isomorphism provided in (\ref{isomorphism}) together with the standard Iwasawa decomposition of \(\operatorname{GL}(\mathbb{C}^n)\) (see, e.g., \cite[Chapter 26]{bump}), one can derive an Iwasawa decomposition for \(\operatorname{GL}(\mathbb{R}^{2n}, J_{\text{cn}})\). By identifying the corresponding linear maps with their matrix representations relative to the canonical bases, we obtain $\operatorname{GL}(\mathbb{R}^{2n}, J_{\text{cn}}) = \operatorname{K}\operatorname{B},$
where
\[
\operatorname{K} = \Biggl\{ \begin{pmatrix} X & -Y \\ Y & X \end{pmatrix} : X + iY \in \operatorname{U}(n) \Biggr\},
\]
\[
\operatorname{A} = \Biggl\{ \operatorname{diag}(t_1, \dots, t_n, t_1, \dots, t_n) : t_i > 0 \Biggr\},
\]
\[
\operatorname{N} = \Biggl\{ \begin{pmatrix} X & -Y \\ Y & X \end{pmatrix} :
X = \begin{bmatrix} 1 & \cdots & 0 \\ \vdots & \ddots & \vdots \\ * & \cdots & 1 \end{bmatrix},\quad
Y = \begin{bmatrix} 0 & \cdots & 0 \\ \vdots & \ddots & \vdots \\ * & \cdots & 0 \end{bmatrix}
\Biggr\},
\]
and we define \(\operatorname{B} = \operatorname{A}\operatorname{N}\).

\subsubsection{Invariants}\label{invariante}

The iterative application of an algebra product to construct multilinear maps---and the study of their invariance under isomorphisms---is a well-established theme in (Lie) algebra theory, with numerous classical results characterizing their structure and properties. In particular, multilinear Lie polynomials and their traces serve as fundamental tools in this context. We recall these notions as a foundation for exploring analogous invariants within our framework.

Given a (Lie) algebra \((\mathbb{R}^{m}, \mu)\), one can construct various multilinear maps and forms that remain invariant under isomorphisms. A representative example is a \textit{multilinear Lie polynomial} \(P_{\mu}\) in \(k\) variables \(X_1, X_2, \dots, X_k\) (where \(k \geq 2\)), which belongs to \(L^{k}(\mathbb{R}^{m}; \mathbb{R}^{m})\). It is well established that \(P_{\mu}\) can be expressed as
\[
P_{\mu}(X_1, X_2, \dots, X_k) = \sum_{ \substack{\sigma \in \operatorname{S}_{k}\\ \sigma(1) = 1 }} a_{\sigma} \, \mu(\dots \mu(\mu(X_{\sigma(1)},X_{\sigma(2)} ), X_{\sigma(3)}) \dots X_{\sigma(k)})
\]
as documented in \cite[\S 5.6.2.]{Reutenauer} and \cite[\S 4.8.1.]{Bahturin}.

An important associated notion is the \textit{\(i\)-th trace of \(P_{\mu}\)} (see, for instance, \cite[Appendix B]{Lee2}). This trace operation, denoted by \(\operatorname{tr}_{i}\), is a mapping \(L^{k+1}(\mathbb{R}^{m};\mathbb{R}^{m}) \to L^{k}(\mathbb{R}^{m}; \mathbb{R})\) that assigns to \((\operatorname{tr}_{i} \mu)(X_1, \dots, X_{k})\) the trace of the linear operator
\[
\mu(X_1, \dots, \underbrace{\square}_{\text{\tiny \(i\)-th entry}}, \dots , X_{k}) : \mathbb{R}^{m} \to \mathbb{R}^{m}.
\]

Furthermore, linear combinations or the tensor product of multilinear forms preserved under isomorphisms serve as mechanisms to construct novel algebraic \textit{invariants}.

To illustrate these ideas, consider the non-associative, non-commutative polynomial \(P(X,Y,Z) = X \cdot (Y \cdot Z)\) and the map \(\mathcal{P}: L^{2}(\mathbb{R}^{m};\mathbb{R}^{m}) \to  L^{3}(\mathbb{R}^{m};\mathbb{R}^{m})\) defined by \(\mathcal{P}(\mu)(v_1,v_2,v_3) = \mu(v_1, \mu(v_2,v_3))\). It can be verified that \(\mathcal{P}\) is a \(\operatorname{GL}(\mathbb{R}^{m})\)-equivariant continuous map, meaning that \(\mathcal{P}(g \cdot \mu) = g \cdot \mathcal{P}(\mu)\). This property facilitates the natural construction of the transformations \(\operatorname{tr}_1 \circ \mathcal{P}, \operatorname{tr}_3 \circ \mathcal{P}: L^{2}(\mathbb{R}^{m} ;\mathbb{R}^{m}) \to L^{2}(\mathbb{R}^{m};\mathbb{R})\), both of which inherit \(\operatorname{GL}( \mathbb{R}^{m} )\)-equivariance and continuity.

When \(\mathfrak{g} = (\mathbb{R}^{m},\mu)\) is a Lie algebra, it follows that \((\operatorname{tr}_1 \circ \mathcal{P}) \mu\) vanishes identically, whereas \((\operatorname{tr}_3 \circ \mathcal{P})\mu\) corresponds to the classical \textit{Cartan-Killing form} of \(\mathfrak{g}\). Moreover, the expression \((\operatorname{tr}_3 \circ \mathcal{P})\mu +  c \, (\operatorname{tr}_2\mu) \otimes (\operatorname{tr}_2\mu)\), for some scalar \(c\), defines the so-called \textit{modified Cartan-Killing form}, as formulated in \cite[\S IV.]{nesterenko}.

Beyond these classical constructions, the presence of an additional structure, such as a complex structure \(J\), naturally leads to new families of multilinear maps that incorporate this extra geometric data. To illustrate this, let us consider a particular example: let $P$ be the following nonassociative, noncommutative monomial in $k$ variables
\[
P(X_1, X_2, \dots, X_k) = (\dots (((X_1 \cdot X_2) \cdot X_3) \dots) \cdot X_k),
\]
Now, we modify $P$ by inserting powers of \(J_{\tiny{\mbox{cn}}}\) at each stage and define the \textit{modified nonassociative monomial}
\[
Q(X_1, X_2, \dots, X_k) = J_{\tiny{\mbox{cn}}}^{a_{k-1}}\Big( \dots J_{\tiny{\mbox{cn}}}^{a_3} \Big( J_{\tiny{\mbox{cn}}}^{a_2}\Big( J_{\tiny{\mbox{cn}}}^{a_1}\big( J_{\tiny{\mbox{cn}}}^{b_1}X_1 \cdot J_{\tiny{\mbox{cn}}}^{b_2}X_2\big) \cdot J_{\tiny{\mbox{cn}}}^{b_3}X_3\Big) \dots J_{\tiny{\mbox{cn}}}^{b_k}X_k \Big),
\]
where \(a_1, a_2, \dots, a_{k-1}, b_1, b_2, \dots, b_k\) are integers and the powers of \( J_{\tiny{\mbox{cn}}} \) act on each variable as well as on each intermediate product.

If \(\mu: \mathbb{R}^{2n} \times \mathbb{R}^{2n} \to \mathbb{R}^{2n}\) is a bilinear map defining a product on \(\mathbb{R}^{2n}\), then substituting every instance of the formal product \(\cdot\) in \(Q\) with the product given by \(\mu\) results in the mapping
\[
Q_{\mu}(X_1, X_2, \dots, X_k).
\]
Due to the bilinearity of \(\mu\), this construction yields a well-defined \(k\)-linear map on $\mathbb{R}^{2n}$. In this way, we define a \(\operatorname{GL}(\mathbb{R}^{2n},J)\)-equivariant continuous map
\[
\mathcal{Q}: C^{2}_{ J_{\tiny{\mbox{cn}}} }(\mathbb{R}^{2n};\mathbb{R}^{2n}) \rightarrow L^{k}(\mathbb{R}^{2n};\mathbb{R}^{2n}),
\]
given by \(\mu \mapsto Q_{\mu}\).

With this motivation in mind, constructing invariants of Lie algebras with complex structures becomes straightforward. For instance, we define the following function: let \(c_0, \dots, c_7 \in \mathbb{R}\) be arbitrary constants, and consider
\[
\varphi_{\alpha_{0},\ldots,\alpha_{7}} : C^{2}_{{J_{\tiny{\mbox{cn}}}}}(\mathbb{R}^{2n};\mathbb{R}^{2n}) \rightarrow L^2(\mathbb{R}^{2n},\mathbb{R}^{2n})
\]
given by
\begin{eqnarray}\label{funcioninv}
\varphi_{\alpha_{0},\ldots,\alpha_{7}} (\mu)(X,Y) &=&
\begin{array}{l}
\alpha_{ 0}\mu(X,Y) +\\
\alpha_{ 1}J_{\tiny{\mbox{cn}}}\mu(X,Y)+
\alpha_{ 2}\mu(J_{\tiny{\mbox{cn}}}X,Y)+
\alpha_{ 3}\mu(X,J_{\tiny{\mbox{cn}}}Y)+\\
\alpha_{ 4}J_{\tiny{\mbox{cn}}}\mu(J_{\tiny{\mbox{cn}}} X,Y)+
\alpha_{ 5}J_{\tiny{\mbox{cn}}}\mu(X,J_{\tiny{\mbox{cn}}} Y)+
\alpha_{ 6}\mu(J_{\tiny{\mbox{cn}}} X,J_{\tiny{\mbox{cn}}} Y)+\\
\alpha_{ 7} J_{\tiny{\mbox{cn}}} \mu(J_{\tiny{\mbox{cn}}} X,J_{\tiny{\mbox{cn}}} Y).
\end{array}
\end{eqnarray}

We further define \(\psi_{\alpha,\beta}:=\varphi_{1,\alpha,\beta,\beta,0,\ldots,0}\) and \(\theta_{\beta}:=\varphi_{0,1,\beta,\beta,0,\ldots,0}\).

The function \(\varphi_{\alpha_{0},\ldots,\alpha_{7}}\) is continuous and \(\operatorname{GL}(\mathbb{R}^{2n}, J_{\tiny{\mbox{cn}}})\)-equivariant. Moreover, it assigns a second algebra structure on \(\mathbb{R}^{2n}\) to each algebra law in \(C^{2}_{J_{\tiny{\mbox{cn}}}}(\mathbb{R}^{2n};\mathbb{R}^{2n})\).

Furthermore, the functions \(\psi_{\alpha,\beta}\) and \(\theta_{\beta}\) assign skew-symmetric products to each element of \( C^{2}_{{J_{\tiny{\mbox{cn}}}}}(\mathbb{R}^{2n};\mathbb{R}^{2n}) \). In particular, we have the relation
\begin{eqnarray}\label{defotrivial}
\psi_{t,-t}(\mu) = \mu + t J \cdot \mu = g_{t} \cdot \mu,
\end{eqnarray}
where \( g_{t}^{-1} = \operatorname{Id} - t J \).

\subsection{Classification of complex structures on four-dimensional real Lie algebras.}

The classification of invariant complex structures on four-dimensional solvable Lie groups has been investigated by Gabriela Ovando \cite{ovando1} and Joanne Snow \cite{snow}. Their works offer a systematic approach to describing these structures, though in both cases, the results include specific examples where complete parameterizations of all possible complex structures were obtained.

In \cite[Proposition 3.2]{ovando2}, Gabriela Ovando provides a full classification of complex subalgebras within the complexification of four-dimensional solvable real Lie algebras that admit complex structures. Building on this work, \cite{LauretRodriguezValencia} presents a list of explicit complex structures, although it contains certain errors and omissions (for instance, $J_1$ and $J_2$ are not complex structures on the Lie algebra $\mathfrak{rr}'_{3,\gamma}\times \mathbb{R}$ and a complex structure on the Lie algebra $\mathfrak{r}'_2$ is missing).

Similarly, in \cite{Rezaei-AghdamSephid}, Adel Rezaei-Aghdam and Mehdi Sephid give a complete parametrization of complex structures on all four-dimensional real Lie algebras (not just the solvable ones). However, their classification also contains inaccuracies and redundancies (e.g., the complex structures $J_1$ and $J_3$, respectively $J_2$ and $J_4$, on $A_{4,10}$ are equivalent respectively, and a similar situation happens with $A_{3,8}$ and the two-parameter family of complex structures on $A_{4,12}$ is missing).

Following the notation used in \cite[Tables 1 and 2]{LauretRodriguezValencia}, we provide a corrected list of four-dimensional real Lie algebras endowed with a complex structure in the following proposition.

\begin{proposition} Any four-dimensional Lie algebra endowed with a complex structure is holomorphically isomorphic to exactly one of the following $(\mathbb{R}^{4}, \mu_{i}, J_{\tiny{\mbox{cn}}})$, where $\mu_{i}$ is defined in Table 1 or Table 2.

\begin{table}[h]
\centering
\begin{tabular}{ll}
$(\mathfrak{a}_{4}, J)$: & $\mu_{0}:=\, $ The four dimensional abelian Lie algebra\\
$(\mathfrak{rh}_{3}\times \mathbb{R},J)$: & $\mu_{1}:=\left\{ [e_1, e_3] = e_2 \right.$\\
%%%%%%%%%%%%%%%%%%%%%%%%%%%%%%%%%%%%%%%%%%%%%%%%%%%%%%%%%%%%%%%%%%%%%%%%%
$(\mathfrak{r}_{2}\times \mathfrak{r}_{2},J)$: & $\mu_{2}:=\left\{  [e_1, e_3] = e_3, [e_2, e_4] = e_4 \right.$\\
%%%%%%%%%%%%%%%%%%%%%%%%%%%%%%%%%%%%%%%%%%%%%%%%%%%%%%%%%%%%%%%%%%%%%%%%%
$(\mathfrak{r}_{2}',J_1(s,t))$, $t\neq0$: & $ \mu_{3}(a,b) :=\left\{ \begin{array}{l}[e_1, e_2] = e_2, [e_1, e_4] = e_4, \\{[e_2, e_3]} = a e_2+b e_4, [e_3, e_4] = b e_2-a e_4  \end{array}\right.$\\
                                          & \mbox{with } $ a=-{\frac {s}{t}},  b=-{\frac {{s}^{2}+{t}^{2}}{t}} $\\
$(\mathfrak{r}_{2}',J_2)$: & $\mu_{4}:=\left\{\begin{array}{l} [e_1, e_3] = e_3, [e_1, e_4] = e_4, \\ {[e_2, e_3]} = e_4, [e_2, e_4] = -e_3 \end{array}\right.$\\
$(\mathfrak{r}_{2}',J_3)$: & $\mu_{5}:=\left\{\begin{array}{l} [e_1, e_2] = e_2, [e_1, e_3] = e_4, [e_1, e_4] = e_4, \\ {[e_2, e_3]} = -e_4, [e_3, e_4] = -e_2 \end{array}\right.$\\
%%%%%%%%%%%%%%%%%%%%%%%%%%%%%%%%%%%%%%%%%%%%%%%%%%%%%%%%%%%%%%%%%%%%%%%%%
$(\mathfrak{rr}_{3,0} \times \mathbb{R},J)$: & $\mu_{6}(0,0)$ {\tiny See Table 2}\\
%%%%%%%%%%%%%%%%%%%%%%%%%%%%%%%%%%%%%%%%%%%%%%%%%%%%%%%%%%%%%%%%%%%%%%%%%
$(\mathfrak{rr}_{3,1} \times \mathbb{R},J)$: & $\mu_{3}(0,0)$ \\
%:=\left\{  [e_1, e2] = e2, [e1, e_4] = e_4 \right.$
%%%%%%%%%%%%%%%%%%%%%%%%%%%%%%%%%%%%%%%%%%%%%%%%%%%%%%%%%%%%%%%%%%%%%%%%%
$(\mathfrak{rr}'_{3,a} \times \mathbb{R},J_1)$, $a>0$: & $\mu_{7}(a):=\left\{  [e_1, e_2] = a e_2+e_4, [e_1, e_4] = -e_2+a e_4 \right.$\\
$(\mathfrak{rr}'_{3,a} \times \mathbb{R},J_2)$, $a>0$: & $\mu_{7}(-a)$\\
%superchequeado
%:=\left\{ [e_1, e2] = -a e2+e_4, [e1, e_4] = -e2-a e_4 \right.$
$(\mathfrak{rr}'_{3,0},J) \times \mathbb{R}$: & $\mu_{7}(0)$\\
%:=\left\{  [e1, e2] = e_4, [e1, e_4] = -e2 \right.$
%%%%%%%%%%%%%%%%%%%%%%%%%%%%%%%%%%%%%%%%%%%%%%%%%%%%%%%%%%%%%%%%%%%%%%%%%
$(\mathfrak{ so}(3,\mathbb{R}) \times \mathbb{R}, J(a))$ :  & $\mu_{8,+}(a):=\left\{ [e_1, e_2] = e_3, [e_1, e_3] = -e_2+ae_4, [e_2, e_3] = e_1 \right.$\\
$(\mathfrak{ sl}(2,\mathbb{R}) \times \mathbb{R}, J(a))$ :  & $\mu_{8,-}(a):=\left\{ [e_1, e_2] = e_3, [e_1, e_3] = +e_2+a e_4, [e_2, e_3] = e_1  \right.$\\
\end{tabular}\newline
\caption{Complex structures on four-dimensional real Lie algebras.}
\end{table}

\begin{table}
\centering
\begin{tabular}{ll}
$(\mathfrak{r}_{4,1},J)$: & $\mu_{9}:=\left\{ [e_1, e_2] = e_2, [e_1, e_3] = e_3+e_4, [e_1, e_4] = e_4 \right.$\\
%%%%%%%%%%%%%%%%%%%%%%%%%%%%%%%%%%%%%%%%%%%%%%%%%%%%%%%%%%%%%%%%%%%%%%%%%
$(\mathfrak{r}'_{4,\gamma,\delta},J_1)$, $\delta>0$: & $\mu_{6}(a,b):=\left\{ [e_1, e_2] = a e_2+b e_4, [e_1, e_3] = e_3, [e_1, e_4] = -b e_2+a e_4 \right.$\\
                                      & with $a=\gamma$, $b=-\delta$\\
$(\mathfrak{r}'_{4,\gamma,\delta},J_2)$, $\delta>0$: & $\mu_{6}(a,b)$ with $a=\gamma$, $b=\delta$\\
%%%%%%%%%%%%%%%%%%%%%%%%%%%%%%%%%%%%%%%%%%%%%%%%%%%%%%%%%%%%%%%%%%%%%%%%%
$(\mathfrak{r}_{4,\alpha,\alpha},J)$, $-1\leq \alpha < 1$, $\alpha \neq0$ : & $\mu_{6}(a,0):=\left\{ [e_1, e_2] = ae_2, [e_1, e_3] = e_3, [e_1, e_4] = ae_4 \right.$\\
                                                                            & with $a=\alpha$\\
$(\mathfrak{r}_{4,\alpha,1},J)$, $-1<\alpha \leq1$, $\alpha \neq0$: & $\mu_{6}(a,0)$ with $a = \sfrac{1}{\alpha}$\\
%:=\left\{ [e1, e2] = \sfrac{1}{a} e2, [e1, e_3] = e_3, [e1, e_4] = \sfrac{1}{a}e_4
%%%%%%%%%%%%%%%%%%%%%%%%%%%%%%%%%%%%%%%%%%%%%%%%%%%%%%%%%%%%%%%%%%%%%%%%%
$(\mathfrak{d}_{4,\lambda},J_1)$, $\lambda> \sfrac{1}{2}$, $a\neq1$ &  $\mu_{10}(a) :=\left\{ \begin{array}{l} [{\it e_1},{\it e_2}]=(a-1){\it e_2},[{\it e_1},{\it e_3}]={\it e_3},
                                                                        [{\it e_1},{\it e_4}]=a{\it e_4},\\ {[{\it e_2},{\it e_3}]}={\it e_4} \end{array}\right.$\\
                                                                    & with $a=\sfrac{1}{\lambda}$ \\
$(\mathfrak{d}_{4,\lambda},J_2)$, $\lambda> \sfrac{1}{2}$, $a\neq1$ & $\mu_{10}(a)$ with $a=\sfrac{1}{(1-\lambda)}$\\
%%%%%%%%%%%%%%%%%%%%%%%%%%%%%%%%%%%%%%%%%%%%%%%%%%%%%%%%%%%%%%%%%%%%%%%%%
$(\mathfrak{d}_{4,1},J)$, : & $\mu_{10}(1)$\\
%%%%%%%%%%%%%%%%%%%%%%%%%%%%%%%%%%%%%%%%%%%%%%%%%%%%%%%%%%%%%%%%%%%%%%%%%
$(\mathfrak{d}_{4,\sfrac{1}{2}},J_1)$, : & $\mu_{10}(2)$\\
$(\mathfrak{d}_{4,\sfrac{1}{2}},J_2)$, : & $\mu_{11}(+):=\left\{ \begin{array}{l} [e_1, e_2] = \tfrac{1}{2}e_2, [e_1, e_3] = e_3, [e_1, e_4] = \tfrac{1}{2}e_4, \\ {[e_2, e_4]} = e_3 \end{array}\right.$\\
$(\mathfrak{d}_{4,\sfrac{1}{2}},J_3)$, : & $\mu_{11}(-):=\left\{ \begin{array}{l} [e_1, e_2] = -\tfrac{1}{2}e_2, [e_1, e_3] = -e_3, [e_1, e_4] = -\tfrac{1}{2}e_4, \\ {[e_2, e_4]} = e_3 \end{array} \right.$\\
%%%%%%%%%%%%%%%%%%%%%%%%%%%%%%%%%%%%%%%%%%%%%%%%%%%%%%%%%%%%%%%%%%%%%%%%%
$(\mathfrak{d}_{4},J_1)$ : & $\mu_{10}(0)$\\
$(\mathfrak{d}_{4},J_2)$ : & $\mu_{12}:=\left\{ [e_1, e_2] = -e_2, [e_1, e_3] = e_3+e_4, [e_2, e_3] = e_4 \right.$\\
%:=\left\{ [e1, e_2] = -e2, [e1, e_3] = e_3, [e2, e_3] = e_4 \right.$
%%%%%%%%%%%%%%%%%%%%%%%%%%%%%%%%%%%%%%%%%%%%%%%%%%%%%%%%%%%%%%%%%%%%%%%%%
$(\mathfrak{d}'_{4,\delta},J_1)$, $\delta>0$ : & $\mu_{13,-}(a)$ with $a=-\sfrac{1}{2\delta}$\\
$(\mathfrak{d}'_{4,\delta},J_2)$, $\delta>0$ : & $\mu_{13,+}(a)$ with $a=-\sfrac{1}{2\delta}$\\
$(\mathfrak{d}'_{4,\delta},J_3)$, $\delta>0$ : & $\mu_{13,-}(a):=\left\{\begin{array}{l} [e_1, e_2] = ae_2-e_4, [e_1, e_3] = 2ae_3, \\ {[e_1, e_4]} = e_2+ae_4, [e_2, e_4] = e_3  \end{array} \right.$\\
                                               &  with $a=\sfrac{1}{2\delta}$\\
$(\mathfrak{d}'_{4,\delta},J_4)$, $\delta>0$ : & $\mu_{13,+}(a):=\left\{\begin{array}{l} [e_1, e_2] = ae_2+e_4, [e_1, e_3] = 2ae_3, \\ {[e_1, e_4]} = -e_2+ae_4,  {[e_2, e_4]} = e_3 \end{array} \right.$\\
                                               &  with $a=\sfrac{1}{2\delta}$\\
%%%%%%%%%%%%%%%%%%%%%%%%%%%%%%%%%%%%%%%%%%%%%%%%%%%%%%%%%%%%%%%%%%%%%%%%%

$(\mathfrak{d}'_{4,0},J_1)$, : & $\mu_{13,+}(0)$\\
$(\mathfrak{d}'_{4,0},J_2)$, : & $\mu_{13,-}(0)$\\
%%%%%%%%%%%%%%%%%%%%%%%%%%%%%%%%%%%%%%%%%%%%%%%%%%%%%%%%%%%%%%%%%%%%%%%%%
$(\mathfrak{h}_{4},J)$, : & $\mu_{14}:=\left\{\begin{array}{l} [e_1, e_2] = e_2, [e_1, e_3] = e_2+e_3,{[e_1, e_4]} = 2e_4, \\ {[e_2, e_3]} = e_4 \end{array} \right.$\\
\end{tabular}\newline
\caption{Complex structures on four-dimensional real Lie algebras }
\end{table}

\end{proposition}

\begin{remark}
It is important to note that $(\mathfrak{rr}_{3,1} \times \mathbb{R}, J)$ is holomorphically isomorphic to $(\mathbb{R}^{4}, \mu_{3}(a, 0), J_{\tiny{\mbox{cn}}})$ for any $a \in \mathbb{R}$.
\end{remark}

\section{Degenerations of 4-dimensional Lie algebras endowed with complex structures.}

In this section, we describe the procedure for determining the \textit{Hasse diagram} of degenerations associated with four-dimensional Lie algebras endowed with complex structures. Our approach begins with Corollary \ref{dercomp} and Proposition \ref{orden}, which provide a systematic way to organize pairs \((\mathbb{R}^{4}, \mu, J_{\tiny{\mbox{cn}}})\) based on the dimensions of their derivation algebra and holomorphic derivation algebra (see Table 3). Specifically, if a Lie algebra endowed with a complex structure, \((\mathbb{R}^{4}, \mu, J_{\tiny{\mbox{cn}}})\), properly degenerates to another, \((\mathbb{R}^{4}, \lambda, J_{\tiny{\mbox{cn}}})\), under the action of \(\operatorname{GL}(\mathbb{R}^{2n}, J_{\tiny{\mbox{cn}}})\), then the underlying Lie algebra \((\mathbb{R}^{4}, \mu)\) also degenerates to \((\mathbb{R}^{4}, \lambda)\) in the classical sense. As a consequence, we obtain the inequalities
\[
\operatorname{Dim}(\operatorname{Der}_{J_{\tiny{\mbox{cn}}}}(\mathbb{R}^{2n}, \mu)) < \operatorname{Dim}(\operatorname{Der}_{J_{\tiny{\mbox{cn}}}}(\mathbb{R}^{2n}, \lambda))
\]
and
\[
\operatorname{Dim}(\operatorname{Der}(\mathbb{R}^{2n}, \mu)) \leq \operatorname{Dim}(\operatorname{Der}(\mathbb{R}^{2n}, \lambda)).
\]
Once the ordering has been established, we justify the non-degenerations using the previously introduced invariants. Additionally, we make use of the classification of degenerations of four-dimensional real Lie algebras provided in \cite{nesterenko}. For completeness, we revisit this classification in the appendix, focusing on the non-degenerations of those Lie algebras that admit a complex structure, which are the ones relevant to our study.

\begin{landscape}

\begin{table}[p]
\centering
\begin{center}
    \begin{tabular}{c : c: p{4cm} : p{3.0cm} : p{5.5cm} : p{4.15cm} }

      \multirow{2}{*}{ \begin{tabular}{l}$\operatorname{Dim}$ \\ $\operatorname{Der}_{J_{\tiny{\mbox{cn}}}}$ \end{tabular}} &
      \multirow{2}{*}{ \begin{tabular}{l}$\operatorname{Dim}$ \\ $\operatorname{Der}$ \end{tabular} } &
     \multicolumn{4}{c}{$4$-dimensional Lie Algebras endowed with complex strucuture}\\

     \multicolumn{2}{c:}{}&
     \multicolumn{2}{c:}{}&
     Unimodular &
     Abelian complex str.\\
    \hdashline
     \multirow{1}{*}{$0$} & 4   &   & &  & $(\mathfrak{r}_2\times \mathfrak{r}_{2},J)$,  $(\mathfrak{r}^{'}_{2},J_2)$  \\
    \hdashline
     \multirow{3}{*}{$1$} & 4  &  &  &$(\mathfrak{sl}(2,\mathbb{R})\times \mathbb{R}, J(t))$ , $(\mathfrak{so}(3,\mathbb{R})\times \mathbb{R}, J(t))$ &  $(\mathfrak{d}_{4,1},J)$ \\
                          & 5  &  $(\mathfrak{ d}^{'}_{4,\delta},J_{k_{4}})$, $(\mathfrak{ d}_{4,\lambda}, J_{k_{2}})$, $(\mathfrak{ h}_{ 4},J)$
                               &
                               & $(\mathfrak{ d}^{'}_{4,0},J_{k_{2}})$, $(\mathfrak{d}_4,J_2)$ &  \\
                          & 7  & $(\mathfrak{ d}_{4,\frac{1}{2}}, J_2)$,$(\mathfrak{ d}_{4,\frac{1}{2}}, J_3)$ &  & &  \\
   \hdashline
     \multirow{5}{*}{$2$} & 4  &
                               &$(\mathfrak{r}^{'}_{2 }, J_{1}(a,b))$, $(\mathfrak{r}^{'}_{2 }, J_{3})$
                               &
                               &$(\mathfrak{r}^{'}_{2 }, J_{1}(0,1))$ \\
                          & 5  &
                               &
                               &$(\mathfrak{d}_4,J_1)$
                               & \\
                          & 6  & $(\mathfrak{r}^{'}_{4,c,d},J_{k_{2}})$
                               & $(\mathfrak{rr}^{'}_{3,a}\times\mathbb{R},J_{k_{2}})$,
                               & $(\mathfrak{r}^{'}_{4,-\frac{1}{2},d},J_{k_{2}})$, $(\mathfrak{rr}^{'}_{3,0}\times\mathbb{R},J)$
                               &  \\
                          & 7  &  $(\mathfrak{ d}_{4,\frac{1}{2}}, J_1)$ &  & & \\
                          & 8  &  $(\mathfrak{ r}_{4,\alpha,\alpha },J)$, $(\mathfrak{r}_{4,\alpha,1},J)$, $(\mathfrak{ r}_{4,1},J)$
                               &  $(\mathfrak{rr}_{3,1 }\times \mathbb{R},J)$
                               &  $(\mathfrak{ r}_{4,-\frac{1}{2},-\frac{1}{2}},J)$
                               &  $(\mathfrak{r}_2\times \mathbb{R}^{2},J)$ \\
   \hdashline
     \multirow{3}{*}{$4$} & 4  &    & $(\mathfrak{r}^{'}_{2 }, J_{1}(0,-1))$  &  & \\
                          & 10  &  \multicolumn{4}{c}{$(\mathfrak{h}_{3}\times \mathbb{R},J)$}  \\
                          & 12  &  $(\mathfrak{ r}_{4,1,1},J)$ &  & & \\
      \hdashline
   $8$ & $16$ & \multicolumn{4}{c}{$(\mathfrak{a}_{4},J)$}\\
    \end{tabular}
\end{center}
\caption{Dimension of  complex derivations. Here $k_{j}=1,2,\ldots,j$, $(a,b) \neq (0,\pm1)$ and $\lambda \neq 1, \frac{1}{2}$}
\end{table}
\end{landscape}

In the remainder of this section, we provide illustrative examples to show how to obtain the classification of orbit closures in the variety of 4-dimensional Lie algebras endowed with complex structures

\begin{example}
The following proposition illustrates how to exploit the parametrization of a $B$-orbit to degenerate a Lie algebra with a complex structure into another.
\end{example}

\begin{proposition}
The Lie algebra endowed with the complex structure  $(\mathfrak{ d}^{'}_{4,\delta},J_{1})$ degenerates to $(\mathfrak{r}^{'}_{4,\frac{1}{2}, \delta},J_{2})$ under the action of
$\operatorname{GL}(\mathbb{R}^{4}, J_{\operatorname{cn}})$.
\end{proposition}

\begin{proof}
  By the observation after Proposition \ref{Borel}, we need to find a Lie algebra law $\eta$ such that $\eta \in \overline{\operatorname{B}\cdot \mu_{13,-}\left( -\frac{1}{2\delta} \right) } \cap \operatorname{GL}(\mathbb{R}^{4}, J_{\operatorname{cn}}) \cdot \mu_{6}(\frac{1}{2},\delta)$.

It is so easy to write explicitly the elements of the $\operatorname{B}$-orbit of $\mu_{13,-}\left( -\frac{1}{2\delta} \right)$. If $g = \left(
          \begin{array}{cc}
            A & -B \\
            B & A \\
          \end{array}
        \right)
$ with $A= \left(
             \begin{array}{cc}
               t_1 & 0 \\
               x   & t_2 \\
             \end{array}
           \right)
$ and $B= \left(
            \begin{array}{cc}
              0 & 0 \\
              y & 0 \\
            \end{array}
          \right)
$, with $t_1,t_2,x,y \in \mathbb{R}$ and $t_1,t_2 >0$, then $ g \cdot \mu_{13,-}\left( -\frac{1}{2\delta} \right)$ is of the form:

$$
\begin{array}{l}
[{\it e_1},{\it e_2}]=-{\frac { \left( {t_{{2}}}^{2}+2\,{y}^{2}\delta \right) }{2{t_{{2}}}^{2}t_{{1}}\delta}}{\it e_2}
+{\frac {y }{{t_{{2}}}^{2}}}{\it e_3}
+{\frac { \left( -{t_{{2}}}^{2}+xy \right) }{{t_{{2}}}^{2}t_{{1}}}}{\it e_4},\\
{[{\it e_1},{\it e_3}]}=-{\frac { \left( 2\,{t_{{2}}}^{2}x\delta-{t_{{2}}}^{2}y+2\,{y}^{3}\delta+2\,y{x}^{2}\delta \right) }{2{t_{{2}}}^{2}{t_{{1}}}^{2}\delta}}{\it e_2}
+{\frac { \left({y}^{2}\delta+{x}^{2}\delta-{t_{{2}}}^{2} \right) }{{t_{{2}}}^{2}t_{{1}}\delta}}{\it e_3}
+{\frac { \left( -{t_{{2}}}^{2}x-2\,{t_{{2}}}^{2}y\delta+2\,x{y}^{2}\delta+2\,{x}^{3}\delta \right) }{2{t_{{2}}}^{2}{t_{{1}}}^{2}\delta}}{\it e_4},\\
{[{\it e_1},{\it e_4}]}={\frac { \left( {t_{{2}}}^{2}+xy \right) }{{t_{{2}}}^{2}t_{{1}}}}{\it e_2}
-{\frac {x}{{t_{{2}}}^{2}}}{\it e_3}
-{\frac { \left( {t_{{2}}}^{2}+2\,{x}^{2}\delta \right) }{2{t_{{2}}}^{2}t_{{1}}\delta}}{\it e_4},\\
{[{\it e_2},{\it e_3}]}={\frac {xy }{{t_{{2}}}^{2}t_{{1}}}}{\it e_2}
-{\frac {x}{{t_{{2}}}^{2}}}{\it e_3}
-{\frac {{x}^{2}}{{t_{{2}}}^{2}t_{{1}}}}{\it e_4},
[{\it e_2},{\it e_4}]=-{\frac {y}{{t_{{2}}}^{2}}}{\it e_2}
+{\frac {t_{{1}}}{{t_{{2}}}^{2}}}{\it e_3}
+{\frac {x}{{t_{{2}}}^{2}}}{\it e_4},
[{\it e_3},{\it e_4}]=-{\frac {{y}^{2}}{{t_{{2}}}^{2}t_{{1}}}}{\it e_2}
+{\frac {y}{{t_{{2}}}^{2}}}{\it e_3}
+{\frac {xy}{{t_{{2}}}^{2}t_{{1}}}}{\it e_4}
\end{array}
$$
It follows that  $B \cdot \mu_{13,-}\left( -\frac{1}{2\delta} \right)$ is a subset of the subspace $W$ of $C^{2}_{ J_{\tiny{\mbox{cn}}} }(\mathbb{R}^{4};\mathbb{R}^{4})$ given by
$$
W =
\left\{
\begin{array}{l}
[{\it e_1},{\it e_2}]=c_{{1}}{\it e_2}+c_{{2}}{\it e_3}+c_{{3}}{\it e_4},
[{\it e_1},{\it e_3}]=c_{{4}}{\it e_2}+c_{{5}}{\it e_3}+c_{{6}}{\it e_4},\\
{[{\it e_1},{\it e_4}]}=c_{{7}}{\it e_2}+c_{{8}}{\it e_3}+c_{{9}}{\it e_4},
[{\it e_2},{\it e_3}]=c_{{10}}{\it e_2}+c_{{8}}{\it e_3}+c_{{12}}{\it e_4},\\
{[{\it e_2},{\it e_4}]}=c_{{13}}{\it e_2}+c_{{14}}{\it e_3}+c_{{15}}{\it e_4},\\
{[{\it e_3},{\it e_4}]}= \left(c_{{1}} -c_{{9}}+c_{{12}} \right) {\it e_2}+c_{{2}}{\it e_3}+ \left( c_{{3}}+c_{{7}}-c_{{10}} \right) {\it e_4}
\end{array} : c_i \in \mathbb{R}
\right\}
$$
 It is straightforward to see that the Lie algebra laws in $W$ that are in the $\operatorname{GL}(\mathbb{R}^{4},J_{\tiny{\mbox{cn}}} )$-orbit of $\mu_{6}(\frac{1}{2},\delta)$ satisfy the condition $[{\it e_2},{\it e_3}]=[{\it e_2},{\it e_4}]=[{\it e_3},{\it e_4}]=0$. With the purpose of getting closer to the $\operatorname{GL}(\mathbb{R}^{4},J_{\tiny{\mbox{cn}}} )$-orbit of $\mu_{6}(\frac{1}{2},\delta)$, we can attempt to restrict our attention to the Lie algebra laws in the $\operatorname{B}$-orbit of $\mu_{13,-}\left( -\frac{1}{2\delta} \right)$ that satisfy $[{\it e_2},{\it e_3}]$ and $[{\it e_3},{\it e_4}]$ are equal to zero, or equivalently, $x=0$ and $y=0$:

$$
\begin{array}{l}
[{\it e_1},{\it e_2}]=-{\frac {1}{2t_{{1}}\delta}}{\it e_2}
-{\frac {1}{t_{{1}}}}{\it e_4},
[{\it e_1},{\it e_3}]=-{\frac {1}{t_{{1}}\delta}}{\it e_3},
[{\it e_1},{\it e_4}]={\frac {1}{t_{{1}}}}{\it e_2}
-{\frac {1}{2t_{{1}}\delta}}{\it e_4},
[{\it e_2},{\it e_4}]={\frac {t_{{1}}}{{t_{{2}}}^{2}}}{\it e_3}.
\end{array}
$$
This motivates to let $g(u) = \operatorname{diag}(-\frac{1}{\delta} , \exp(u), -\frac{1}{\delta} , \exp(u) )$ with $u \in \mathbb{R}$, and so
$\xi_u = g(u) \cdot \mu_{13,-}\left( -\frac{1}{2\delta} \right)$ is the Lie algebra law:
$$
[{\it e_1},{\it e_2}]=\frac{1}{2}\,{\it e_2}+\delta\,{\it e_4},
[{\it e_1},{\it e_3}]={\it e_3},
[{\it e_1},{\it e_4}]=-\delta\,{\it e_2} + \frac{1}{2}\,{\it e_4},
[{\it e_2},{\it e_4}]=-{\frac { 1}{\delta\, {{\rm e}^{2u}}}}{\it e_3},
$$
It is clear that $g(u) \in \operatorname{GL}(\mathbb{R}^{4}, J_{\tiny{\mbox{cn}}} )$ and $g(u) \cdot \mu_{13,-}\left( -\frac{1}{2\delta} \right)$  tends to $\mu_{6}(\frac{1}{2},\delta)$ as $u$ tends to $\infty$; which is much stronger than what we wanted to prove.
\end{proof}

\begin{example}
Here we describe how to use the $\operatorname{GL}(\mathbb{R}^{2n}, J_{\tiny{\mbox{cn}}})$-invariants defined above to prove a non-degeneracy result.
\end{example}

\begin{proposition}
If $k \neq 2$ or $(c, d) \neq \left( \frac{1}{2}, \delta \right)$, then the Lie algebra endowed with the complex structure $(\mathfrak{d}'_{4, \delta}, J_1)$ does not degenerate to $(\mathfrak{r}'_{4, c, d}, J_k)$ under the action of $\operatorname{GL}(\mathbb{R}^4, J_{\operatorname{cn}})$.
\end{proposition}

\begin{proof}
Suppose, for the sake of contradiction, that $\mu_6(c, (-1)^k d) \in \overline{\operatorname{GL}(\mathbb{R}^4, J_{\tiny{\mbox{cn}}}) \cdot \mu_{13, -}\left( -\frac{1}{2\delta} \right)}$.
First, consider the $\operatorname{GL}(\mathbb{R}^4, J_{\tiny{\mbox{cn}}})$-equivariant continuous function $\psi_{\alpha, \beta}$, where $\psi_{\alpha, \beta}$ is the function defined before Equation (\ref{funcioninv}).
Therefore, the skew-symmetric algebra $\mathfrak{A} = (\mathbb{R}^4, \psi_{\alpha, \beta}(\mu_{13, -}\left( -\frac{1}{2\delta} \right)))$ degenerates to $\mathfrak{B} = (\mathbb{R}^4, \psi_{\alpha, \beta}(\mu_6(c, (-1)^k d)))$ in the usual sense of degenerations of  algebras, where

$$
\mathfrak{A} =
\left\{
\begin{array}{l}
[e_1, e_2] = \frac{-1 + 2\alpha\delta + 2\beta\delta}{2\delta} e_2 - \frac{2\delta + \alpha + \beta}{2\delta} e_4, {[e_1, e_3]} = \frac{\alpha}{\delta} e_1 - \frac{1}{\delta} e_3, \\
{[e_1, e_4]} = \frac{2\delta + \alpha + \beta}{2\delta} e_2 + \frac{-1 + 2\alpha\delta + 2\beta\delta}{2\delta} e_4, \\
{[e_2, e_3]} = -\frac{\beta}{2\delta} e_2 - \beta e_4, {[e_2, e_4]} = -\alpha e_1 + e_3, {[e_3, e_4]} = -\beta e_2 + \frac{\beta}{2\delta} e_4.
\end{array}
\right.
$$

and

$$
\mathfrak{B} =
\left\{
\begin{array}{l}
[e_1, e_2] = \left(c - (-1)^k \alpha  d - (-1)^k \beta d\right) e_2 + \left((-1)^k d + \alpha c + \beta c\right) e_4, \\
{[e_1, e_3]} = -\alpha e_1 + e_3, \\
{[e_1, e_4]} = \left((-1)^{k+1} d - \beta c - \alpha c\right) e_2 + \left(c - (-1)^k \alpha d -  (-1)^k \beta d\right) e_4, \\
{[e_2, e_3]} = \beta c e_2 +  (-1)^k \beta d e_4, {[e_3, e_4]} =  (-1)^k \beta d e_2 - \beta c e_4.
\end{array}
\right.
$$

Taking $\alpha = \beta = \frac{1}{2\delta}$, we obtain that $\mathfrak{A}$ is unimodular, and thus $\mathfrak{B}$ is also unimodular. Therefore, we have the system of equations:

$$
\left\{
\begin{array}{l}
\frac{\delta + 2c\delta + 2(-1)^{k+1} d}{\delta} = 0, \\
\frac{-1 + 2c}{2\delta} = 0.
\end{array}
\right.
$$

From this, we conclude that $c = \frac{1}{2}$ and $d = (-1)^k \delta$. Recall that by the definition of $(\mathfrak{r}'_{4, c, d}, J_k)$, the real number $d$ is positive, so $k = 2$ and $(c, d) = \left( \frac{1}{2}, \delta \right)$, which contradicts our hypothesis.
\end{proof}

\begin{remark}
It is proved in \cite{nesterenko} that the Lie algebra $\mathfrak{d}'_{4, \delta}$ degenerates to $\mathfrak{r}'_{4, c, d}$ if and only if $c = \frac{1}{2}$ and $d = \frac{1}{\delta}$.
Therefore, in the proof above, we could have focused solely on proving that $k = 2$.
\end{remark}

As an additional illustration of the techniques employed, we now present another example.

\begin{proposition}
The Lie algebra endowed with the complex structure $(\mathfrak{h}_4, J)$ does not degenerate to $(\mathfrak{r}_{4, \frac{1}{2}, \frac{1}{2}}, J)$ under the action of $\operatorname{GL}(\mathbb{R}^4, J_{\operatorname{cn}})$.
\end{proposition}

\begin{proof}
Suppose, for the purpose of deriving a contradiction, that $\mu_6\left(\frac{1}{2}, 0\right) \in \overline{\operatorname{GL}(\mathbb{R}^4, J_{\tiny{\mbox{cn}}}) \cdot \mu_{14}}$. Evaluating the $\operatorname{GL}(\mathbb{R}^4, J_{\tiny{\mbox{cn}}})$-equivariant continuous function $\psi_{1, 0}$ at both $\mu_{14}$ and $\mu_6\left(\frac{1}{2}, 0\right)$, it follows that the skew-symmetric algebra $\mathfrak{A} = \left( \mathbb{R}^4, \psi_{1, 0}(\mu_{14}) \right)$ degenerates to $\mathfrak{B} = \left( \mathbb{R}^4, \psi_{1, 0}(\mu_6\left(\frac{1}{2}, 0\right)) \right)$, where

$$
\mathfrak{A} =
\left\{
\begin{array}{l}
[e_1, e_2] = e_2 + e_4, {[e_1, e_3]} = -e_1 + e_2 + e_3 + e_4, {[e_1, e_4]} = -2e_2 + 2e_4, {[e_2, e_3]} = -e_2 + e_4.
\end{array}
\right.
$$

and

$$
\mathfrak{B} =
\left\{
\begin{array}{l}
{[e_1, e_2]} = \frac{1}{2} e_2 + \frac{1}{2} e_4, {[e_1, e_3]} = -e_1 + e_3, {[e_1, e_4]} = -\frac{1}{2} e_2 + \frac{1}{2} e_4.
\end{array}
\right.
$$

Now, consider the modified Killing form $\kappa_g : C^2(\mathbb{R}^4; \mathbb{R}^4) \rightarrow \Sigma^2(\mathbb{R}^4; \mathbb{R})$. We have that $\kappa_g(\mathfrak{A}) \xrightarrow{\text{deg}} \kappa_g(\mathfrak{B})$, where $\kappa_g(\mathfrak{A})$ is the symmetric tensor:

$$
(2 + 16g) e_1^\ast \odot e_1^\ast + 2(4 + 8g) e_1^\ast \odot e_3^\ast + (2 + 4g) e_3^\ast \odot e_3^\ast,
$$
and $\kappa_g(\mathfrak{B})$ is:

$$
(1 + 4g) e_1^\ast \odot e_1^\ast + 2(1 + 2g) e_1^\ast \odot e_3^\ast + (1 + g) e_3^\ast \odot e_3^\ast.
$$

Taking $g = -\frac{1}{2}$, we obtain a contradiction, since the dimension of the radical of $\kappa_{-\frac{1}{2}}(\mathfrak{A}) = -6 \, e_1^\ast \odot e_1^\ast$ is 3, while the dimension of the radical of $\kappa_{-\frac{1}{2}}(\mathfrak{B}) = - e_1^\ast \odot e_1^\ast + \frac{1}{2} e_3^\ast \odot e_3^\ast$ is 2.
\end{proof}

\begin{remark}
It is important to emphasize that the Lie algebra $ \mathfrak{h}_4 $ degenerates to the Lie algebra $\mathfrak{r}_{4, \frac{1}{2}, \frac{1}{2}}$ with respect to the action of the group $\operatorname{GL}(\mathbb{R}^4)$.
\end{remark}

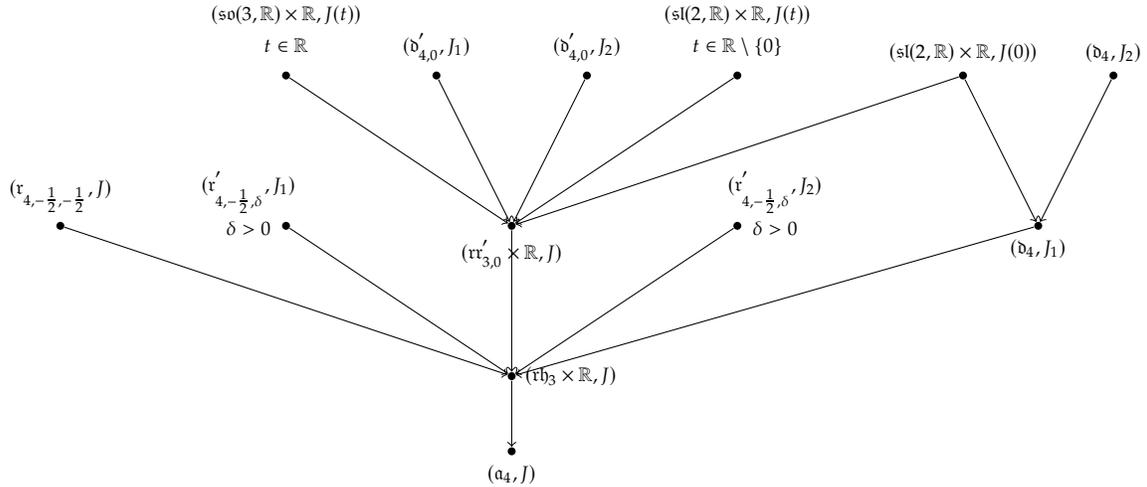
\begin{figure}[h!]
\begin{center}
\begin{tikzpicture}
\node[circle,fill,inner sep=0pt,minimum size=3pt,label=above:{\scriptsize$(\mathfrak{d}^{'}_{4,0},J_1)$}] (d4p1) at (-1,+0) {};
\node[circle,fill,inner sep=0pt,minimum size=3pt,label=above:{\scriptsize$(\mathfrak{d}^{'}_{4,0},J_2)$}] (d4p2) at (+1,+0) {};
%\node[circle,fill,inner sep=0pt,minimum size=3pt,label=above:{}] (so3a) at (-2,+0) {};
\node[circle,fill,inner sep=0pt,minimum size=3pt,label=above:{\begin{tabular}{c} \scriptsize$(\mathfrak{so}(3,\mathbb{R})\times \mathbb{R}, J(t))$ \\ \scriptsize$t \in \mathbb{R}$ \end{tabular}}] (so3) at (-3,+0) {};
%\node[circle,fill,inner sep=0pt,minimum size=3pt,label=left:{}]  (so3b) at (-4,+0) {};
%\node[circle,fill,inner sep=0pt,minimum size=3pt,label=right:{}] (sl2a) at (+2,+0) {};
\node[circle,fill,inner sep=0pt,minimum size=3pt,label=above:{\begin{tabular}{c} \scriptsize$(\mathfrak{sl}(2,\mathbb{R})\times \mathbb{R}, J(t))$ \\ \scriptsize$t \in \mathbb{R}\setminus\{0\}$ \end{tabular}} ] (sl2)  at (+3,+0) {};
%\node[circle,fill,inner sep=0pt,minimum size=3pt,label=left:{ }] (sl2b) at (+4,+0) {};
\node[circle,fill,inner sep=0pt,minimum size=3pt,label=above:{\scriptsize$(\mathfrak{sl}(2,\mathbb{R})\times \mathbb{R}, J(0))$}](sl20)  at (+6,+0) {};
\node[circle,fill,inner sep=0pt,minimum size=3pt,label=above:{\scriptsize$(\mathfrak{d}_4,J_2)$ }](d42)  at (+8,+0) {};
\node[circle,fill,inner sep=0pt,minimum size=3pt,label=below:{\scriptsize$(\mathfrak{rr}^{'}_{3,0}\times \mathbb{R},J)$}] (r3p)  at (+0,-2) {};
%\node[circle,fill,inner sep=0pt,minimum size=3pt,label=right:{}] (r4p2a) at (-2,-2) {};
\node[circle,fill,inner sep=0pt,minimum size=3pt,label=above:{}] (r4p1) at (-3,-2) {};
\node[ label=above:{\begin{tabular}{c} \scriptsize$(\mathfrak{r}^{'}_{4,-\tfrac{1}{2},\delta},J_1)$ \\ \scriptsize$ \delta >0$ \end{tabular}}] (r4p1label) at (-3.5,-2.5) {};
%\node[circle,fill,inner sep=0pt,minimum size=3pt,label=above:{}] (r4p2b) at (-4,-2) {};
%\node[circle,fill,inner sep=0pt,minimum size=3pt,label=left:{}]  (r4p1a) at (+2,-2) {};
\node[circle,fill,inner sep=0pt,minimum size=3pt,label=right:{}]  (r4p2) at (+3,-2) {};
\node[ label=above:{\begin{tabular}{c} \scriptsize$(\mathfrak{r}^{'}_{4,-\tfrac{1}{2},\delta},J_2)$ \\ \scriptsize$ \delta >0$ \end{tabular}}] (r4p2label) at (3.5,-2.5) {};
%\node[circle,fill,inner sep=0pt,minimum size=3pt,label=right:{}] (r4p1b) at (+4,-2) {};
\node[circle,fill,inner sep=0pt,minimum size=3pt,label=above:{\scriptsize$(\mathfrak{r}_{4,-\tfrac{1}{2},-\tfrac{1}{2}},J)$ }] (r4) at (-6,-2) {};
\node[circle,fill,inner sep=0pt,minimum size=3pt,label=below:{\scriptsize$(\mathfrak{d}_4,J_1)$ }] (d41)at (+7,-2) {};
\node[circle,fill,inner sep=0pt,minimum size=3pt,label=right:{\scriptsize$(\mathfrak{rh}_{3}\times \mathbb{R},J)$ }] (r3) at (+0,-4) {};
\node[circle,fill,inner sep=0pt,minimum size=3pt,label=below:{\scriptsize$(\mathfrak{a}_4,J)$ }] (a4) at (+0,-5) {};
%\draw[thick] (so3a)edge (so3b)
%\draw[thick] (sl2a) edge (sl2) (sl2) edge (sl2b);
%\draw[thick] (r4p2a) edge (r4p2b) (r4p1a) edge (r4p1b);
\draw [->] (d4p1) edge (r3p)  (d4p2) edge (r3p) (so3) edge (r3p); %(so3a) edge (r3p) (so3b) edge (r3p)
%\draw [->] (sl2a) edge (r3p) (sl2b) edge (r3p);
\draw [->] (sl2)edge (r3p)  (sl20)edge (r3p) ;
%\draw [->] (r4p2a)edge (r3)(r4p2b)edge (r3) (r4p1a)edge (r3)(r4p1b)edge (r3);
\draw [->] (r4p2)edge (r3) (r4p1)edge (r3);
\draw [->] (d41) edge (r3);
\draw [->] (r3p) edge (r3) (r4) edge (r3);
\draw [->] (d42) edge (d41) (sl20)edge (d41);
\draw [->] (r3)  edge (a4);
\end{tikzpicture}
\caption{Degenerations of $4$-dimensional unimodular Lie algebras endowed with complex structure.}\label{unimodu}
\end{center}
\end{figure}

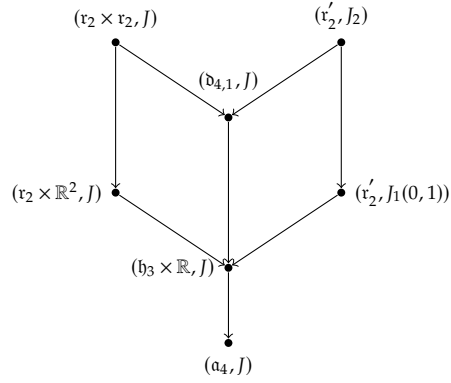
\begin{figure}[h!]
\begin{center}
\begin{tikzpicture}
\node[circle,fill,inner sep=0pt,minimum size=3pt,label=above:{\scriptsize$(\mathfrak{r}_2\times \mathfrak{r}_{2},J)$}] (r2r2) at (-1.5,0) {};
\node[circle,fill,inner sep=0pt,minimum size=3pt,label=above:{\scriptsize$(\mathfrak{r}^{'}_{2},J_2)$}] (r2p3) at (+1.5,0) {};
\node[ label=above:{\scriptsize$(\mathfrak{d}_{4,1},J)$}] (d41label) at (0,-0.9) {};
\node[circle,fill,inner sep=0pt,minimum size=3pt, label=above:{}] (d41) at (0,-1) {};
\node[circle,fill,inner sep=0pt,minimum size=3pt,label=left:{\scriptsize$(\mathfrak{r}_2\times \mathbb{R}^{2},J)$}] (r2z) at (-1.5,-2) {};
\node[circle,fill,inner sep=0pt,minimum size=3pt,label=right:{\scriptsize$(\mathfrak{r}^{'}_{2},J_1(0,1))$}] (r2p1) at (+1.5,-2) {};
\node[circle,fill,inner sep=0pt,minimum size=3pt,label=left:{\scriptsize$(\mathfrak{h}_{3}\times \mathbb{R},J)$ }] (h3) at (0,-3) {};
\node[circle,fill,inner sep=0pt,minimum size=3pt,label=below:{\scriptsize$(\mathfrak{a}_{4},J)$ }] (a4) at (0,-4) {};
\draw [->] (r2r2) edge (d41) (r2r2) edge (r2z) (r2p3) edge (d41) (r2p3) edge  (r2p1);
\draw [->] (d41) edge (h3) (r2z) edge (h3) (r2p1) edge (h3)   (h3) edge (a4);
\end{tikzpicture}
\caption{Degenerations of $4$-dimensional Lie algebras endowed with abelian complex structure.}\label{abelian}
\end{center}
\end{figure}

By using arguments similar to those used above we obtain the following theorem:

\begin{theorem}\label{principal}
  The Hasse diagram of degenerations in $\mathcal{L}_{J_{\tiny{\mbox{cn}}}}(\mathbb{R}^{4})$ is summarized in Figure \ref{unimodu} through Figure 5.

\begin{figure}[h!]
\begin{center}
\begin{tikzpicture}

\node[circle,fill,inner sep=0pt,minimum size=3pt,label=left:{$(\mathfrak{r}_{2}',J_3)$}] (a1) at (0,+0) {};
\node[circle,fill,inner sep=0pt,minimum size=3pt,label=above:{\begin{tabular}{c} \scriptsize$(\mathfrak{r}_{2}',J_1(a,b))$ \\ \scriptsize$ (a,b)\neq (0,\pm1) $ \end{tabular}}] (a2) at (2,+0) {};
\node[circle,fill,inner sep=0pt,minimum size=3pt,label=above:{$(\mathfrak{r}_{2}',J_2)$}] (a3) at (4,+0) {};
\node[circle,fill,inner sep=0pt,minimum size=3pt,label=above:{\begin{tabular}{c} \scriptsize$(\mathfrak{rr}'_{3,a}\times \mathbb{R},J_k)$ \\ \scriptsize$a>0,\, k=1,2$ \end{tabular}}] (a4) at (6,+0) {};
\node[circle,fill,inner sep=0pt,minimum size=3pt,label=right:{$(\mathfrak{rr}_{3,1}\times \mathbb{R},J)$}] (a5) at (8,+0) {};
\node[circle,fill,inner sep=0pt,minimum size=3pt,label=left:{$(\mathfrak{r}_{2}',J_1(0,-1))$}] (a6) at (0,-2) {};
\node[circle,fill,inner sep=0pt,minimum size=3pt,label=right:{$(\mathfrak{rh}_{3} \times \mathbb{R},J)$}] (a7) at (4,-2) {};
\node[circle,fill,inner sep=0pt,minimum size=3pt,label=below:{$(\mathfrak{a}_{4}, J)$}] (a8) at (2,-4) {};

\draw [->] (a1) edge (a7) (a2) edge (a7) (a3) edge (a7) (a4) edge (a7) (a5) edge (a7);
\draw [->] (a1) edge (a6);
\draw [->] (a4) edge (a8) (a6) edge (a8);

\end{tikzpicture}$ $\newline
{Figure 3.}
\end{center}
\end{figure}

\begin{figure}[h!]
\begin{center}
\begin{tikzpicture}

\node[circle,fill,inner sep=0pt,minimum size=3pt,label=above:{\begin{tabular}{c} \scriptsize$(\mathfrak{d}'_{4,\delta},J_1)$ \\ \scriptsize$ \delta>0 $ \end{tabular}}] (a1) at (0,+0) {};
\node[circle,fill,inner sep=0pt,minimum size=3pt,label=above:{\begin{tabular}{c} \scriptsize$(\mathfrak{d}'_{4,\delta},J_4)$ \\ \scriptsize$ \delta>0 $ \end{tabular}}] (a2) at (2,+0) {};
\node[circle,fill,inner sep=0pt,minimum size=3pt,label=above:{\begin{tabular}{c} \scriptsize$(\mathfrak{d}'_{4,\delta},J_2)$ \\ \scriptsize$ \delta>0 $ \end{tabular}}] (a3) at (4,+0) {};
\node[circle,fill,inner sep=0pt,minimum size=3pt,label=above:{\begin{tabular}{c} \scriptsize$(\mathfrak{d}'_{4,\delta},J_3)$ \\ \scriptsize$ \delta>0 $ \end{tabular}}] (a4) at (6,+0) {};
\node[circle,fill,inner sep=0pt,minimum size=3pt,label=above:{\begin{tabular}{c} \scriptsize$(\mathfrak{d}_{4,a},J_k)$ \\ \scriptsize$a>\frac{1}{2}, a\neq1$\\ \scriptsize$k=1,2 $ \end{tabular}}] (a5) at (8,+0) {};
%%%%
\node[circle,fill,inner sep=0pt,minimum size=3pt,label=left:{\begin{tabular}{c} \scriptsize$(\mathfrak{r}'_{4,c,d},J_2)$ \\ \scriptsize$ c\neq -\frac{1}{2}$ \\ \scriptsize$  d>0$ \end{tabular}}] (a6) at (1,-2) {};
\node[circle,fill,inner sep=0pt,minimum size=3pt,label=left:{\begin{tabular}{c} \scriptsize$(\mathfrak{r}'_{4,c,d},J_1)$ \\ \scriptsize$  c\neq -\frac{1}{2}$ \\ \scriptsize$  d>0 $ \end{tabular}}] (a7) at (5,-2) {};
\node[circle,fill,inner sep=0pt,minimum size=3pt,label=above:{\begin{tabular}{c} \scriptsize$(\mathfrak{r}_{4,\alpha,\alpha},J)$ \\\scriptsize$ \alpha \in \mathbb{R} \setminus\{0\}$ \\ \scriptsize$ \alpha \neq \pm \frac{1}{2},1 $ \end{tabular}}] (a8) at (9,-2) {};
\node[circle,fill,inner sep=0pt,minimum size=3pt,label=right:{$(\mathfrak{rh}_{3}\times \mathbb{R},J)$}] (a9) at (5,-4) {};
\node[circle,fill,inner sep=0pt,minimum size=3pt,label=below:{$(\mathfrak{a}_{4}, J)$}] (a10) at (5,-6) {};
%%%%
\node[label={[rotate=90]above:{\scriptsize $(c,d)=(\frac{1}{2},\delta)$}}] (a11) at (1.25,-1) {};
\node[label={[rotate=90]above:{\scriptsize $(c,d)=(\frac{1}{2},\delta)$}}] (a11) at (5.25,-1) {};

%\node[label={[label distance=0.5cm,text depth=-1ex,rotate=-90]right:a long text}] at (2.5,1) {};

\draw [->] (a1) edge (a6) (a2) edge (a6);
\draw [->] (a3) edge (a7) (a4) edge (a7);
\draw [->] (a5) edge (a9);
\draw [->] (a6) edge (a9) (a7) edge (a9) (a8) edge (a9);
\draw [->] (a9) edge (a10);

\end{tikzpicture}
$ $\newline
{Figure 4.}
\end{center}
\end{figure}

\begin{figure}[h!]
\begin{center}
\begin{tikzpicture}

\node[circle,fill,inner sep=0pt,minimum size=3pt,label=above:{$(\mathfrak{d}_{4,\sfrac{1}{2}},J_2)$}] (a1) at (0,+0) {};
\node[circle,fill,inner sep=0pt,minimum size=3pt,label=above:{$(\mathfrak{d}_{4,\sfrac{1}{2}},J_3)$}] (a2) at (2,+0) {};
\node[circle,fill,inner sep=0pt,minimum size=3pt,label=above:{$(\mathfrak{h}_{4},J)$}] (a3) at (4,+0) {};
%%%%%%
\node[circle,fill,inner sep=0pt,minimum size=3pt,label=left:{$(\mathfrak{r}_{4,\sfrac{1}{2},\sfrac{1}{2}},J)$}] (a4) at (1,-2) {};
\node[circle,fill,inner sep=0pt,minimum size=3pt,label=right:{$(\mathfrak{d}_{4,\sfrac{1}{2}},J_1)$}] (a5) at (3,-2) {};
%\node[label={[rotate=90]above:{$(\mathfrak{r}_{4,\sfrac{1}{2},\sfrac{1}{2}},J)$}}] (a4a) at (1.3,-1) {};
%\node[label={[rotate=90]above:{$(\mathfrak{d}_{4,\sfrac{1}{2}},J_1)$}}] (a5a) at (3.3,-1) {};
%%%%%%
\node[circle,fill,inner sep=0pt,minimum size=3pt,label=right:{$(\mathfrak{r}_{4,1},J)$}] (a6) at (5,-2) {};
\node[circle,fill,inner sep=0pt,minimum size=3pt,label=left:{$(\mathfrak{rh}_{3}\times \mathbb{R},J)$}] (a7) at (3,-4) {};
\node[circle,fill,inner sep=0pt,minimum size=3pt,label=right:{$(\mathfrak{r}_{4,1,1}, J )$}] (a8) at (5,-4) {};
\node[circle,fill,inner sep=0pt,minimum size=3pt,label=below:{$(\mathfrak{a}_4 , J )$}] (a9) at (3,-6) {};

\draw [->]  (a1) edge (a4) (a2) edge (a4) (a2) edge (a5) (a3) edge (a5) ;
\draw [->]  (a4) edge (a7) (a5) edge (a7) (a6) edge (a8);
\draw [->]  (a7) edge (a9) (a8) edge (a9);
\draw [->]  (a6) edge (a7);

\end{tikzpicture}$ $\newline
{Figure 5.}
\end{center}
\end{figure}

\end{theorem}

In the diagrams presented in Theorem \ref{principal}, we employ a slight abuse of notation: certain labeled vertices correspond to one or two families of Lie algebras that admit a complex structure. For example, the node labeled as ${\begin{tabular}{c} \scriptsize$(\mathfrak{rr}'_{3,a}\times \mathbb{R},J_k)$ \\ \scriptsize$a>0,\, k=1,2$ \end{tabular}}$ encapsulates the one-parameter families $(\mathfrak{rr}'_{3,a}\times \mathbb{R},J_1)$ and $(\mathfrak{rr}'_{3,a}\times \mathbb{R},J_2)$ for $a>0$. Similarly, the vertex labeled as ${\begin{tabular}{c} \scriptsize$(\mathfrak{r}_{2}',J_1(a,b))$ \\ \scriptsize$ (a,b)\neq (0,\pm1) $ \end{tabular}}$
represents the two-parameter family $(\mathfrak{r}_{2}',J_1(a,b))$ under the condition $(a,b)\neq (0,\pm1)$. The directed edge connecting this vertex to the node labeled $(\mathfrak{rh}_{3} \times \mathbb{R},J)$ signifies that every pair in the aforementioned family degenerates into $(\mathfrak{rh}_{3} \times \mathbb{R},J)$. In Figure 4, the edges labeled with the condition $(c,d)=(\frac{1}{2},\delta)$ indicates that  $(\mathfrak{d}'_{4,\delta},J_1)$ or $(\mathfrak{d}'_{4,\delta},J_)$ degenerates to the pair $(\mathfrak{r}'_{4,c,d},J_2)$ if and only if the condition is satisfied.

For the sake of clarity, only the essential degenerations required for Theorem \ref{principal} are explicitly depicted, while those following from transitivity have been omitted. A comprehensive listing of all fundamental degenerations can be found in the Appendix.

\section{Applications}

Let \((\mathbb{R}^{2n}, \mu, J_{\tiny{\mbox{cn}}})\) be a Lie algebra equipped with a complex structure. A positive-definite inner product \(\langle \cdot, \cdot \rangle\) on \(\mathbb{R}^{2n}\) is said to be \textit{compatible} with \(J_{\tiny{\mbox{cn}}}\) if \(J_{\tiny{\mbox{cn}}}\) acts as an isometry, which is equivalent to being a skew-symmetric operator with respect to \(\langle \cdot, \cdot \rangle\). In this case, the pair \((J_{\tiny{\mbox{cn}}}, \langle \cdot, \cdot \rangle)\) defines an \textit{Hermitian structure} on the Lie algebra \((\mathbb{R}^{2n}, \mu)\).

We adopt an approach similar to that developed by Jens Heber in \cite{Heber}, where the focus is placed on varying \textit{Lie algebra structures rather than scalar products} when analyzing left-invariant metrics on Lie groups. Given a Lie algebra structure \(\mu\) on \(\mathbb{R}^{n}\) and the canonical inner product \(\langle \cdot, \cdot \rangle_{\tiny{\mbox{cn}}}\) of $\mathbb{R}^{n}$, an invertible linear map \( g: \mathbb{R}^{n} \to \mathbb{R}^{n} \) naturally leads to an equivalence between the metric Lie algebras \((\mathbb{R}^{n},\mu, g^{\ast}\langle \cdot, \cdot \rangle_{\tiny{\mbox{cn}}})\) and \((\mathbb{R}^{n}, g\cdot\mu, \langle \cdot, \cdot \rangle_{\tiny{\mbox{cn}}})\). Here, \( g \) simultaneously acts as a Lie algebra isomorphism between \((\mathbb{R}^{n},\mu)\) and its transformed counterpart \((\mathbb{R}^{n}, g\cdot \mu) \), while preserving the metric structure by being a linear isometry between \((\mathbb{R}^{n}, g^{\ast}\langle \cdot,  \cdot \rangle_{\tiny{\mbox{cn}}})\) and \((\mathbb{R}^{n},\langle \cdot, \cdot \rangle_{\tiny{\mbox{cn}}})\).

In our context, the action of \(\operatorname{GL}(\mathbb{R}^{2n}, J_{\tiny{\mbox{cn}}})\) on \(\mu\) classifies all positive-definite inner products that are compatible with \(J_{\tiny{\mbox{cn}}}\). Specifically, for any Lie bracket of the form \(g \cdot \mu\), where \( g \in \operatorname{GL}(\mathbb{R}^{2n}, J_{\tiny{\mbox{cn}}})\), the complex structure \(J_{\tiny{\mbox{cn}}}\) remains well-defined on \((\mathbb{R}^{2n}, g \cdot \mu)\), naturally giving rise to a Hermitian structure \((J_{\tiny{\mbox{cn}}}, \langle \cdot, \cdot \rangle_{\tiny{\mbox{cn}}})\). The map \( g \) then induces a compatible inner product on \((\mathbb{R}^{2n}, \mu, J_{\tiny{\mbox{cn}}})\) via the pullback \( \langle \cdot, \cdot \rangle = g^* \langle \cdot, \cdot \rangle_{\tiny{\mbox{cn}}}\), ensuring that \( g \) acts as a \textit{complex isometry} between \((\mathbb{R}^{2n}, \mu, J_{\tiny{\mbox{cn}}}, \langle \cdot, \cdot \rangle)\) and \((\mathbb{R}^{2n}, g \cdot \mu, J_{\tiny{\mbox{cn}}}, \langle \cdot, \cdot \rangle_{\tiny{\mbox{cn}}})\).

Conversely, let \((V, J)\) be a real vector space endowed with a complex structure \(J\), and let \(\langle \cdot, \cdot \rangle\) be any positive-definite inner product on \(V\) that is compatible with \(J\). Viewing \((V, J)\) as a complex vector space via the action
\[
(a + b \sqrt{-1}) \cdot v = a v + b J v, \quad \text{for } v \in V, \quad a, b \in \mathbb{R},
\]
we can define a positive-definite Hermitian inner product
\[
\langle \! \langle \triangle , \triangledown \rangle \! \rangle = \langle \triangle , \triangledown \rangle + \sqrt{-1}  \langle J \triangle , \triangledown \rangle.
\]
By applying the Gram-Schmidt process, we obtain an ordered orthonormal basis \((u_1, \dots, u_n)\) for \((V, J, \langle \! \langle \triangle , \triangledown \rangle \! \rangle)\), which induces an ordered orthonormal basis for the Euclidean vector space \((V, \langle \cdot, \cdot \rangle)\) of the form \((u_1, \dots, u_n, u_{n+1}, \dots, u_{2n})\), where \(u_{n+i} = J u_i\) for \(i = 1, \dots, n\).

Therefore, if \(\langle \cdot, \cdot \rangle\) is compatible with \((\mathbb{R}^{2n}, \mu, J_{\tiny{\mbox{cn}}})\) and \((u_1, \dots, u_n, u_{n+1}, \dots, u_{2n})\) is the ordered basis described above, the linear map \(g: \mathbb{R}^{2n} \to \mathbb{R}^{2n}\) defined by \(g u_i = e_i\) for \(i = 1, \dots, 2n\) belongs to \(\operatorname{GL}(\mathbb{R}^{2n}, J_{\tiny{\mbox{cn}}})\). Moreover, \(g\) is a complex isometry between \((\mathbb{R}^{2n}, \mu, J_{\tiny{\mbox{cn}}}, \langle \cdot, \cdot \rangle)\) and \((\mathbb{R}^{2n}, g \cdot \mu, J_{\tiny{\mbox{cn}}}, \langle \cdot, \cdot \rangle_{\tiny{\mbox{cn}}})\).

The concepts introduced in the previous sections provide a framework for analyzing the evolution of left-invariant Hermitian structures on Lie groups under geometric flows, such as the \textit{Hermitian Curvature Flow} and the \textit{Chern-Ricci flow}. Moreover, since certain geometric quantities vary continuously as functions of tensors in a given space, our results can also be applied to the study of curvature-related aspects of Hermitian structures on Lie groups (see, for instance, \cite[Section 4]{FCRojitas} for the symplectic case). Additionally, they offer a way to infer properties of the cohomology and deformation theory of these structures (see Section \ref{defo}), among other significant applications.

Let us conclude this section with an application of the previously obtained results. For instance, since the only proper degeneration of \( (\mathfrak{r}'_{2},J_{1}(0,-1)) \) is \( (a_4,J) \), we can easily derive the following result:

\begin{theorem}
The four-dimensional Lie algebra \(\mathfrak{r}'_{2}\), equipped with its bi-invariant complex structure, admits a unique compatible metric, up to isometry and scaling.
\end{theorem}

\section{Deformations}\label{defo}

As we mentioned above, the study of the variety of Lie algebras endowed with complex structures has promising uses in homogeneous geometry. A significant aspect of this exploration is the existence of \textit{open orbits} in $\mathcal{L}_{J_{\text{cn}}}(\mathbb{R}^{2n})$. These open orbits provide a natural setting to identify what we might consider \textit{distinguished} complex structures on Lie algebras.

Consider a Lie algebra $(\mathbb{R}^{2n}, \mu)$  equipped with the complex structure $ J_{\tiny{\mbox{cn}}} $. We define the vector space of $2$-cocycles compatible with this structure as:
$$
Z_{2}^{J_{\text{cn}}}(\mathfrak{g}; \mathfrak{g}) = \left\{ \lambda \in C_{J_{\text{cn}}}^{2}(\mathbb{R}^{2n}; \mathbb{R}^{2n}) \mid \delta_{\mu} \lambda = 0 \right\},
$$
where the differential $\delta_{\mu}$, the \textit{linearization} of the Jacobiator, is given by:
\[
\delta_{\mu} \lambda(X_1, X_2, X_3) \coloneqq \sum_{\sigma \in S_3} \operatorname{sign}(\sigma) \left( \mu\big(\lambda(X_{\sigma(1)},X_{\sigma(2)}),X_{\sigma(3)}\big) + \lambda \big(\mu(X_{\sigma(1)},X_{\sigma(2)}),X_{\sigma(3)}\big) \right)
\]

Now, define $B_{2}^{J_{\text{cn}}}(\mathfrak{g}; \mathfrak{g})$ as the subspace of $2$-coboundaries, represented by elements of the form:
$$
\delta_{\mu} A = A \mu(\cdot, \cdot) - \mu(A \cdot, \cdot) - \mu(\cdot, A \cdot),
$$
for some $A \in \mathfrak{gl}(\mathbb{R}^{2n}, J_{\text{cn}})$. This subspace reflects the tangent space of the $\operatorname{GL}(\mathbb{R}^{2n}, J_{\text{cn}})$-orbit of $\mu$.

If the second cohomology group $H_{2}^{J_{\text{cn}}}(\mathfrak{g}; \mathfrak{g}) = Z_{2}^{J_{\text{cn}}}(\mathfrak{g}; \mathfrak{g}) / B_{2}^{J_{\text{cn}}}(\mathfrak{g}; \mathfrak{g})$ vanishes, then by the results in \cite[\S 24]{NijenhuisRichardson}, the orbit of $\mu$ under the action of $\operatorname{GL}(\mathbb{R}^{2n}, J_{\text{cn}})$ forms an open subset of the variety $\mathcal{L}_{J_{\text{cn}}}(\mathbb{R}^{2n})$. For example, it follows that the orbit $\operatorname{GL}(\mathbb{R}^{4}) \cdot \mu_4$, corresponding to the complex structure $J_2$ on $\mathfrak{r}'_2$, is an open set in $\mathcal{L}_{J_{\text{cn}}}(\mathbb{R}^{4})$ (which implies that $J_2$ is rigid in the sense of \cite{GiganteTomassini}).

This observation has practical implications: any Lie algebra with complex structure \( (\mathbb{R}^{2n}, \mu', J_{\text{cn}}) \) that is a small deformation of \( (\mathbb{R}^{4}, \mu_{4}, J_{\text{cn}}) \) must be holomorphically isomorphic to it. But what exactly is a \textit{deformation} of a Lie algebra endowed with a complex structure?

It is a deformation in the classical sense of Lie algebra deformations, with the additional requirement that the integrability condition imposed by the complex structure is preserved. More precisely, if \( \lambda \in C^{2}(\mathbb{R}^{2n}; \mathbb{R}^{2n}) \) is such that \( (\mathbb{R}^{2n}, \mu + \lambda, J_{\tiny{\mbox{cn}}}) \) defines a Lie algebra with complex structure, then we say that this is a deformation of \( (\mathbb{R}^{2n}, \mu, J_{\tiny{\mbox{cn}}}) \). The above definition is equivalent to requiring that \( \lambda \) satisfies the following conditions:
\begin{itemize}
\item (Complex compatibility) \( \lambda \in C_{J_{\text{cn}}}^{2}(\mathbb{R}^{2n}; \mathbb{R}^{2n}) \), ensuring that \( J_{\tiny{\mbox{cn}}} \) remains integrable.
\item (Deformation equation) The new bracket \( \mu + \lambda \) satisfies the Jacobi identity, which is equivalent to the well-known deformation equation:
   \[
   \delta_{\mu} \lambda + \frac{1}{2}[\lambda, \lambda] = 0.
   \]
\end{itemize}
The deformation is considered \textit{trivial} if it  \( (\mathbb{R}^{2n}, \mu + \lambda, J_{\tiny{\mbox{cn}}}) \) is holomorphically isomorphic to \( (\mathbb{R}^{2n}, \mu , J_{\tiny{\mbox{cn}}}) \).

The goal is to study solutions of this equation for \( \lambda \) close to zero. A common approach is to express \( \lambda \) as a formal power series, leading to the widely studied \textit{formal deformations}. Another type of solution arises by requiring that each term in the deformation equation vanishes independently, resulting in the so-called \textit{linear deformations} or \textit{infinitesimal deformations}. These correspond to 2-cocycles with respect to the adjoint representation of \( \mathfrak{g} := (\mathbb{R}^{2n},\mu) \) in the vector space $C_{J_{\text{cn}}}^{2}(\mathbb{R}^{2n}; \mathbb{R}^{2n})$  and satisfy the condition that \( (\mathbb{R}^{2n}, \lambda) \) is a Lie algebra. In this case, \( \mu + t \lambda \) defines a Lie algebra law for any \( t \in \mathbb{R} \). For instance, it follows from Equation (\ref{defotrivial}) that $(\mathbb{R}^{2n},\psi_{t,-t}(\mu),J_{\tiny{\mbox{cn}}})$ defines a linear deformation of $(\mathbb{R}^{2m},\mu,J_{\tiny{\mbox{cn}}})$ which is trivial.

Let \( (\mathbb{R}^{2n},\mu+\lambda_{t}, J_{\tiny{\mbox{cn}}}) \) be a continuous family with \( \lambda_0 = 0 \) and such that \( (\mathbb{R}^{2n},\mu+\lambda_t) \) is isomorphic to \( (\mathbb{R}^{2n},\mu) \) for any small \( t \).
In this case, even though the underlying Lie algebra remains unchanged, the family \( (\mathbb{R}^{2n},\mu+\lambda_{t}, J_{\tiny{\mbox{cn}}}) \) still represents a genuine deformation of the complex structure. This aligns with the classical perspective of deformations in complex geometry (Fr\"olicher-Kodaira-Nijenhuis-Spencer theory) where variations of the complex structure occur without altering the differentiable or topological structure of the space. Examples of this type include the one-parameter families of complex structures on \( \mathfrak{so}(3,\mathbb{R}) \times \mathbb{R} \) and \( \mathfrak{sl}(2,\mathbb{R}) \times \mathbb{R} \), arising from linear deformations of \( \mu_{8,+}(0) \) and \( \mu_{8,-}(0) \), respectively.

Finally, we note that the same analysis can be carried out in the subvariety of Lie algebras with abelian complex structures, where the space \( C_{J_{\text{cn}}}^{2}(\mathbb{R}^{2n}) \) must be replaced by the vector space \( C_{J_{\text{cn}}, \text{ab}}^{2}(\mathbb{R}^{2n}) \), defined as

\[
C_{J_{\text{cn}}, \text{ab}}^{2}(\mathbb{R}^{2n}) = \left\{\mu \in C^{2}(\mathbb{R}^{2n};\mathbb{R}^{2n})  : \mu(JX,JY) = \mu(X,Y) \right\}.
\]

%{
%\clearpage

%}


\begin{thebibliography}{BML}




\bibitem{ABDGH}
Romina Arroyo, Laura Barberis, Ver\'onica Diaz, Yamile Godoy, Isabel Hern\'{a}ndez:
\emph{Classification of nilpotent almost abelian Lie groups admitting left-invariant complex or symplectic structures.}
\textbf{arXiv:2406.06819 [math.DG]}



\bibitem{Bahturin}
Yuri Bahturin: %\newline
\emph{Identical Relations in Lie Algebras.} %\newline
De Gruyter Expositions in Mathematics %\newline
\textbf{68}.
Walter de Gruyter GmbH, Berlin/Boston
Second Edition. (2021)

\bibitem{Borel1}
Armand Borel: %\newline
\emph{Groupes Lineaires Algebriques.} %\newline
Annals of Mathematics
\textbf{64}, Number 1 (July 1956), 20--82


\bibitem{Borel2}
Armand Borel: %\newline
\emph{Linear Algebraic Groups.} %\newline
Graduate Texts in Mathematics
\textbf{126}.
Springer-Verlag New York (1991)

\bibitem{bump}
Daniel Bump: %\newline
\emph{Lie Groups.} %\newline_3
Graduate Texts in Mathematics \textbf{{225}},
Springer Science+Business Media, New York. Second Edition. (2013)


\bibitem{CaoZheng}
Kexiang Cao, Fangyang Zheng:
\emph{Fino-Vezzoni conjecture on Lie algebras with abelian ideals of codimension two.}
Mathematische Zeitschrift,
\textbf{Vol. 307}, Issue 2 (June 2024), Article: 31.







\bibitem{FCRojitas}
Edison Alberto Fern\'{a}ndez-Culma, Nadina Rojas: %\newline
\emph{Classification of Orbit Closures in the Variety of $4$-Dimensional Symplectic Lie Algebras.}
Algebras and Representation Theory
\textbf{Vol. 27}, Issue 2  (April 2024), 1013--1032.

\bibitem{FinoMainenti}
Anna Fino, Asia Mainenti:
\emph{A note on p-K\"{a}hler structures on compact quotients of Lie groups.}
Annali di Matematica Pura ed Applicata (1923 -)
\textbf{Vol. 203}, Issue 5 (October 2024), 2111--2124.





\bibitem{GiganteTomassini}
Giuliana Gigante, Giuseppe Tomassini: %\newline
\emph{ Deformations of Complex Structures on a Real Lie Algebra.}
In: Complex Analysis and Geometry. University Series in Mathematics (USMA).
Edited by: Vincenzo Ancona and Alessandro Silva.
Springer Science+Business Media, LLC. New York (1993).


\bibitem{GrunewaldOHalloran}
Fritz Grunewald and Joyce O'Halloran: %\newline
\emph{Varieties of nilpotent Lie algebras of dimension less than six.} %\newline
Journal of Algebra %\newline
\textbf{Vol. 112}, Issue 2 (February 1988), Pages 315--325.

\bibitem{Heber}
Jens Heber: %\newline
\emph{Noncompact homogeneous Einstein spaces.}
Inventiones mathematicae
\textbf{Vol. 133}, Issue 2 (July 1998), 279--352.



%\bibitem{HilgertNeeb}
%Joachim Hilgert and Karl-Hermann Neeb: %\newline
%\emph{Structure and Geometry of Lie Groups.}
%Springer Monographs in Mathematics,
%Springer New York Dordrecht Heidelberg London (2012)


\bibitem{Jablonski1}
Michael Jablonski: %\newline
\emph{Distinguished orbits of reductive groups.} %\newline
Rocky Mountain Journal of Mathematics
\textbf{Vol. 42}, Number 5 (October 2012), 1521--1549.

\bibitem{Kaygorodov}
Ivan Kaygorodov: %\newline
\emph{Non-associative algebraic structures: classification and structure.}
Communications in Mathematics
\textbf{Vol. 32 }, Number 3 (2024), 1--62.




\bibitem{LatorreUgarteVillacampa}
Adela Latorre, Luis Ugarte, Raquel Villacampa:
\emph{Complex structures on nilpotent Lie algebras with one-dimensional center.}
Journal of Algebra,
\textbf{Vol. 614} (January 2023), 271--306.



\bibitem{LauretRodriguezValencia}
Jorge Lauret, Edwin Alejandro Rodr\'iguez-Valencia: %\newline
\emph{On the Chern-Ricci flow and its solitons for Lie groups.}
Mathematische Nachrichten
\textbf{Vol. 288}, Issue 13 (September 2015), 1512--1526.


\bibitem{Lee2}
John M. Lee: %\newline
\emph{Introduction to Riemannian Manifolds.}
Graduate Texts in Mathematics
\textbf{176}, Springer International Publishing AG.
Second Edition. (2018)

\bibitem{Milnor}
John Milnor: %\newline
\emph{Curvatures of left invariant metrics on Lie groups.}
Advances in Mathematics
\textbf{Vol. 21}, Issue 3 (September 1976), 293--329.



\bibitem{nesterenko}
Maryna Nesterenko and Roman Popovych: %\newline
\emph{Contractions of low-dimensional Lie algebras}
Journal of Mathematical Physics
\textbf{Vol. 47}, Issue 12, 123515 (December 2006), 45 pages.



\bibitem{NijenhuisRichardson}
Albert Nijenhuis and Roger W. Richardson: %\newline
\emph{Cohomology and deformations in graded Lie algebras.}
Bulletin of the American Mathematical Society
\textbf{Vol. 72}, Number 1 (1966), 1--29.


\bibitem{ovando1}
Gabriela Ovando: %\newline
\emph{Invariant complex structures on solvable real Lie groups.}
manuscripta mathematica
\textbf{Vol. 103}, Issue 1 (September 2000), 19--30.

\bibitem{ovando2}
Gabriela Ovando: %\newline
\emph{Complex, Symplectic and K\"{a}hler Structures on Four Dimensional Lie Groups}
Revista de la Uni\'on Matem\'atica Argentina
\textbf{Vol. 45}, N\'umero 2 (2004), 55--67.


\bibitem{Reutenauer}
Christophe Reutenauer: %\newline
\emph{Free Lie algebras.} %\newline
London Mathematical Society Monographs %\newline
Clarendon Press Oxford (1993)



\bibitem{Rezaei-AghdamSephid}
Adel Rezaei-Aghdam, Mehdi Sephid: %\newline
\emph{Complex and bi-Hermitian structures on four-dimensional real Lie algebras.}
Journal of Physics A: Mathematical and Theoretical
\textbf{Vol. 43}, Number 32 (August 2010) 325210, 14pp.


\bibitem{snow}
Joanne Erdman Snow: %\newline
\emph{Invariant Complex Structures on Four-Dimensional Solvable Real Lie Groups.}
manuscripta mathematica
\textbf{Vol. 66}, Issue 1 (December 1990), 397--412.









\end{thebibliography}
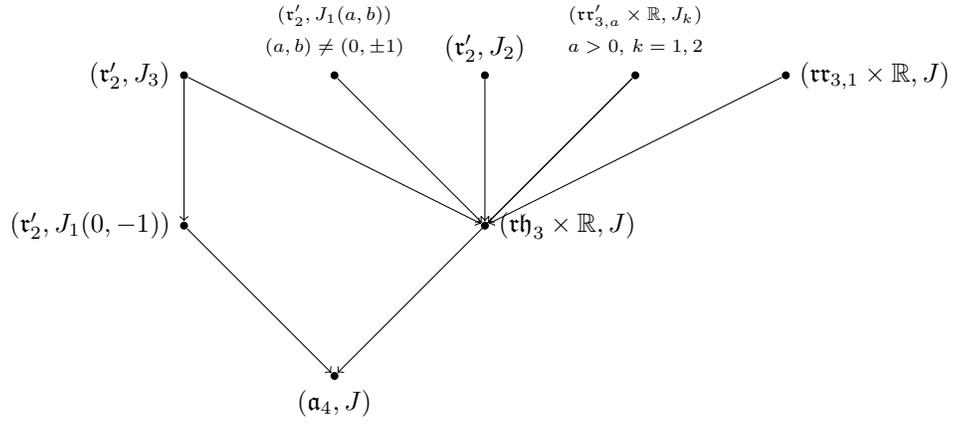
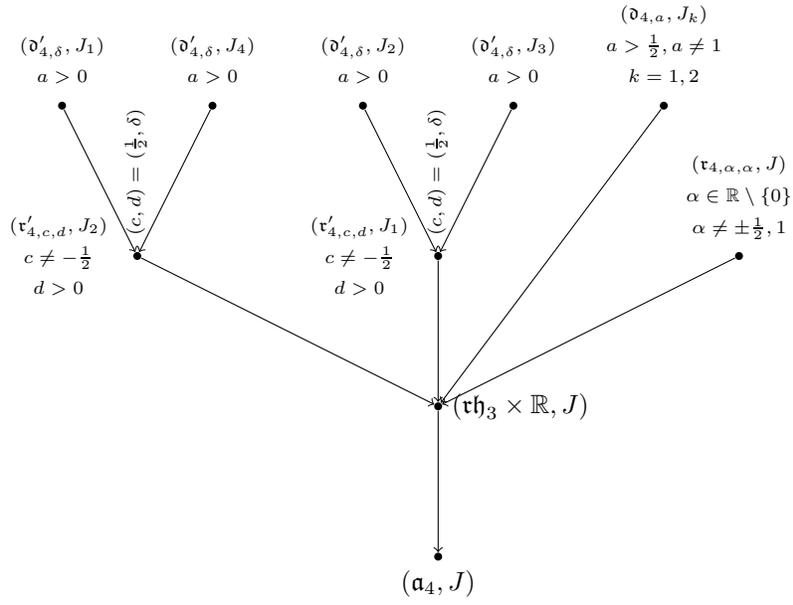
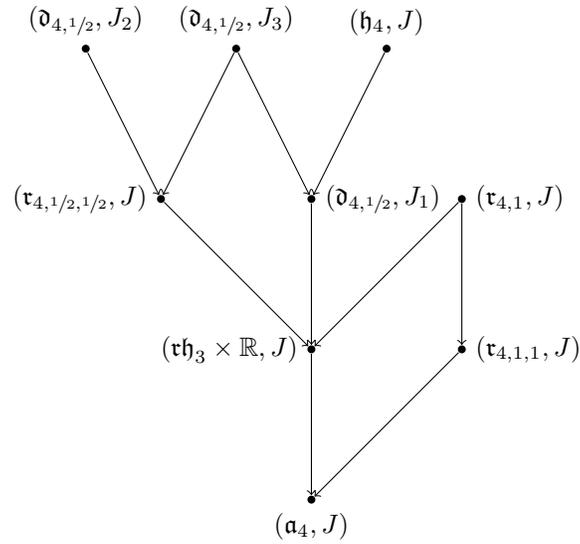
\end{document}

\section{Appendix}

\subsection{Invariants}

Let $C^{2}(\mathbb{R}^{2n};\mathbb{R}^{2n})$ and $C^{2}_{ J_{\tiny{\mbox{cn}}}}(\mathbb{R}^{2n};\mathbb{R}^{2n})$ denote the vector spaces:
\begin{eqnarray*}
C^{2}(\mathbb{R}^{2n};\mathbb{R}^{2n}) &:=&\left\{\mu:\mathbb{R}^{2n} \times \mathbb{R}^{2n} \rightarrow \mathbb{R}^{2n} : \mu \mbox{ is a bilinear map and skew-symmetric}\right\}\\
C^{2}_{ J_{\tiny{\mbox{cn}}} }(\mathbb{R}^{2n};\mathbb{R}^{2n}) &:=& \left\{\mu:\mathbb{R}^{2n} \times \mathbb{R}^{2n} \rightarrow \mathbb{R}^{2n} : \mu \in C^{2}(\mathbb{R}^{2n};\mathbb{R}^{2n}) \mbox{ and } \operatorname{N}_{J_{\tiny{\mbox{cn}}}}(\mu) \equiv 0 \right\}
\end{eqnarray*}

Here, ${J_{\tiny{\mbox{cn}}}}: \mathbb{R}^{2n} \rightarrow \mathbb{R}^{2n}$ is the linear map defined by ${J_{\tiny{\mbox{cn}}}} e_1 = e_2, {J_{\tiny{\mbox{cn}}}} e_3 = e_4, \ldots , {J_{\tiny{\mbox{cn}}}} e_{2n-1} = e_{2n}$
and ${J_{\tiny{\mbox{cn}}}}^{2} = -\operatorname{Id}$, and $\operatorname{N}_{{J_{\tiny{\mbox{cn}}}}}(\mu)$ is the Nijenhuis tensor associated to ${J_{\tiny{\mbox{cn}}}}$:
$$
\operatorname{N}_{{J_{\tiny{\mbox{cn}}}}}(\mu) (X,Y) = \mu(X,Y) + {J_{\tiny{\mbox{cn}}}}\mu({J_{\tiny{\mbox{cn}}}}X,Y) +{J_{\tiny{\mbox{cn}}}}\mu (X,{J_{\tiny{\mbox{cn}}}}Y) - \mu({J_{\tiny{\mbox{cn}}}}X,{J_{\tiny{\mbox{cn}}}}Y).
$$

Given $\alpha_{0},\ldots,\alpha_{7} \in \mathbb{R}$, let $\varphi_{\alpha_{0},\ldots,\alpha_{7}} : C^{2}_{{J_{\tiny{\mbox{cn}}}}}(\mathbb{R}^{2n};\mathbb{R}^{2n}) \rightarrow L^2(\mathbb{R}^{2n},\mathbb{R}^{2n})$ the
\(\operatorname{GL}(\mathbb{R}^{2n},J_{\tiny{\mbox{cn}}}  )\)-equivariant function to be defined by
$$
\varphi_{\alpha_{0},\ldots,\alpha_{7}} (\mu)(X,Y) =
\begin{array}{l}
\alpha_{ 0}\mu(X,Y) +\\
\alpha_{ 1}J\mu(X,Y)+
\alpha_{ 2}\mu(JX,Y)+
\alpha_{ 3}\mu(X,JY)+\\
\alpha_{ 4}J\mu(J X,Y)+
\alpha_{ 5}J\mu(X,J Y)+
\alpha_{ 6}\mu(J X,J Y)+\\
\alpha_{ 7}J\mu(J X,J Y),
\end{array}
$$
and set $\psi_{\alpha,\beta}:=\varphi_{1,\alpha,\beta,\beta,0,\ldots,0}$, $\theta_{\beta}:=\varphi_{0,1,\beta,\beta,0,\ldots,0}$.

Let $f,g \in \mathbb{R}$. We also consider the $\operatorname{GL}(\mathbb{R}^4)$-equivariant continuous functions $\kappa_{g}: C^{2}(\mathbb{R}^{m};\mathbb{R}^{m}) \rightarrow \Sigma^2(\mathbb{R}^{m}; \mathbb{R})$ and  $\Psi_{f,g}: C^{2}(\mathbb{R}^{m};\mathbb{R}^{m}) \rightarrow \Sigma^{3}(\mathbb{R}^{m};\mathbb{R})$ given by:
$$
\kappa_{g}(\lambda)(X,Y):= \operatorname{Trace}(\lambda(X, \lambda(Y , \square))) + g \operatorname{Trace}( \lambda(X, \square))\operatorname{Trace}( \lambda(Y, \square)),
$$
respectively
\begin{eqnarray*}
\Psi_{f,g}(\mu)(X_1,X_2,X_3) &=& \sum_{\sigma \in \mathfrak{S}_{3}}\operatorname{Trace}(\mu(X_{\sigma(1)}, \mu(X_{\sigma(2)}, \mu(X_{\sigma(3)}, \square )))) \\
&&  +f \, \kappa(\mu)(X_1,X_2) \operatorname{Trace}(\mu(X_3, \square)) +  f\, \kappa(\mu)(X_2,X_3) \operatorname{Trace}(\mu(X_1, \square)) \\
&& +  f\, \kappa(\mu)(X_3,X_1) \operatorname{Trace}(\mu(X_2, \square)) \\
&&+ g \operatorname{Trace}(\mu(X_1, \square)) \operatorname{Trace}(\mu(X_2, \square)) \operatorname{Trace}(\mu(X_3, \square))
\end{eqnarray*} where $\kappa(\mu)$ is the Killing form of $(\mathbb{R}^4,\mu)$.

The symmetric bilinear form $\kappa_{g}(\lambda)$ is just the \textit{modified Killing form} of $(\mathbb{R}^4, \lambda)$. Note that the sum $\displaystyle \sum_{\sigma \in \mathfrak{S}_{3}}\operatorname{Trace}(\mu(X_{\sigma(1)}, \mu(X_{\sigma(2)}, \mu(X_{\sigma(3)}, \square ))))$ produces six terms, which can be grouped into two distinct expressions appearing three times each, leading to the simplification $$ 3 \operatorname{Trace}(\mu(X_{1}, \mu(X_{2}, \mu(X_{3}, \square )))) + 3 \operatorname{Trace}(\mu(X_{2}, \mu(X_{1}, \mu(X_{3}, \square ))))$$.

\begin{center}
{\Huge Companion file}
\end{center}

\section{Invariants}

\subsection{$\mathfrak{rh}_{3}\times \mathbb{R}$}
It is unimodular and its (left)-invariant complex structure is abelian
$$
\left\{
\begin{array}{l}
[{\it e1},{\it e3}]={\it e2}
\end{array}
\right.
$$

\begin{enumerate}
\begin{multicols}{2}
  \item $\operatorname{Dim} \operatorname{Der} = 10 $
  \item $\operatorname{Dim} \operatorname{Der}_{J} = 4$
\end{multicols}
\end{enumerate}

%%%%%%%%%%%%%%%%%%%%%%%%%%%%%%%%%%%%%%%%%%%%%%%%%%%%%%%%%%%%%%%%%%%%%%%%%%%%%%%%%%%%%%%%%%%%%%%%%%%%%%%%%%%%%%%%%%%%%%%%%%%%%%%%%%%%%%%%%%%%%%%%%%%%

\subsection{$\mathfrak{ r}_{2 } \times \mathfrak{ r}_{2 }$}
$$
\left\{
\begin{array}{l}
[[{\it e1},{\it e3}]={\it e3},[{\it e2},{\it e4}]={\it e4}]
\end{array}
\right.
$$

It is non-unimodular and its (left)-invariant complex structure is abelian

\begin{enumerate}
\begin{multicols}{2}
  \item $\operatorname{Dim} \operatorname{Der} = { 4 }$
  \item $\operatorname{Dim} \operatorname{Der}_{J} = { 0 }$
\end{multicols}
\end{enumerate}

%%%%%%%%%%%%%%%%%%%%%%%%%%%%%%%%%%%%%%%%%%%%%%%%%%%%%%%%%%%%%%%%%%%%%%%%%%%%%%%%%%%%%%%%%%%%%%%%%%%%%%%%%%%%%%%%%%%%%%%%%%%%%%%%%%%%%%%%%%%%%%%%%%%%

\subsection{ $(\mathfrak{r}^{'}_{2 }, J_{1}(a,b))$ with $b\neq0$ and $a\neq0$}
$$
\left\{
\begin{array}{l}
[{\it e1},{\it e2}]={\it e2},
[{\it e1},{\it e4}]={\it e4},
[{\it e2},{\it e3}]=a{\it e2}+b{\it e4},
[{\it e3},{\it e4}]=b{\it e2}-a{\it e4}
\end{array}
\right.
$$

It is non-unimodular and its (left)-invariant complex structure is non-abelian

\begin{enumerate}
\begin{multicols}{2}
  \item $\operatorname{Dim} \operatorname{Der} = { 4 }$
  \item $\operatorname{Dim} \operatorname{Der}_{J} = {  2 }$
\end{multicols}
%Hay algo con 4a^2+2b^2+2b = 0
\end{enumerate}

%%%%%%%%%%%%%%%%%%%%%%%%%%%%%%%%%%%%%%%%%%%%%%%%%%%%%%%%%%%%%%%%%%%%%%%%%%%%%%%%%%%%%%%%%%%%%%%%%%%%%%%%%%%%%%%%%%%%%%%%%%%%%%%%%%%%%%%%%%%%%%%%%%%%

\subsection{$(\mathfrak{r}^{'}_{2 }, J_{1}(a,b))$ with $a= 0$ and $b \neq \pm 1,0$}
$$
\left\{
\begin{array}{l}
[{\it e1},{\it e2}]={\it e2},
[{\it e1},{\it e4}]={\it e4},
[{\it e2},{\it e3}]=b{\it e4},
[{\it e3},{\it e4}]=b{\it e2}
\end{array}
\right.
$$

It is non-unimodular and its (left)-invariant complex structure is non-abelian

\begin{enumerate}
\begin{multicols}{2}
  \item $\operatorname{Dim} \operatorname{Der} = { 4 }$
  \item $\operatorname{Dim} \operatorname{Der}_{J} = {  2 }$
\end{multicols}
\end{enumerate}

%%%%%%%%%%%%%%%%%%%%%%%%%%%%%%%%%%%%%%%%%%%%%%%%%%%%%%%%%%%%%%%%%%%%%%%%%%%%%%%%%%%%%%%%%%%%%%%%%%%%%%%%%%%%%%%%%%%%%%%%%%%%%%%%%%%%%%%%%%%%%%%%%%%%

\subsection{$(\mathfrak{r}^{'}_{2 }, J_{1}(a,b))$ with $a= 0$ and $b = 1$ }
$$
\left\{
\begin{array}{l}
[{\it e1},{\it e2}]={\it e2},
[{\it e1},{\it e4}]={\it e4},
[{\it e2},{\it e3}]={\it e4},
[{\it e3},{\it e4}]={\it e2}
\end{array}
\right.
$$

It is non-unimodular and its (left)-invariant complex structure is abelian

\begin{enumerate}
\begin{multicols}{2}
  \item $\operatorname{Dim} \operatorname{Der} = { 4 }$
  \item $\operatorname{Dim} \operatorname{Der}_{J} = { 2 }$
\end{multicols}
\end{enumerate}

%%%%%%%%%%%%%%%%%%%%%%%%%%%%%%%%%%%%%%%%%%%%%%%%%%%%%%%%%%%%%%%%%%%%%%%%%%%%%%%%%%%%%%%%%%%%%%%%%%%%%%%%%%%%%%%%%%%%%%%%%%%%%%%%%%%%%%%%%%%%%%%%%%%%

\subsection{$(\mathfrak{r}^{'}_{2 }, J_{1}(a,b))$ with $a= 0$ and $b = -1$ }
$$
\left\{
\begin{array}{l}
[{\it e1},{\it e2}]={\it e2},
[{\it e1},{\it e4}]={\it e4},
[{\it e2},{\it e3}]=-{\it e4},
[{\it e3},{\it e4}]=-{\it e2}
\end{array}
\right.
$$

It is non-unimodular and its (left)-invariant complex structure is bi-invariant.

\begin{enumerate}
\begin{multicols}{2}
  \item $\operatorname{Dim} \operatorname{Der} = {  4 }$
  \item $\operatorname{Dim} \operatorname{Der}_{J} = { 4 }$
\end{multicols}
\end{enumerate}

%%%%%%%%%%%%%%%%%%%%%%%%%%%%%%%%%%%%%%%%%%%%%%%%%%%%%%%%%%%%%%%%%%%%%%%%%%%%%%%%%%%%%%%%%%%%%%%%%%%%%%%%%%%%%%%%%%%%%%%%%%%%%%%%%%%%%%%%%%%%%%%%%%%%

\subsection{$(\mathfrak{ r}^{'}_{2 },J_2)$}
$$
\left\{
[{\it e1},{\it e3}]={\it e3},
[{\it e1},{\it e4}]={\it e4},
[{\it e2},{\it e3}]={\it e4},
[{\it e2},{\it e4}]=-{\it e3}
\right.
$$

It is non-unimodular and its (left)-invariant complex structure is abelian

\begin{enumerate}
\begin{multicols}{2}
  \item $\operatorname{Dim} \operatorname{Der} = { 4 }$
  \item $\operatorname{Dim} \operatorname{Der}_{J} = { 0 }$
\end{multicols}
\end{enumerate}

%%%%%%%%%%%%%%%%%%%%%%%%%%%%%%%%%%%%%%%%%%%%%%%%%%%%%%%%%%%%%%%%%%%%%%%%%%%%%%%%%%%%%%%%%%%%%%%%%%%%%%%%%%%%%%%%%%%%%%%%%%%%%%%%%%%%%%%%%%%%%%%%%%%%

\subsection{$(\mathfrak{ r}^{'}_{2 },J_3)$}
$$
\left\{
[{\it e1},{\it e2}]={\it e2},
[{\it e1},{\it e3}]={\it e4},
[{\it e1},{\it e4}]={\it e4},
[{\it e2},{\it e3}]=-{\it e4},
[{\it e3},{\it e4}]=-{\it e2}
\right.
$$

It is non-unimodular and its (left)-invariant complex structure is non-abelian

\begin{enumerate}
\begin{multicols}{2}
  \item $\operatorname{Dim} \operatorname{Der} = { 4 }$
  \item $\operatorname{Dim} \operatorname{Der}_{J} = { 2 }$
\end{multicols}
\end{enumerate}
%%%%%%%%%%%%%%%%%%%%%%%%%%%%%%%%%%%%%%%%%%%%%%%%%%%%%%%%%%%%%%%%%%%%%%%%%%%%%%%%%%%%%%%%%%%%%%%%%%%%%%%%%%%%%%%%%%%%%%%%%%%%%%%%%%%%%%%%%%%%%%%%%%%%

\subsection{$\mathfrak{rr}_{3,0 } \times \mathbb{R} \cong \mathfrak{r}_{2} \times \mathbb{R}^{2}$}
$$
\left\{
\begin{array}{l}
[{\it e1},{\it e3}]={\it e3}
\end{array}
\right.
$$

It is non unimodular and its (left)-invariant complex structure is abelian

\begin{enumerate}
\begin{multicols}{2}
  \item $\operatorname{Dim} \operatorname{Der} = 8$
  \item $\operatorname{Dim} \operatorname{Der}_{J} = 2$
\end{multicols}
\end{enumerate}

%%%%%%%%%%%%%%%%%%%%%%%%%%%%%%%%%%%%%%%%%%%%%%%%%%%%%%%%%%%%%%%%%%%%%%%%%%%%%%%%%%%%%%%%%%%%%%%%%%%%%%%%%%%%%%%%%%%%%%%%%%%%%%%%%%%%%%%%%%%%%%%%%%%%

\subsection{$\mathfrak{rr}_{3,1 }\times \mathbb{R}$}
$$
\left\{
\begin{array}{l}
[{\it e1},{\it e2}]={\it e2},[{\it e1},{\it e4}]={\it e4}
\end{array}
\right.
$$

It is non-unimodular and its (left)-invariant complex structure is non-abelian

\begin{enumerate}
\begin{multicols}{2}
  \item $\operatorname{Dim} \operatorname{Der} = 8$
  \item $\operatorname{Dim} \operatorname{Der}_{J} = 2$
\end{multicols}
\end{enumerate}

%%%%%%%%%%%%%%%%%%%%%%%%%%%%%%%%%%%%%%%%%%%%%%%%%%%%%%%%%%%%%%%%%%%%%%%%%%%%%%%%%%%%%%%%%%%%%%%%%%%%%%%%%%%%%%%%%%%%%%%%%%%%%%%%%%%%%%%%%%%%%%%%%%%%

\subsection{$(\mathfrak{ rr}^{'}_{3,a } \times \mathbb{R} , J_{1})$, $a>0$}
$$
\left\{
\begin{array}{l}
[{\it e1},{\it e2}]=a{\it e2}+{\it e4},[{\it e1},{\it e4}]=-{\it e2}+a{\it e4}
\end{array}
\right.
$$

It is non-unimodular and its (left)-invariant complex structure is non-abelian

\begin{enumerate}
\begin{multicols}{2}
  \item $\operatorname{Dim} \operatorname{Der} = { 6}$
  \item $\operatorname{Dim} \operatorname{Der}_{J} = { 2 }$
\end{multicols}
\end{enumerate}

%%%%%%%%%%%%%%%%%%%%%%%%%%%%%%%%%%%%%%%%%%%%%%%%%%%%%%%%%%%%%%%%%%%%%%%%%%%%%%%%%%%%%%%%%%%%%%%%%%%%%%%%%%%%%%%%%%%%%%%%%%%%%%%%%%%%%%%%%%%%%%%%%%%%

\subsection{$(\mathfrak{ rr}^{'}_{3,a }\times \mathbb{R} , J_{2})$, $a>0$}
$$
\left\{
\begin{array}{l}
[{\it e1},{\it e2}]=-a{\it e2}+{\it e4},
[{\it e1},{\it e4}]=-{\it e2}-a{\it e4}
\end{array}
\right.
$$

It is non-unimodular and its (left)-invariant complex structure is non-abelian

\begin{enumerate}
\begin{multicols}{2}
  \item $\operatorname{Dim} \operatorname{Der} = { 6 }$
  \item $\operatorname{Dim} \operatorname{Der}_{J} = { 2 }$
\end{multicols}
\end{enumerate}

%%%%%%%%%%%%%%%%%%%%%%%%%%%%%%%%%%%%%%%%%%%%%%%%%%%%%%%%%%%%%%%%%%%%%%%%%%%%%%%%%%%%%%%%%%%%%%%%%%%%%%%%%%%%%%%%%%%%%%%%%%%%%%%%%%%%%%%%%%%%%%%%%%%%

\subsection{$\mathfrak{ rr}^{'}_{3,0 }\times \mathbb{R}$}
$$
\left\{
\begin{array}{l}
[{\it e1},{\it e2}]={\it e4},
[{\it e1},{\it e4}]=-{\it e2}
\end{array}
\right.
$$

It is unimodular and its (left)-invariant complex structure is non-abelian

\begin{enumerate}
\begin{multicols}{2}
  \item $\operatorname{Dim} \operatorname{Der} = { 6}$
  \item $\operatorname{Dim} \operatorname{Der}_{J} = { 2}$
\end{multicols}
\end{enumerate}

%%%%%%%%%%%%%%%%%%%%%%%%%%%%%%%%%%%%%%%%%%%%%%%%%%%%%%%%%%%%%%%%%%%%%%%%%%%%%%%%%%%%%%%%%%%%%%%%%%%%%%%%%%%%%%%%%%%%%%%%%%%%%%%%%%%%%%%%%%%%%%%%%%%%

\subsection{$(\mathfrak{ so}(3,\mathbb{R}) \times \mathbb{R}, J(t))$}
$$
[{\it e1},{\it e2}]={\it e3},
[{\it e1},{\it e3}]=-{\it e2}+t{\it e4},
[{\it e2},{\it e3}]={\it e1}
$$

It is unimodular and its (left)-invariant complex structure is non-abelian

\begin{enumerate}
\begin{multicols}{2}
  \item $\operatorname{Dim} \operatorname{Der} = { 4 }$
  \item $\operatorname{Dim} \operatorname{Der}_{J} = { 1 }$
\end{multicols}
\end{enumerate}

%%%%%%%%%%%%%%%%%%%%%%%%%%%%%%%%%%%%%%%%%%%%%%%%%%%%%%%%%%%%%%%%%%%%%%%%%%%%%%%%%%%%%%%%%%%%%%%%%%%%%%%%%%%%%%%%%%%%%%%%%%%%%%%%%%%%%%%%%%%%%%%%%%%%

\subsection{$(\mathfrak{ sl}(2,\mathbb{R}) \times \mathbb{R}, J(t))$}
$$
[{\it e1},{\it e2}]={\it e3},
[{\it e1},{\it e3}]={\it e2}+t{\it e4},
[{\it e2},{\it e3}]={\it e1}
$$

It is unimodular and its (left)-invariant complex structure is non-abelian

\begin{enumerate}
\begin{multicols}{2}
  \item $\operatorname{Dim} \operatorname{Der} = { 4 }$
  \item $\operatorname{Dim} \operatorname{Der}_{J} = { 1 }$
\end{multicols}
\end{enumerate}

%%%%%%%%%%%%%%%%%%%%%%%%%%%%%%%%%%%%%%%%%%%%%%%%%%%%%%%%%%%%%%%%%%%%%%%%%%%%%%%%%%%%%%%%%%%%%%%%%%%%%%%%%%%%%%%%%%%%%%%%%%%%%%%%%%%%%%%%%%%%%%%%%%%%

\subsection{$\mathfrak{ r}_{4,1 }$}
$$
[{\it e1},{\it e2}]={\it e2},
[{\it e1},{\it e3}]={\it e3}+{\it e4},
[{\it e1},{\it e4}]={\it e4}
$$

It is non-unimodular and its (left)-invariant complex structure is non-abelian

\begin{enumerate}
\begin{multicols}{2}
  \item $\operatorname{Dim} \operatorname{Der} = { 8 }$
  \item $\operatorname{Dim} \operatorname{Der}_{J} = { 2 }$
\end{multicols}
\end{enumerate}

%%%%%%%%%%%%%%%%%%%%%%%%%%%%%%%%%%%%%%%%%%%%%%%%%%%%%%%%%%%%%%%%%%%%%%%%%%%%%%%%%%%%%%%%%%%%%%%%%%%%%%%%%%%%%%%%%%%%%%%%%%%%%%%%%%%%%%%%%%%%%%%%%%%%

%a \in \mathbb{R}^{\ast}
\subsection{$\mathfrak{ r}_{4,a,a }$ with $a\in \mathbb{R}$, $a\neq 0, -\frac{1}{2}, 1$}
$$
\left\{
\begin{array}{l}
[{\it e1},{\it e2}]=a{\it e2},
[{\it e1},{\it e3}]={\it e3},
[{\it e1},{\it e4}]=a{\it e4}
\end{array}
\right.
$$

It is non-unimodular and its (left)-invariant complex structure is non-abelian

\begin{enumerate}
\begin{multicols}{2}
  \item $\operatorname{Dim} \operatorname{Der} = { 8 }$
  \item $\operatorname{Dim} \operatorname{Der}_{J} = { 2 }$
\end{multicols}
\end{enumerate}

%%%%%%%%%%%%%%%%%%%%%%%%%%%%%%%%%%%%%%%%%%%%%%%%%%%%%%%%%%%%%%%%%%%%%%%%%%%%%%%%%%%%%%%%%%%%%%%%%%%%%%%%%%%%%%%%%%%%%%%%%%%%%%%%%%%%%%%%%%%%%%%%%%%%

\subsection{$\mathfrak{ r}_{4,a,a }$ with $a = -1/2$}
$$
[{\it e1},{\it e2}]=-\frac{1}{2}\,{\it e2},
[{\it e1},{\it e3}]={\it e3},
[{\it e1},{\it e4}]=-\frac{1}{2}\,{\it e4}
$$

It is unimodular and its (left)-invariant complex structure is non-abelian

\begin{enumerate}
\begin{multicols}{2}
  \item $\operatorname{Dim} \operatorname{Der} = { 8 }$
  \item $\operatorname{Dim} \operatorname{Der}_{J} = { 2 }$
\end{multicols}
\end{enumerate}

%%%%%%%%%%%%%%%%%%%%%%%%%%%%%%%%%%%%%%%%%%%%%%%%%%%%%%%%%%%%%%%%%%%%%%%%%%%%%%%%%%%%%%%%%%%%%%%%%%%%%%%%%%%%%%%%%%%%%%%%%%%%%%%%%%%%%%%%%%%%%%%%%%%%

%\subsection{$\mathfrak{ r}_{4,a,a }$ with $a = 1/2$}
%$$
%\left\{
%\begin{array}{l}
%[{\it e1},{\it e2}]={\it e2},[{\it e1},{\it e3}]=1/2\,{\it e3},[{\it
%e1},{\it e4}]=1/2\,{\it e4}
%\end{array}
%\right.
%$$
%
%It is non-unimodular and its (left)-invariant complex structure is non-abelian
%
%\begin{enumerate}
%\begin{multicols}{2}
%  \item $\operatorname{Dim} \operatorname{Der} = { 8 }$
%  \item $\operatorname{Dim} \operatorname{Der}_{J} = { 2 }$
%\end{multicols}
%\end{enumerate}

%%%%%%%%%%%%%%%%%%%%%%%%%%%%%%%%%%%%%%%%%%%%%%%%%%%%%%%%%%%%%%%%%%%%%%%%%%%%%%%%%%%%%%%%%%%%%%%%%%%%%%%%%%%%%%%%%%%%%%%%%%%%%%%%%%%%%%%%%%%%%%%%%%%%

\subsection{$\mathfrak{ r}_{4,a,a }$ with $a = 1$}
$$
[{\it e1},{\it e2}]={\it e2},
[{\it e1},{\it e3}]={\it e3},
[{\it e1},{\it e4}]={\it e4}
$$

It is non-unimodular and its (left)-invariant complex structure is non-abelian

\begin{enumerate}
\begin{multicols}{2}
  \item $\operatorname{Dim} \operatorname{Der} = { 12 }$
  \item $\operatorname{Dim} \operatorname{Der}_{J} = { 4 }$
\end{multicols}
\end{enumerate}

%%%%%%%%%%%%%%%%%%%%%%%%%%%%%%%%%%%%%%%%%%%%%%%%%%%%%%%%%%%%%%%%%%%%%%%%%%%%%%%%%%%%%%%%%%%%%%%%%%%%%%%%%%%%%%%%%%%%%%%%%%%%%%%%%%%%%%%%%%%%%%%%%%%%

\subsection{$(\mathfrak{ r}^{'}_{4,a,b },J_1)$ with $a\neq -\frac{1}{2}$, $b>0$}
$$
\left\{
\begin{array}{l}
[{\it e1},{\it e2}]=a{\it e2}+b{\it e4},
[{\it e1},{\it e3}]={\it e3},
[{\it e1},{\it e4}]=-b{\it e2}+a{\it e4}
\end{array}
\right.
$$

It is non-unimodular and its (left)-invariant complex structure is non-abelian

\begin{enumerate}
\begin{multicols}{2}
  \item $\operatorname{Dim} \operatorname{Der} = {6 }$
  \item $\operatorname{Dim} \operatorname{Der}_{J} = { 2}$
\end{multicols}
\end{enumerate}

%%%%%%%%%%%%%%%%%%%%%%%%%%%%%%%%%%%%%%%%%%%%%%%%%%%%%%%%%%%%%%%%%%%%%%%%%%%%%%%%%%%%%%%%%%%%%%%%%%%%%%%%%%%%%%%%%%%%%%%%%%%%%%%%%%%%%%%%%%%%%%%%%%%%

\subsection{$(\mathfrak{ r}^{'}_{4,a,b },J_1)$ with $a=-\frac{1}{2}$, $b>0$}
$$
\left\{
\begin{array}{l}
[{\it e1},{\it e2}]=-\frac{1}{2}\,{\it e2}+b{\it e4},
[{\it e1},{\it e3}]={\it e3},
[{\it e1},{\it e4}]=-b{\it e2}-\frac{1}{2}\,{\it e4}
\end{array}
\right.
$$

It is unimodular and its (left)-invariant complex structure is non-abelian

\begin{enumerate}
\begin{multicols}{2}
  \item $\operatorname{Dim} \operatorname{Der} = { 6 }$
  \item $\operatorname{Dim} \operatorname{Der}_{J} = {  2}$
\end{multicols}
\end{enumerate}

%%%%%%%%%%%%%%%%%%%%%%%%%%%%%%%%%%%%%%%%%%%%%%%%%%%%%%%%%%%%%%%%%%%%%%%%%%%%%%%%%%%%%%%%%%%%%%%%%%%%%%%%%%%%%%%%%%%%%%%%%%%%%%%%%%%%%%%%%%%%%%%%%%%%

\subsection{$(\mathfrak{ r}^{'}_{4,a,b },J_2)$ with $a\neq-\frac{1}{2}$, $b>0$}
$$
\left\{
\begin{array}{l}
[{\it e1},{\it e2}]=a{\it e2}-b{\it e4},
[{\it e1},{\it e3}]={\it e3},
[{\it e1},{\it e4}]=b{\it e2}+a{\it e4}
\end{array}
\right.
$$

It is non-unimodular and its (left)-invariant complex structure is non-abelian

\begin{enumerate}
\begin{multicols}{2}
  \item $\operatorname{Dim} \operatorname{Der} = { 6}$
  \item $\operatorname{Dim} \operatorname{Der}_{J} = { 2 }$
\end{multicols}
\end{enumerate}

%%%%%%%%%%%%%%%%%%%%%%%%%%%%%%%%%%%%%%%%%%%%%%%%%%%%%%%%%%%%%%%%%%%%%%%%%%%%%%%%%%%%%%%%%%%%%%%%%%%%%%%%%%%%%%%%%%%%%%%%%%%%%%%%%%%%%%%%%%%%%%%%%%%%

\subsection{$(\mathfrak{ r}^{'}_{4,\alpha,\beta },J_2)$ with $\alpha=- \frac{1}{2}$}
$$
\left\{
\begin{array}{l}
[{\it e1},{\it e2}]=-\frac{1}{2}\,{\it e2}-b{\it e4},
[{\it e1},{\it e3}]={\it e3},
[{\it e1},{\it e4}]=b{\it e2}-\frac{1}{2}\,{\it e4}

\end{array}
\right.
$$

It is unimodular and its (left)-invariant complex structure is non-abelian

\begin{enumerate}
\begin{multicols}{2}
  \item $\operatorname{Dim} \operatorname{Der} = { 6 }$
  \item $\operatorname{Dim} \operatorname{Der}_{J} = { 2 }$
\end{multicols}
\end{enumerate}

%%%%%%%%%%%%%%%%%%%%%%%%%%%%%%%%%%%%%%%%%%%%%%%%%%%%%%%%%%%%%%%%%%%%%%%%%%%%%%%%%%%%%%%%%%%%%%%%%%%%%%%%%%%%%%%%%%%%%%%%%%%%%%%%%%%%%%%%%%%%%%%%%%%%

\subsection{$(\mathfrak{d }_{ 4},J_1)$}
$$
\left\{
\begin{array}{l}
[{\it e1},{\it e2}]=-{\it e2},
[{\it e1},{\it e3}]={\it e3},
[{\it e2},{\it e3}]={\it e4}
\end{array}
\right.
$$

It is unimodular and its (left)-invariant complex structure is non-abelian

\begin{enumerate}
\begin{multicols}{2}
  \item $\operatorname{Dim} \operatorname{Der} = { 5 }$
  \item $\operatorname{Dim} \operatorname{Der}_{J} = { 2 }$
\end{multicols}
\end{enumerate}

%%%%%%%%%%%%%%%%%%%%%%%%%%%%%%%%%%%%%%%%%%%%%%%%%%%%%%%%%%%%%%%%%%%%%%%%%%%%%%%%%%%%%%%%%%%%%%%%%%%%%%%%%%%%%%%%%%%%%%%%%%%%%%%%%%%%%%%%%%%%%%%%%%%%
\subsection{$(\mathfrak{d }_{ 4},J_2)$}
$$
\left\{
\begin{array}{l}
[{\it e1},{\it e2}]=-{\it e2},
[{\it e1},{\it e3}]={\it e3}+{\it e4},
[{\it e2},{\it e3}]={\it e4}
\end{array}
\right.
$$

It is unimodular and its (left)-invariant complex structure is non-abelian

\begin{enumerate}
\begin{multicols}{2}
  \item $\operatorname{Dim} \operatorname{Der} = { 5 }$
  \item $\operatorname{Dim} \operatorname{Der}_{J} = { 1 }$
\end{multicols}
\end{enumerate}

%%%%%%%%%%%%%%%%%%%%%%%%%%%%%%%%%%%%%%%%%%%%%%%%%%%%%%%%%%%%%%%%%%%%%%%%%%%%%%%%%%%%%%%%%%%%%%%%%%%%%%%%%%%%%%%%%%%%%%%%%%%%%%%%%%%%%%%%%%%%%%%%%%%%

\subsection{$(\mathfrak{ d}_{4, \lambda}$ with $ \lambda \neq 1, \frac{1}{2} , J_1)$ and $ \lambda >1/2$}
$$
\left\{
\begin{array}{l}
[[{\it e1},{\it e2}]={\frac {\lambda}{1-\lambda}}{\it e2},[{\it e1},
{\it e3}]={\it e3},[{\it e1},{\it e4}]={\frac {1}{1-\lambda}}{\it e4},[
{\it e2},{\it e3}]={\it e4}]

\end{array}
\right.
$$

It is non-unimodular and its (left)-invariant complex structure is non-abelian

\begin{enumerate}
\begin{multicols}{2}
  \item $\operatorname{Dim} \operatorname{Der} = { 5 }$
  \item $\operatorname{Dim} \operatorname{Der}_{J} = {  1 }$
\end{multicols}
\end{enumerate}

%%%%%%%%%%%%%%%%%%%%%%%%%%%%%%%%%%%%%%%%%%%%%%%%%%%%%%%%%%%%%%%%%%%%%%%%%%%%%%%%%%%%%%%%%%%%%%%%%%%%%%%%%%%%%%%%%%%%%%%%%%%%%%%%%%%%%%%%%%%%%%%%%%%%

\subsection{$(\mathfrak{ d}_{4, \lambda}$ with $ \lambda \neq 1, \frac{1}{2} , J_2)$ and $\lambda >1/2$}
$$
\left\{
\begin{array}{l}
[[{\it e1},{\it e2}]={\frac { \left( 1-\lambda \right) }{
\lambda}}{\it e2},[{\it e1},{\it e3}]={\it e3},[{\it e1},{\it e4}]={\frac {1}{\lambda}}{
\it e4},[{\it e2},{\it e3}]={\it e4}]
\end{array}
\right.
$$

It is non-unimodular and its (left)-invariant complex structure is non-abelian

\begin{enumerate}
\begin{multicols}{2}
  \item $\operatorname{Dim} \operatorname{Der} = { 5 }$
  \item $\operatorname{Dim} \operatorname{Der}_{J} = { 1 }$
\end{multicols}
\end{enumerate}

%%%%%%%%%%%%%%%%%%%%%%%%%%%%%%%%%%%%%%%%%%%%%%%%%%%%%%%%%%%%%%%%%%%%%%%%%%%%%%%%%%%%%%%%%%%%%%%%%%%%%%%%%%%%%%%%%%%%%%%%%%%%%%%%%%%%%%%%%%%%%%%%%%%%

\subsection{$(\mathfrak{ d}_{4,1 },J)$}
$$
\left\{
\begin{array}{l}
[{\it e1},{\it e3}]={\it e3},
[{\it e1},{\it e4}]={\it e4},
[{\it e2},{\it e3}]={\it e4}
\end{array}
\right.
$$

It is non-unimodular and its (left)-invariant complex structure is abelian

\begin{enumerate}
\begin{multicols}{2}
  \item $\operatorname{Dim} \operatorname{Der} = { 5 }$
  \item $\operatorname{Dim} \operatorname{Der}_{J} = { 1 }$
\end{multicols}
\end{enumerate}

%%%%%%%%%%%%%%%%%%%%%%%%%%%%%%%%%%%%%%%%%%%%%%%%%%%%%%%%%%%%%%%%%%%%%%%%%%%%%%%%%%%%%%%%%%%%%%%%%%%%%%%%%%%%%%%%%%%%%%%%%%%%%%%%%%%%%%%%%%%%%%%%%%%%

\subsection{$(\mathfrak{ d}_{4,\frac{1}{2} },J_1)$}
$$
\left\{
\begin{array}{l}
[{\it e1},{\it e2}]={\it e2},
[{\it e1},{\it e3}]={\it e3},
[{\it e1},{\it e4}]=2\,{\it e4},
[{\it e2},{\it e3}]={\it e4}
\end{array}
\right.
$$

It is non-unimodular and its (left)-invariant complex structure is non-abelian

\begin{enumerate}
\begin{multicols}{2}
  \item $\operatorname{Dim} \operatorname{Der} = { 7 }$
  \item $\operatorname{Dim} \operatorname{Der}_{J} = { 2 }$
\end{multicols}
\end{enumerate}

%%%%%%%%%%%%%%%%%%%%%%%%%%%%%%%%%%%%%%%%%%%%%%%%%%%%%%%%%%%%%%%%%%%%%%%%%%%%%%%%%%%%%%%%%%%%%%%%%%%%%%%%%%%%%%%%%%%%%%%%%%%%%%%%%%%%%%%%%%%%%%%%%%%%

\subsection{$(\mathfrak{ d}_{4,\frac{1}{2} },J_2)$}
$$
\left\{
\begin{array}{l}
[{\it e1},{\it e2}]=\frac{1}{2}\,{\it e2},
[{\it e1},{\it e3}]={\it e3},
[{\it e1},{\it e4}]=\frac{1}{2}\,{\it e4},
[{\it e2},{\it e4}]={\it e3}
\end{array}
\right.
$$

It is non-unimodular and its (left)-invariant complex structure is non-abelian

\begin{enumerate}
\begin{multicols}{2}
  \item $\operatorname{Dim} \operatorname{Der} = { 7 }$
  \item $\operatorname{Dim} \operatorname{Der}_{J} = { 1 }$
\end{multicols}
\end{enumerate}

%%%%%%%%%%%%%%%%%%%%%%%%%%%%%%%%%%%%%%%%%%%%%%%%%%%%%%%%%%%%%%%%%%%%%%%%%%%%%%%%%%%%%%%%%%%%%%%%%%%%%%%%%%%%%%%%%%%%%%%%%%%%%%%%%%%%%%%%%%%%%%%%%%%%

\subsection{$(\mathfrak{ d}_{4,\frac{1}{2} },J_3)$}
$$
\left\{
\begin{array}{l}
[{\it e1},{\it e2}]=-\frac{1}{2}\,{\it e2},
[{\it e1},{\it e3}]=-{\it e3},
[{\it e1},{\it e4}]=-\frac{1}{2}\,{\it e4},
[{\it e2},{\it e4}]={\it e3}
\end{array}
\right.
$$

It is non-unimodular and its (left)-invariant complex structure is non-abelian

\begin{enumerate}
\begin{multicols}{2}
  \item $\operatorname{Dim} \operatorname{Der} = { 7 }$
  \item $\operatorname{Dim} \operatorname{Der}_{J} = { 1 }$
\end{multicols}
\end{enumerate}

%%%%%%%%%%%%%%%%%%%%%%%%%%%%%%%%%%%%%%%%%%%%%%%%%%%%%%%%%%%%%%%%%%%%%%%%%%%%%%%%%%%%%%%%%%%%%%%%%%%%%%%%%%%%%%%%%%%%%%%%%%%%%%%%%%%%%%%%%%%%%%%%%%%%

\subsection{$(\mathfrak{ d}^{'}_{4,a},J_1)$ with $a>0$}
$$
\left\{
\begin{array}{l}
[{\it e1},{\it e2}]=a{\it e2}+{\it e4},
[{\it e1},{\it e3}]=2\,a{\it e3},
[{\it e1},{\it e4}]=-{\it e2}+a{\it e4},
[{\it e2},{\it e4}]={\it e3}
\end{array}
\right.
$$

It is non-unimodular and its (left)-invariant complex structure is non-abelian

\begin{enumerate}
\begin{multicols}{2}
  \item $\operatorname{Dim} \operatorname{Der} = { 5 }$
  \item $\operatorname{Dim} \operatorname{Der}_{J} = { 1 }$
\end{multicols}
\end{enumerate}

%%%%%%%%%%%%%%%%%%%%%%%%%%%%%%%%%%%%%%%%%%%%%%%%%%%%%%%%%%%%%%%%%%%%%%%%%%%%%%%%%%%%%%%%%%%%%%%%%%%%%%%%%%%%%%%%%%%%%%%%%%%%%%%%%%%%%%%%%%%%%%%%%%%%

\subsection{$(\mathfrak{ d}^{'}_{4,a},J_2)$ with $a>0$}
$$
\left\{
\begin{array}{l}
[{\it e1},{\it e2}]=-a{\it e2}+{\it e4},
[{\it e1},{\it e3}]=-2\,a{\it e3},
[{\it e1},{\it e4}]=-{\it e2}-a{\it e4},
[{\it e2},{\it e4}]={\it e3}
\end{array}
\right.
$$

It is non-unimodular and its (left)-invariant complex structure is non-abelian

\begin{enumerate}
\begin{multicols}{2}
  \item $\operatorname{Dim} \operatorname{Der} = { 5 }$
  \item $\operatorname{Dim} \operatorname{Der}_{J} = { 1 }$
\end{multicols}
\end{enumerate}

%%%%%%%%%%%%%%%%%%%%%%%%%%%%%%%%%%%%%%%%%%%%%%%%%%%%%%%%%%%%%%%%%%%%%%%%%%%%%%%%%%%%%%%%%%%%%%%%%%%%%%%%%%%%%%%%%%%%%%%%%%%%%%%%%%%%%%%%%%%%%%%%%%%%

\subsection{$(\mathfrak{ d}^{'}_{4,a},J_3)$ with $a>0$}
$$
\left\{
\begin{array}{l}
[{\it e1},{\it e2}]=a{\it e2}-{\it e4},
[{\it e1},{\it e3}]=2\,a{\it e3},
[{\it e1},{\it e4}]={\it e2}+a{\it e4},
[{\it e2},{\it e4}]={\it e3}
\end{array}
\right.
$$

It is non-unimodular and its (left)-invariant complex structure is non-abelian

\begin{enumerate}
\begin{multicols}{2}
  \item $\operatorname{Dim} \operatorname{Der} = { 5 }$
  \item $\operatorname{Dim} \operatorname{Der}_{J} = { 1 }$
\end{multicols}
\end{enumerate}

%%%%%%%%%%%%%%%%%%%%%%%%%%%%%%%%%%%%%%%%%%%%%%%%%%%%%%%%%%%%%%%%%%%%%%%%%%%%%%%%%%%%%%%%%%%%%%%%%%%%%%%%%%%%%%%%%%%%%%%%%%%%%%%%%%%%%%%%%%%%%%%%%%%%

%%%%%%%%%%%%%%%%%%%%%%%%%%%%%%%%%%%%%%%%%%%%%%%%%%%%%%%%%%%%%%%%%%%%%%%%%%%%%%%%%%%%%%%%%%%%%%%%%%%%%%%%%%%%%%%%%%%%%%%%%%%%%%%%%%%%%%%%%%%%%%%%%%%%

\subsection{$(\mathfrak{ d}^{'}_{4,a},J_4)$ with $a>0$}
$$
\left\{
\begin{array}{l}
[{\it e1},{\it e2}]=-a{\it e2}-{\it e4},
[{\it e1},{\it e3}]=-2\,a{\it e3},
[{\it e1},{\it e4}]={\it e2}-a{\it e4},
[{\it e2},{\it e4}]={\it e3}
\end{array}
\right.
$$

It is non-unimodular and its (left)-invariant complex structure is non-abelian

\begin{enumerate}
\begin{multicols}{2}
  \item $\operatorname{Dim} \operatorname{Der} = { 5 }$
  \item $\operatorname{Dim} \operatorname{Der}_{J} = { 1 }$
\end{multicols}
\end{enumerate}

\subsection{$(\mathfrak{ d}^{'}_{4,0 },J_1)$}
$$
\left\{
\begin{array}{l}
[{\it e1},{\it e2}]={\it e4},
[{\it e1},{\it e4}]=-{\it e2},
[{\it e2},{\it e4}]={\it e3}
\end{array}
\right.
$$

It is unimodular and its (left)-invariant complex structure is non-abelian

\begin{enumerate}
\begin{multicols}{2}
  \item $\operatorname{Dim} \operatorname{Der} = { 5 }$
  \item $\operatorname{Dim} \operatorname{Der}_{J} = { 1 }$
\end{multicols}
\end{enumerate}

%%%%%%%%%%%%%%%%%%%%%%%%%%%%%%%%%%%%%%%%%%%%%%%%%%%%%%%%%%%%%%%%%%%%%%%%%%%%%%%%%%%%%%%%%%%%%%%%%%%%%%%%%%%%%%%%%%%%%%%%%%%%%%%%%%%%%%%%%%%%%%%%%%%%

\subsection{$(\mathfrak{ d}^{'}_{4,0 },J_2)$}
$$
\left\{
\begin{array}{l}
[{\it e1},{\it e2}]=-{\it e4},
[{\it e1},{\it e4}]={\it e2},
[{\it e2},{\it e4}]={\it e3}
\end{array}
\right.
$$

It is unimodular and its (left)-invariant complex structure is non-abelian

\begin{enumerate}
\begin{multicols}{2}
  \item $\operatorname{Dim} \operatorname{Der} = { 5 }$
  \item $\operatorname{Dim} \operatorname{Der}_{J} = { 1 }$
\end{multicols}
\end{enumerate}

%%%%%%%%%%%%%%%%%%%%%%%%%%%%%%%%%%%%%%%%%%%%%%%%%%%%%%%%%%%%%%%%%%%%%%%%%%%%%%%%%%%%%%%%%%%%%%%%%%%%%%%%%%%%%%%%%%%%%%%%%%%%%%%%%%%%%%%%%%%%%%%%%%%%

\subsection{$(\mathfrak{ h}_{ 4},J)$}
$$
\left\{
\begin{array}{l}
[{\it e1},{\it e2}]={\it e2},
[{\it e1},{\it e3}]={\it e2}+{\it e3},
[{\it e1},{\it e4}]=2\,{\it e4},
[{\it e2},{\it e3}]={\it e4}
\end{array}
\right.
$$

It is non-unimodular and its (left)-invariant complex structure is non-abelian

\begin{enumerate}
\begin{multicols}{2}
  \item $\operatorname{Dim} \operatorname{Der} = { 5 }$
  \item $\operatorname{Dim} \operatorname{Der}_{J} = { 1 }$
\end{multicols}
\end{enumerate}

\section{Unimodular case}

\begin{table}[ht]
\centering
%\operatorname{Dim}
\begin{center}
    \begin{tabular}{c : c: p{6cm} : p{5cm} }
     $\operatorname{Dim} \operatorname{Der}_{J_{\tiny{\mbox{cn}}}}$ & $\operatorname{Dim} \operatorname{Der}$ & \multicolumn{2}{c}{Unimodular Lie Algebra + Complex structure}\\
    \hdashline
     \multirow{2}{*}{$1$} & 4  &  $(\mathfrak{sl}(2,\mathbb{R})\times \mathbb{R}, J(t))$ , $(\mathfrak{so}(3,\mathbb{R})\times \mathbb{R}, J(t))$ & \\
                          & 5  &  & $(\mathfrak{d}_4,J_2)$, $(\mathfrak{d}^{'}_{4,0},J_1)$, $(\mathfrak{d}^{'}_{4,0},J_2)$\\
   \hdashline
     \multirow{3}{*}{$2$} & 5  &   & $(\mathfrak{d}_4,J_1)$\\
                          & 6  &  $(\mathfrak{r}^{'}_{4,-1/2,\delta},J_1)$,$(\mathfrak{r}^{'}_{4,-1/2,\delta},J_2)$, &  $(\mathfrak{rr}^{'}_{3,0}\times \mathbb{R},J)$\\
                          & 8  &  $(\mathfrak{r}_{4,-1/2,-1/2},J)$ & \\
  \hdashline
   $4$ & $10$ & $(\mathfrak{rh}_{3}\times \mathbb{R},J)$ & \\
   \hdashline
   $8$ & $16$ & $(\mathfrak{a}_{4},J)$ & \\
    \end{tabular}
\end{center}
{\scriptsize Table 1: Dimension of  complex derivations.}
\end{table}

\begin{center}
\begin{tikzpicture}
\node[circle,fill,inner sep=0pt,minimum size=3pt,label=above:{\scriptsize$(\mathfrak{d}^{'}_{4,0},J_1)$}] (d4p1) at (-1,+0) {};
\node[circle,fill,inner sep=0pt,minimum size=3pt,label=above:{\scriptsize$(\mathfrak{d}^{'}_{4,0},J_2)$}] (d4p2) at (+1,+0) {};

%\node[circle,fill,inner sep=0pt,minimum size=3pt,label=above:{}] (so3a) at (-2,+0) {};
\node[circle,fill,inner sep=0pt,minimum size=3pt,label=above:{\begin{tabular}{c} \scriptsize$(\mathfrak{so}(3,\mathbb{R})\times \mathbb{R}, J(t))$ \\ \scriptsize$t \in \mathbb{R}$ \end{tabular}}] (so3) at (-3,+0) {};
%\node[circle,fill,inner sep=0pt,minimum size=3pt,label=left:{}]  (so3b) at (-4,+0) {};

%\node[circle,fill,inner sep=0pt,minimum size=3pt,label=right:{}] (sl2a) at (+2,+0) {};
\node[circle,fill,inner sep=0pt,minimum size=3pt,label=above:{\begin{tabular}{c} \scriptsize$(\mathfrak{sl}(2,\mathbb{R})\times \mathbb{R}, J(t))$ \\ \scriptsize$t \in \mathbb{R}\setminus\{0\}$ \end{tabular}} ] (sl2)  at (+3,+0) {};
%\node[circle,fill,inner sep=0pt,minimum size=3pt,label=left:{ }] (sl2b) at (+4,+0) {};

\node[circle,fill,inner sep=0pt,minimum size=3pt,label=above:{\scriptsize$(\mathfrak{sl}(2,\mathbb{R})\times \mathbb{R}, J(0))$}](sl20)  at (+6,+0) {};

\node[circle,fill,inner sep=0pt,minimum size=3pt,label=above:{\scriptsize$(\mathfrak{d}_4,J_2)$ }](d42)  at (+8,+0) {};

\node[circle,fill,inner sep=0pt,minimum size=3pt,label=below:{\scriptsize$(\mathfrak{rr}^{'}_{3,0}\times \mathbb{R},J)$}] (r3p)  at (+0,-2) {};

%\node[circle,fill,inner sep=0pt,minimum size=3pt,label=right:{}] (r4p2a) at (-2,-2) {};
\node[circle,fill,inner sep=0pt,minimum size=3pt,label=above:{\begin{tabular}{c} \scriptsize$(\mathfrak{r}^{'}_{4,-\tfrac{1}{2},\delta},J_1)$ \\ \scriptsize$ \delta >0$ \end{tabular}}] (r4p1) at (-3,-2) {};
%\node[circle,fill,inner sep=0pt,minimum size=3pt,label=above:{}] (r4p2b) at (-4,-2) {};

%\node[circle,fill,inner sep=0pt,minimum size=3pt,label=left:{}]  (r4p1a) at (+2,-2) {};
\node[circle,fill,inner sep=0pt,minimum size=3pt,label=right:{\begin{tabular}{c} \scriptsize$(\mathfrak{r}^{'}_{4,-\tfrac{1}{2},\delta},J_2)$ \\ \scriptsize$ \delta >0$ \end{tabular}}]  (r4p2) at (+3,-2) {};
%\node[circle,fill,inner sep=0pt,minimum size=3pt,label=right:{}] (r4p1b) at (+4,-2) {};

\node[circle,fill,inner sep=0pt,minimum size=3pt,label=above:{\scriptsize$(\mathfrak{r}_{4,-\tfrac{1}{2},-\tfrac{1}{2}},J)$ }] (r4) at (-6,-2) {};
\node[circle,fill,inner sep=0pt,minimum size=3pt,label=below:{\scriptsize$(\mathfrak{d}_4,J_1)$ }] (d41)at (+7,-2) {};
\node[circle,fill,inner sep=0pt,minimum size=3pt,label=right:{\scriptsize$(\mathfrak{rh}_{3}\times \mathbb{R},J)$ }] (r3) at (+0,-4) {};
\node[circle,fill,inner sep=0pt,minimum size=3pt,label=below:{\scriptsize$(\mathfrak{a}_4,J)$ }] (a4) at (+0,-5) {};
%\draw[thick] (so3a)edge (so3b)
%\draw[thick] (sl2a) edge (sl2) (sl2) edge (sl2b);
%\draw[thick] (r4p2a) edge (r4p2b) (r4p1a) edge (r4p1b);
\draw [->] (d4p1) edge (r3p)  (d4p2) edge (r3p) (so3) edge (r3p); %(so3a) edge (r3p) (so3b) edge (r3p)
%\draw [->] (sl2a) edge (r3p) (sl2b) edge (r3p);
\draw [->] (sl2)edge (r3p)  (sl20)edge (r3p) ;
%\draw [->] (r4p2a)edge (r3)(r4p2b)edge (r3) (r4p1a)edge (r3)(r4p1b)edge (r3);
\draw [->] (r4p2)edge (r3) (r4p1)edge (r3);
\draw [->] (d41) edge (r3);
\draw [->] (r3p) edge (r3) (r4) edge (r3);
\draw [->] (d42) edge (d41) (sl20)edge (d41);
\draw [->] (r3)  edge (a4);

\end{tikzpicture}
\end{center}

\subsection{$ (\mathfrak{sl}(2,\mathbb{R})\times \mathbb{R}, J( t)) $ with $t \in \mathbb{R}$}\label{sl2}

\begin{enumerate}

\item  For any $t\in \mathbb{R} \setminus\{0\}$,
$
{(\mathfrak{sl}(2,\mathbb{R})\times \mathbb{R}, J( t))}
\centernot{ \xrightarrow{\text{\,\,deg\,\,}} }
{(\mathfrak{d}_4,J_1)}
$\newline

Suppose for the sake of contradiction that $
{(\mathfrak{sl}(2,\mathbb{R})\times \mathbb{R}, J( t))}
%\centernot
{ \xrightarrow{\text{\,\,deg\,\,}} }
{(\mathfrak{d}_4,J_2)}
$. First, we consider the $\operatorname{GL}(\mathbb{R}^4, {J_{\tiny{\mbox{cn}}}})$-equivariant continuous function $\psi_{\frac{1}{t},0}: C^{2}_{J_{\tiny{\mbox{cn}}}}(\mathbb{R}^{4};\mathbb{R}^{4}) \rightarrow C^{2}(\mathbb{R}^{4};\mathbb{R}^{4})$. Therefore $$\mathfrak{A}:=(\mathbb{R}^{4}, \psi_{\frac{1}{t},0}(\mu_{8,-}(t)) { \xrightarrow{\text{\,\,deg\,\,}} } \mathfrak{B}:=(\mathbb{R}^{4}, \psi_{\frac{1}{t},0}( \mu_{12} )),$$
where:
$$
\mathfrak{A} =
\left\{
\begin{array}{l}
[{\it e1},{\it e2}]=-{\frac {1}{t}}{\it e1}+{\it e3},
[{\it e1},{\it e3}]={\frac { \left( {t}^{2}+1 \right)}{t}} {\it e4},
[{\it e2},{\it e3}]={\it e1}+{\frac {1}{t}}{\it e3}
\end{array}
\right.
$$
and
$$
\mathfrak{B} =
\left\{
\begin{array}{l}
[{\it e1},{\it e2}]=-{\it e2}-{\frac {1}{t}}{\it e4},
[{\it e1},{\it e3}]=-{\frac {1}{t}}{\it e1}+{\it e3},
[{\it e2},{\it e3}]=-{\frac {1}{t}}{\it e2}+{\it e4}
\end{array}
\right.
$$
Consider the $\operatorname{GL}(\mathbb{R}^4)$-equivariant continuous function $\kappa: C^{2}(\mathbb{R}^{4};\mathbb{R}^{4}) \rightarrow \Sigma^2(\mathbb{R}^{4}; \mathbb{R})$ defined by
$$
\kappa(\lambda)(X,Y):= \operatorname{Trace}(\lambda(X, \lambda(Y , \square)));
$$
the symmetric bilinear form $\kappa(\lambda)$ is just the Killing form of $(\mathbb{R}^4, \lambda)$. Therefore $\kappa(\mathfrak{A}) {\xrightarrow{\text{\,\,deg\,\,}} } \kappa(\mathfrak{B}) $. The matrix of $\kappa(\mathfrak{A})$ with respect to the canonical basis is:
$$
\operatorname{diag}(0,\frac{2}{t^2}-2,0,0)
$$
while the matrix of $\kappa(\mathfrak{B})$ with respect to the canonical basis is:
$$
\operatorname{diag}(2,0,\frac{2}{t^2},0)
$$
which is a contradiction because the dimension of the radical of $\kappa(\mathfrak{A})$ is $3$ or $4$ and the dimension of the radical of $\kappa(\mathfrak{B})$ is $2$.

\item\label{sl2d4}
$
{(\mathfrak{sl}(2,\mathbb{R})\times \mathbb{R}, J( 0 ))}
\xrightarrow{\text{\,\,deg\,\,}}
{(\mathfrak{d}_4,J_1)}
$\newline

Consider the curve $g_t = \left(
          \begin{array}{cc}
            A(t) & -B(t) \\
            B(t) & A(t) \\
          \end{array}
        \right)$ where $$A(t)= \left(
                                 \begin{array}{cc}
                                   1 &0 \\
                                   0 & e^t \\
                                 \end{array}
                               \right)
                             \qquad
                             \mbox{and}
                             \qquad
                         B(t)=\left(
                                 \begin{array}{cc}
                                   0 &0 \\
                                   e^t & 0 \\
                                 \end{array}
                                  \right).$$
Note that $g_t \in \operatorname{GL}(\mathbb{R}^4, {J_{\tiny{\mbox{cn}}}})$ and its action on $g_t \cdot \mu_{8,-}(0)$ gives
$$
\left\{
\begin{array}{l}
[{\it e1},{\it e2}]=-{\it e2}+{{\rm e}^{-t}}{\it e3},
[{\it e1},{\it e3}]={\it e3},
[{\it e2},{\it e3}]={{\rm e}^{-t}}{\it e1}+{\it e4}
\end{array}
\right.
$$
Taking the limit $g_t \cdot \mu_{8,-}(0)$ as $t \to \infty$, we obtain $\mu_{10}(0)$.

\item\label{sl2r4prima} For any $t\in \mathbb{R}$,
$
{(\mathfrak{sl}(2,\mathbb{R})\times \mathbb{R}, J(t))}
\centernot{ \xrightarrow{\text{\,\,deg\,\,}} }
{(\mathfrak{r}^{'}_{4,-1/2,\delta}, J_k)}
$ with $\delta>0$ and $k=1,2$. \newline

In fact, this holds because $\mathfrak{sl}(2,\mathbb{R})\times \mathbb{R} \centernot{ \xrightarrow{\text{\,\,deg\,\,}} } \mathfrak{r}^{'}_{4,-1/2,\delta} $.  Let us now explore this in detail. Suppose for the sake of contradiction that $\mathfrak{sl}(2,\mathbb{R}) \times \mathbb{R} \xrightarrow{\text{\,\,deg\,\,}} \mathfrak{r}^{'}_{4,-1/2,\delta} $. First, we consider the $\operatorname{GL}(\mathbb{R}^4)$-equivariant continuous function $\Psi: C^{2}(\mathbb{R}^{4};\mathbb{R}^{4}) \rightarrow \sum^{3}(\mathbb{R}^{4};\mathbb{R})$ given by $$\Psi(\mu)(X_1,X_2,X_3)=\sum_{\sigma \in \mathfrak{S}_{3}}\operatorname{Trace}(\mu(X_{\sigma(1)}, \mu(X_{\sigma(2)}, \mu(X_{\sigma(3)}, \square )))).$$ Therefore $$ \Psi(\mu_{8,-}(a)) \xrightarrow{\text{\,\,deg\,\,}} \Psi(\mu_{6}(-\frac{1}{2},\delta)), $$
where:
$$
 \Psi(\mu_{8,-}(a)) \equiv 0
$$
and
$$
\Psi(\mu_{6}(-\frac{1}{2},\delta))= \frac{9}{2}(1+4\delta^2)e_1^{\ast}\otimes e_1^{\ast}\otimes e_1^{\ast}
$$
which is a contradiction.

We can give another proof by using $(\alpha,\beta,\gamma)$-derivations since $\operatorname{Dim}\operatorname{Der}(-1,1,1)(\mathfrak{sl}(2,\mathbb{R})\times \mathbb{R})=6$ and $\operatorname{Dim}\operatorname{Der}(-1,1,1)(\mathfrak{r}^{'}_{4,-1/2,\delta})=4$.

\item\label{sl2r4menosunmedio} For any $t\in \mathbb{R}$,
$
{(\mathfrak{sl}(2,\mathbb{R})\times \mathbb{R}, J(t))}
\centernot{ \xrightarrow{\text{\,\,deg\,\,}} }
{(\mathfrak{r}_{4,-1/2,-1/2}, J)}
$\newline

In fact, this holds because $
{\mathfrak{sl}(2,\mathbb{R})\times \mathbb{R}}
\centernot{ \xrightarrow{\text{\,\,deg\,\,}} }
{\mathfrak{r}_{4,-1/2,-1/2}}
$.  Let us now explore this in detail. Suppose for the sake of contradiction that $
{\mathfrak{sl}(2,\mathbb{R})\times \mathbb{R}}
{ \xrightarrow{\text{\,\,deg\,\,}} }
{\mathfrak{r}_{4,-1/2,-1/2}}
$. First, we consider the $\operatorname{GL}(\mathbb{R}^4)$-equivariant continuous function $\Psi: C^{2}(\mathbb{R}^{4};\mathbb{R}^{4}) \rightarrow L^{3}(\mathbb{R}^{4};\mathbb{R})$ given by $$\Psi(\mu)(X_1,X_2,X_3)=\sum_{\sigma \in \mathfrak{S}_{3}}\operatorname{Trace}(\mu(X_{\sigma(1)}, \mu(X_{\sigma(2)}, \mu(X_{\sigma(3)}, \square )))).$$ Therefore $$ \Psi( \mu_{8,-}(a) ) \xrightarrow{\text{\,\,deg\,\,}}
\Psi( \mu_{6}(-\frac{1}{2},0) ), $$
where:
$$
 \Psi(  \mu_{8,-}(a) ) \equiv 0
$$
and
$$
\Psi( \mu_{6}(-\frac{1}{2},0) ) = \frac{9}{2}e^{\ast}_{1}\otimes e^{\ast}_{1}\otimes e^{\ast}_{1}
$$
which is a contradiction.

We can give another proof by using $(\alpha,\beta,\gamma)$-derivations since $\operatorname{Dim}\operatorname{Der}(-1,1,1)(\mathfrak{sl}(2,\mathbb{R})\times \mathbb{R})=6$ and $\operatorname{Dim}\operatorname{Der}(-1,1,1)(\mathfrak{r}_{4,-1/2,-1/2})=4$.

\item For any $t\in \mathbb{R}$,
$
{(\mathfrak{sl}(2,\mathbb{R})\times \mathbb{R}, J(t))}
\xrightarrow{\text{\,\,deg\,\,}}
{(\mathfrak{rr}^{'}_{3,0}\times \mathbb{R},J)}
$\newline

Consider the curve $g_u = \left(
          \begin{array}{cc}
            A(u) & -B(u) \\
            B(u) & A(u) \\
          \end{array}
        \right)$ where $$A(u)=\left(
                                \begin{array}{cc}
                                  0 & -1 \\
                                  0 & 0 \\
                                \end{array}
                              \right)
                             \qquad
                             \mbox{and}
                             \qquad
                         B(u)=\left(
                                \begin{array}{cc}
                                  0 & 0 \\
                                  e^u & 0 \\
                                \end{array}
                              \right)
                         .$$
Note that $g_u \in \operatorname{GL}(\mathbb{R}^4, {J_{\tiny{\mbox{cn}}}})$ and its action on $\mu_{8,-}(t)$ gives
$$
\left\{
\begin{array}{l}
[{\it e1},{\it e2}]={\it e4},
[{\it e1},{\it e4}]=-{\it e2},
[{\it e2},{\it e4}]=-{{\rm e}^{-2\,u}}{\it e1}-{{\rm e}^{-2\,u}}t{\it e3}
\end{array}
\right.
$$
Taking the limit $g_u \cdot \mu_{8,-}(t)$ as $u \to \infty$, we obtain $\mu_{7}(0)$.

\end{enumerate}

\subsection{$ {(\mathfrak{so}(3,\mathbb{R})\times \mathbb{R}, J(t))} $ with $t \in \mathbb{R}$}

\begin{enumerate}

\item  For any $t\in \mathbb{R}$,
$
{(\mathfrak{so}(3,\mathbb{R})\times \mathbb{R}, J(t))}
\centernot{ \xrightarrow{\text{\,\,deg\,\,}} }
{(\mathfrak{d}_4,J_1)}
$\newline

In fact, this holds because $
{\mathfrak{so}(3,\mathbb{R})\times \mathbb{R}}
\centernot{ \xrightarrow{\text{\,\,deg\,\,}} }
{\mathfrak{d}_4}
$.  Let us now explore this in detail. Suppose for the sake of contradiction that $
{\mathfrak{so}(3,\mathbb{R})\times \mathbb{R}}
{ \xrightarrow{\text{\,\,deg\,\,}} }
{\mathfrak{d}_4}
$. First, consider the \(\operatorname{GL}(\mathbb{R}^4)\)-equivariant continuous function \[ \kappa_{g}: C^{2}(\mathbb{R}^{4};\mathbb{R}^{4}) \rightarrow \Sigma^2(\mathbb{R}^{4}; \mathbb{R}). \]
$$
\kappa(\lambda)(X,Y):= \operatorname{Trace}(\lambda(X, \lambda(Y , \square)));
$$
the symmetric bilinear form $\kappa(\lambda)$ is just the Killing form of $(\mathbb{R}^4, \lambda)$. Therefore $\kappa( \mu_{8,+}(a)) {\xrightarrow{\text{\,\,deg\,\,}} } \kappa( \mu_{12} ) $
where:
$$
\kappa( \mathfrak{so}(3,\mathbb{R})\times \mathbb{R} ) = -2( e_1^\ast \otimes e_1^\ast + e_2^\ast \otimes e_2^\ast + e_3^\ast \otimes e_3^\ast)
$$
and
$$
\kappa( \mathfrak{d}_4 ) = 2 e_1^\ast \otimes e_1^\ast
$$
which is a contradiction because $\Psi( \mathfrak{so}(3,\mathbb{R})\times \mathbb{R} ) $ is negative semi-definite bilinear form while $\Psi( \mathfrak{d}_4 ) $ is non-zero positive semi-definite bilinear form.

\item  For any $t\in \mathbb{R}$,
$
{(\mathfrak{so}(3,\mathbb{R})\times \mathbb{R}, J(t))}
\centernot{ \xrightarrow{\text{\,\,deg\,\,}} }
{(\mathfrak{r}^{'}_{4,-1/2,\delta}, J_k)}
$ with $\delta>0$ and $k=1,2$ \newline

In fact, this holds because $
{\mathfrak{so}(3,\mathbb{R})\times \mathbb{R}}
\centernot{ \xrightarrow{\text{\,\,deg\,\,}} }
{\mathfrak{r}^{'}_{4,-1/2,\delta}}
$.  Let us now explore this in detail. Suppose for the sake of contradiction that $
{\mathfrak{so}(3,\mathbb{R})\times \mathbb{R}}
{ \xrightarrow{\text{\,\,deg\,\,}} }
{\mathfrak{r}^{'}_{4,-1/2,\delta}}
$. First, we consider the $\operatorname{GL}(\mathbb{R}^4)$-equivariant continuous function $\Psi: C^{2}(\mathbb{R}^{4};\mathbb{R}^{4}) \rightarrow L^{3}(\mathbb{R}^{4};\mathbb{R})$ given by $$\Psi(\mu)(X_1,X_2,X_3)=\sum_{\sigma \in \mathfrak{S}_{3}}\operatorname{Trace}(\mu(X_{\sigma(1)}, \mu(X_{\sigma(2)}, \mu(X_{\sigma(3)}, \square )))).$$ Therefore $$ \Psi( \mu_{8,+}(t) ) \xrightarrow{\text{\,\,deg\,\,}}
\Psi( \mu_{6}(-\frac{1}{2},\delta) ), $$
where:
$$
 \Psi( \mu_{8,+}(t) ) \equiv 0
$$
and
$$
\Psi( \mu_{6}(-\frac{1}{2},\delta) )= \frac{9}{2}(1+4\delta^2) e_1^\ast \otimes e_1^\ast \otimes e_1^\ast
$$
which is a contradiction.

We can give another proof by using $(\alpha,\beta,\gamma)$-derivations since $\operatorname{Dim}\operatorname{Der}(-1,1,1)(\mathfrak{so}(3,\mathbb{R})\times \mathbb{R})=6$ and $\operatorname{Dim}\operatorname{Der}(-1,1,1)({\mathfrak{r}^{'}_{4,-1/2,\delta}})= 4 $.

\item For any $t\in \mathbb{R}$,
$
{(\mathfrak{so}(3,\mathbb{R})\times \mathbb{R}, J(t))}
\centernot{ \xrightarrow{\text{\,\,deg\,\,}} }
{(\mathfrak{r}_{4,-1/2,-1/2}, J)}
$\newline

In fact, this holds because $
{\mathfrak{so}(3,\mathbb{R})\times \mathbb{R}}
\centernot{ \xrightarrow{\text{\,\,deg\,\,}} }
{\mathfrak{r}_{4,-1/2,-1/2}}
$.  Let us now explore this in detail. Suppose for the sake of contradiction that $
{\mathfrak{so}(3,\mathbb{R})\times \mathbb{R}}
{ \xrightarrow{\text{\,\,deg\,\,}} }
{\mathfrak{r}_{4,-1/2,-1/2}}
$. First, consider the \(\operatorname{GL}(\mathbb{R}^4)\)-equivariant continuous function \[ \kappa_{g}: C^{2}(\mathbb{R}^{4};\mathbb{R}^{4}) \rightarrow \Sigma^2(\mathbb{R}^{4}; \mathbb{R}). \]
$$
\kappa(\lambda)(X,Y):= \operatorname{Trace}(\lambda(X, \lambda(Y , \square)));
$$
the symmetric bilinear form $\kappa(\lambda)$ is just the Killing form of $(\mathbb{R}^4, \lambda)$.

Therefore $\kappa( \mu_{8,+}(t) ) {\xrightarrow{\text{\,\,deg\,\,}} } \kappa( \mu_{6}(-\frac{1}{2},0 ) ) $
where:
$$
\kappa( \mu_{8,+}(t) ) = -2( e_1^\ast \otimes e_1^\ast + e_2^\ast \otimes e_2^\ast + e_3^\ast \otimes e_3^\ast)
$$
and
$$
\kappa( \mu_{6}(-\frac{1}{2},0) ) = \frac{3}{2}  e_1^\ast \otimes e_1^\ast
$$
which is a contradiction because $\kappa( \mathfrak{so}(3,\mathbb{R})\times \mathbb{R} ) $ is negative semi-definite bilinear form while $\kappa( \mathfrak{r}_{4,-1/2,-1/2}  ) $ is non-zero positive semi-definite bilinear form.

We can give another proof by using $(\alpha,\beta,\gamma)$-derivations since $\operatorname{Dim}\operatorname{Der}(-1,1,1)(\mathfrak{so}(3,\mathbb{R})\times \mathbb{R})=6$ and $\operatorname{Dim}\operatorname{Der}(-1,1,1)(\mathfrak{r}_{4,-1/2,-1/2})=4$.

\item For any $t\in \mathbb{R}$,
$
{(\mathfrak{so}(3,\mathbb{R})\times \mathbb{R}, J(t))}
\xrightarrow{\text{\,\,deg\,\,}}
{(\mathfrak{rr}^{'}_{3,0}\times \mathbb{R},J)}
$\newline

Consider the curve $g_u = \left(
          \begin{array}{cc}
            A(u) & -B(u) \\
            B(u) & A(u) \\
          \end{array}
        \right)$ where $$A(t)=\left(
                                \begin{array}{cc}
                                  0 & -1 \\
                                  0 & 0 \\
                                \end{array}
                              \right)
                             \qquad
                             \mbox{and}
                             \qquad
                         B(t)=\left(
                                \begin{array}{cc}
                                  0 & 0 \\
                                  e^u & 0 \\
                                \end{array}
                              \right)
                         .$$
Note that $g_u \in \operatorname{GL}(\mathbb{R}^4, {J_{\tiny{\mbox{cn}}}})$ and its action on $\mu_{8,+}(a)$ gives
$$
\left\{
\begin{array}{l}
[{\it e1},{\it e2}]={\it e4},
[{\it e1},{\it e4}]=-{\it e2},
[{\it e2},{\it e4}]={{\rm e}^{-2\,t}}{\it e1}-{{\rm e}^{-2\,t}}a{\it e3}
\end{array}
\right.
$$
Taking the limit $g_u \cdot \mu_{8,+}(a)$ as $t \to \infty$, we obtain $\mu_{7}(0)$.

\end{enumerate}

\subsection{$(\mathfrak{d}_4,J_2)$}

\begin{enumerate}

\item
$
(\mathfrak{d}_4,J_2)
\xrightarrow{\text{\,\,deg\,\,}}
{(\mathfrak{d}_4,J_1)}
$\newline

Consider the curve $g_t = \left(
          \begin{array}{cc}
            A(t) & -B(t) \\
            B(t) & A(t) \\
          \end{array}
        \right)$ where $$A(t)= \operatorname{diag}(1,e^{-t})
                             \qquad
                             \mbox{and}
                             \qquad
                         B(t)=0
                         .$$
Note that $g_t \in \operatorname{GL}(\mathbb{R}^4, {J_{\tiny{\mbox{cn}}}})$ and its action on $\mu_{12}$ gives
$$
\left\{
\begin{array}{l}
[{\it e1},{\it e2}]=-{\it e2},
[{\it e1},{\it e3}]={\it e3}+{{\rm e}^{-t}}{\it e4},
[{\it e2},{\it e3}]={\it e4}
\end{array}
\right.
$$
Taking the limit $g_t \cdot \mu_{12}$ as $t \to \infty$, we obtain $\mu_{10}(0)$.

\item
$
(\mathfrak{d}_4,J_2)
\centernot{ \xrightarrow{\text{\,\,deg\,\,}} }
{(\mathfrak{rr}^{'}_{3,0}\times \mathbb{R},J)}
$\newline

In fact, this holds because $
\mathfrak{d}_4
\centernot{ \xrightarrow{\text{\,\,deg\,\,}} }
{\mathfrak{rr}^{'}_{3,0}\times \mathbb{R}}
$.  Let us now explore this in detail. Suppose for the sake of contradiction that $
\mathfrak{d}_4
\centernot{ \xrightarrow{\text{\,\,deg\,\,}} }
{\mathfrak{rr}^{'}_{3,0}\times \mathbb{R}}
$. First, consider the \(\operatorname{GL}(\mathbb{R}^4)\)-equivariant continuous function \[ \kappa_{g}: C^{2}(\mathbb{R}^{4};\mathbb{R}^{4}) \rightarrow \Sigma^2(\mathbb{R}^{4}; \mathbb{R}). \]
$$
\kappa(\lambda)(X,Y):= \operatorname{Trace}(\lambda(X, \lambda(Y , \square)));
$$
the symmetric bilinear form $\kappa(\lambda)$ is just the Killing form of $(\mathbb{R}^4, \lambda)$. Therefore $\kappa( \mu_{12} ) {\xrightarrow{\text{\,\,deg\,\,}} } \kappa( \mu_{7}(0) ) $
where:
$$
\kappa( \mu_{12} ) = 2 e_1^{\ast}\otimes e_{1}^{\ast}
$$
and
$$
\kappa( \mu_{7}(0) ) = -2 e_1^{\ast}\otimes e_{1}^{\ast}
$$
which is a contradiction because $\Psi( \mathfrak{d}_4  ) $ is positive semi-definite bilinear form while $\kappa( {\mathfrak{rr}^{'}_{3,0}\times \mathbb{R}} ) $ is non-zero negative semi-definite bilinear form.

\item
$
(\mathfrak{d}_4,J_2)
\centernot{ \xrightarrow{\text{\,\,deg\,\,}} }
{(\mathfrak{r}^{'}_{4,-1/2,\delta}, J_k)}
$ with $\delta>0$ and $k=1,2$ \newline
See Subsection \ref{sl2} items \ref{sl2d4} and \ref{sl2r4prima}.

\item
$
(\mathfrak{d}_4,J_2)
\centernot{ \xrightarrow{\text{\,\,deg\,\,}} }
{(\mathfrak{r}_{4,-1/2,-1/2}, J)}
$\newline
See Subsection \ref{sl2} items \ref{sl2d4} and \ref{sl2r4menosunmedio}.

\end{enumerate}

\subsection{$(\mathfrak{d}^{'}_{4,0},J_{k})$ with $k=1,2$}

\begin{enumerate}

\item
$
(\mathfrak{d}^{'}_{4,0},J_{k})
\centernot{ \xrightarrow{\text{\,\,deg\,\,}} }
{(\mathfrak{d}_4,J_1)}
$\newline

Because $
\mathfrak{d}^{'}_{4,0}
\centernot{ \xrightarrow{\text{\,\,deg\,\,}} }
{\mathfrak{d}_4}
$, since $\operatorname{Dim}\operatorname{Der}(\mathfrak{d}^{'}_{4,0})=5$ and $\operatorname{Dim}\operatorname{Der}({\mathfrak{d}_4})=5$.

\item
$
(\mathfrak{d}^{'}_{4,0},J_{k})
\xrightarrow{\text{\,\,deg\,\,}}
{(\mathfrak{rr}^{'}_{3,0}\times \mathbb{R},J)}
$\newline

Consider the curve $g_t = \left(
          \begin{array}{cc}
            A(t) & -B(t) \\
            B(t) & A(t) \\
          \end{array}
        \right)$ where $$A(t)= \operatorname{diag}((-1)^{k+1}, e^t)
                             \qquad
                             \mbox{and}
                             \qquad
                         B(t)=0
                         .$$
Note that $g_t \in \operatorname{GL}(\mathbb{R}^4, {J_{\tiny{\mbox{cn}}}})$ and its action on $\mu_{13,\pm}(0)$ gives
$$
\left\{
\begin{array}{l}
[{\it e1},{\it e2}]={\it e4},
[{\it e1},{\it e4}]=-{\it e2},
[{\it e2},{\it e4}]={{\rm e}^{-2\,t}}(-1)^{k+1}{\it e3}
\end{array}
\right.
$$
Taking the limit $g_t \cdot \mu_{13,\pm}(0) $ as $t \to \infty$, we obtain $\mu_{7}(0)$.

\item
$
(\mathfrak{d}^{'}_{4,0},J_{k_{2}})
\centernot{ \xrightarrow{\text{\,\,deg\,\,}} }
{(\mathfrak{r}^{'}_{4,-1/2,\delta}, J_k)}
$ with $\delta>0$ and $k=1,2$ \newline

In fact, this holds because $
\mathfrak{d}^{'}_{4,0}
\centernot{ \xrightarrow{\text{\,\,deg\,\,}} }
{\mathfrak{r}^{'}_{4,-1/2,\delta}}
$.  Let us now explore this in detail. Suppose for the sake of contradiction that $
\mathfrak{d}^{'}_{4,0}
{ \xrightarrow{\text{\,\,deg\,\,}} }
{\mathfrak{r}^{'}_{4,-1/2,\delta}}
$. First, we consider the $\operatorname{GL}(\mathbb{R}^4)$-equivariant continuous function $\Psi: C^{2}(\mathbb{R}^{4};\mathbb{R}^{4}) \rightarrow \Sigma^{3}(\mathbb{R}^{4};\mathbb{R})$ given by $$\Psi(\mu)(X_1,X_2,X_3)=\sum_{\sigma \in \mathfrak{S}_{3}}\operatorname{Trace}(\mu(X_{\sigma(1)}, \mu(X_{\sigma(2)}, \mu(X_{\sigma(3)}, \square )))).$$ Therefore $$ \Psi( \mu_{13,+}(0)) \xrightarrow{\text{\,\,deg\,\,}}
\Psi( \mu_{6}(-\frac{1}{2},\delta) ), $$
where:
$$
 \Psi( \mu_{13,+}(0) ) \equiv 0
$$
and
$$
\Psi( \mu_{6}(-\frac{1}{2},\delta) )= \frac{9}{2}(1+4\delta^2)e_{1}^\ast \otimes e_{1}^\ast \otimes e_{1}^\ast
$$
which is a contradiction.

We can give another proof by using $(\alpha,\beta,\gamma)$-derivations since $\operatorname{Dim}\operatorname{Der}(-1,1,1)( \mathfrak{d}^{'}_{4,0} )= 6$ and $\operatorname{Dim}\operatorname{Der}(-1,1,1)({\mathfrak{r}^{'}_{4,-1/2,\delta}})=4$.

\item
$
(\mathfrak{d}^{'}_{4,0},J_{k_{2}})
\centernot{ \xrightarrow{\text{\,\,deg\,\,}} }
{(\mathfrak{r}_{4,-1/2,-1/2}, J)}
$\newline

In fact, this holds because $
\mathfrak{d}^{'}_{4,0}
\centernot{ \xrightarrow{\text{\,\,deg\,\,}} }
{\mathfrak{r}_{4,-1/2,-1/2}}
$.  Let us now explore this in detail. Suppose for the sake of contradiction that $
\mathfrak{d}^{'}_{4,0}
{ \xrightarrow{\text{\,\,deg\,\,}} }
{\mathfrak{r}_{4,-1/2,-1/2}}
$. First, consider the \(\operatorname{GL}(\mathbb{R}^4)\)-equivariant continuous function \[ \kappa_{g}: C^{2}(\mathbb{R}^{4};\mathbb{R}^{4}) \rightarrow \Sigma^2(\mathbb{R}^{4}; \mathbb{R}). \]
$$
\kappa(\lambda)(X,Y):= \operatorname{Trace}(\lambda(X, \lambda(Y , \square)));
$$
the symmetric bilinear form $\kappa(\lambda)$ is just the Killing form of $(\mathbb{R}^4, \lambda)$.

Therefore $\kappa(  \mu_{13,+}(0) ) {\xrightarrow{\text{\,\,deg\,\,}} } \kappa( \mu_{6}(-\frac{1}{2},0) ) $
where:
$$
\kappa(\mu_{13,+}(0)) = -2 e_1^{\ast} \otimes e_{1}^\ast
$$
and
$$
\kappa(\mu_{6}(-\frac{1}{2},0) ) = \frac{3}{2} e_1^{\ast} \otimes e_{1}^\ast
$$
which is a contradiction because $\kappa( \mu_{13,+}(0) ) $ is negative semi-definite bilinear form while $\kappa( \mu_{6}(-\frac{1}{2},0) ) $ is non-zero positive semi-definite bilinear form.

We can give another proof by using $(\alpha,\beta,\gamma)$-derivations since $\operatorname{Dim}\operatorname{Der}(-1,1,1)(\mathfrak{d}^{'}_{4,0})=6$ and $\operatorname{Dim}\operatorname{Der}(-1,1,1)({\mathfrak{r}_{4,-1/2,-1/2}})=4$.

\end{enumerate}

\subsection{$(\mathfrak{d}_4,J_1)$}

\begin{enumerate}

\item
$
(\mathfrak{d}_4,J_1)
\xrightarrow{\text{\,\,deg\,\,}}
{(\mathfrak{rh}_{3}\times \mathbb{R},J)}
$\newline

Consider the curve $g_t = \left(
          \begin{array}{cc}
            A(t) & -B(t) \\
            B(t) & A(t) \\
          \end{array}
        \right)$ where $$A(t)= \operatorname{diag}(e^t,1)
                             \qquad
                             \mbox{and}
                             \qquad
                         B(t)= \left(
                                 \begin{array}{cc}
                                   0 & 0 \\
                                   -\frac{1}{2}e^{2t} & 0 \\
                                 \end{array}
                               \right)
                         .$$
Note that $g_t \in \operatorname{GL}(\mathbb{R}^4, {J_{\tiny{\mbox{cn}}}})$ and its action on $\mu_{10}(0)$ gives
$$
\left\{
\begin{array}{l}
[{\it e1},{\it e2}]=-{{\rm e}^{-t}}{\it e2},
[{\it e1},{\it e3}]={\it e2}+{{\rm e}^{-t}}{\it e3},
[{\it e2},{\it e3}]={{\rm e}^{-t}}{\it e4}
\end{array}
\right.
$$
Taking the limit $g_t \cdot \mu_{10}(0) $ as $t \to \infty$, we obtain $\mu_{1}$.

\end{enumerate}

\subsection{$(\mathfrak{r}^{'}_{4,-1/2,\delta},J_{k})$ with $\delta>0$ and $k=1,2$ }

\begin{enumerate}

\item
For any $\delta \in \mathbb{R}$, $\delta>0$,
$(\mathfrak{r}^{'}_{4,-1/2,\delta},J_{k_{}})
\xrightarrow{\text{\,\,deg\,\,}}
{(\mathfrak{rh}_{3}\times \mathbb{R},J)}
$\newline

Consider the curve $g_t = \left(
          \begin{array}{cc}
            A(t) & -B(t) \\
            B(t) & A(t) \\
          \end{array}
        \right)$ where $$A(t)= \left(
                                 \begin{array}{cc}
                                   {\rm e}^t & 0 \\
                                   {\frac {  \left( -1 \right)^{k}4\delta}{9+4\,{\delta}^{2}}}{{\rm e}^{2\,t}} & {\frac { \left( -1 \right) ^{k} 4 \delta}{9+4\,{\delta}^{2}}} \\
                                 \end{array}
                               \right)
                             \qquad
                             \mbox{and}
                             \qquad
                         B(t)= \left(
                                 \begin{array}{cc}
                                   0 & 0 \\
                                   {\frac {-6 }{9+4\,{\delta}^{2}}}{{\rm e}^{2\,t}} &  {\frac {-6 }{9+4\,{\delta}^{2}}}  \\
                                 \end{array}
                               \right)
                         .$$
Note that $g_t \in \operatorname{GL}(\mathbb{R}^4, {J_{\tiny{\mbox{cn}}}})$ and its action on $\mu_{6}(-\frac{1}{2},(-1)^k \delta)$ gives
$$
\left\{
\begin{array}{l}
[{\it e1},{\it e2}]=-\frac{1}{2}{{\rm e}^{-t}}{\it e2}+ \left( -1 \right) ^{k}{{\rm e}^{-t}}\delta {\it e4},\\
{[{\it e1},{\it e3}]}={\it e2}+{{\rm e}^{-t}}{\it e3},\\
{[{\it e1},{\it e4}]}=-\left( -1 \right) ^{k}{{\rm e}^{-t}}\delta{\it e2}- \frac{1}{2}{{\rm e}^{-t}}{\it e4}
\end{array}
\right.
$$
Taking the limit $g_t \cdot \mu_{6}(-\frac{1}{2},(-1)^k \delta)$ as $t \to \infty$, we obtain $\mu_1$.

\end{enumerate}

\subsection{${(\mathfrak{r}_{4,-1/2,-1/2}, J)}$ }

\begin{enumerate}
\item
${(\mathfrak{r}_{4,-1/2,-1/2}, J)}
\xrightarrow{\text{\,\,deg\,\,}}
{(\mathfrak{rh}_{3}\times \mathbb{R},J)}
$\newline

Consider the curve $g_t = \left(
          \begin{array}{cc}
            A(t) & -B(t) \\
            B(t) & A(t) \\
          \end{array}
        \right)$ where $$A(t)= \operatorname{diag}(e^t,1)
                             \qquad
                             \mbox{and}
                             \qquad
                         B(t)= \left(
                                 \begin{array}{cc}
                                   0 & 0 \\
                                   -\frac{2}{3}e^{2t} & 0 \\
                                 \end{array}
                               \right)
                         .$$
Note that $g_t \in \operatorname{GL}(\mathbb{R}^4, {J_{\tiny{\mbox{cn}}}})$ and its action on $\mu_{6}(-\frac{1}{2},0)$ gives
$$
\left\{
\begin{array}{l}
[{\it e1},{\it e2}]=-\frac{1}{2} {{\rm e}^{-t}}{\it e2},
[{\it e1},{\it e3}]={\it e2}+{{\rm e}^{-t}}{\it e3},
[{\it e1},{\it e4}]=-\frac{1}{2}{{\rm e}^{-t}}{\it e4}
\end{array}
\right.
$$
Taking the limit $g_t \cdot \mu_{6}(-\frac{1}{2},0) $ as $t \to \infty$, we obtain $\mu_1$.

\end{enumerate}

\subsection{${(\mathfrak{rr}^{'}_{3,0}\times \mathbb{R},J)}$ }

\begin{enumerate}
\item
$
{(\mathfrak{rr}^{'}_{3,0}\times \mathbb{R},J)}
\xrightarrow{\text{\,\,deg\,\,}}
{(\mathfrak{rh}_{3}\times \mathbb{R},J)}
$\newline

Consider the curve $g_t = \left(
          \begin{array}{cc}
            A(t) & -B(t) \\
            B(t) & A(t) \\
          \end{array}
        \right)$ where $$A(t)= \left(
                                 \begin{array}{cc}
                                   e^t & 0 \\
                                   e^{2t} & 1 \\
                                 \end{array}
                               \right)
                             \qquad
                             \mbox{and}
                             \qquad
                         B(t)= 0
                         .$$
Note that $g_t \in \operatorname{GL}(\mathbb{R}^4, {J_{\tiny{\mbox{cn}}}})$ and its action on $\mu_{7}(0)$ gives
$$
\left\{
\begin{array}{l}
[{\it e1},{\it e2}]={\frac {1}{{{\rm e}^{t}}}}{\it e4},
[{\it e1},{\it e3}]={\it e2},
[{\it e1},{\it e4}]=-{\frac {1}{{{\rm e}^{t}}}}{\it e2}
\end{array}
\right.
$$
Taking the limit $g_t \cdot \mu_{7}(0)$ as $t \to \infty$, we obtain $\mu_{1}$.

\end{enumerate}

\subsection{${(\mathfrak{rh}_{3}\times \mathbb{R},J)}$ }

\begin{enumerate}
\item
$
{(\mathfrak{rh}_{3}\times \mathbb{R},J)}
\xrightarrow{\text{\,\,deg\,\,}}
{(\mathfrak{a}_{4},J)}
$\newline

Consider the curve $g_t = \left(
          \begin{array}{cc}
            A(t) & -B(t) \\
            B(t) & A(t) \\
          \end{array}
        \right)$ where $$A(t)= \operatorname{diag}(1,e^{-t})
                             \qquad
                             \mbox{and}
                             \qquad
                         B(t)=0
                         .$$
Note that $g_t \in \operatorname{GL}(\mathbb{R}^4, {J_{\tiny{\mbox{cn}}}})$ and its action on $\mu_{1}$ gives
$$
\left\{
\begin{array}{l}
[{\it e1},{\it e3}]={\frac {1}{{{\rm e}^{t}}}}{\it e2}
\end{array}
\right.
$$
Taking the limit $g_t \cdot \mu_{1}$ as $t \to \infty$, we obtain $\mu_{0}$.

\end{enumerate}

\section{Abelian complex structure case}

\begin{table}[ht]
\centering
%\operatorname{Dim}
\begin{center}
    \begin{tabular}{c : c: p{6cm} : p{5cm} }
     $\operatorname{Dim} \operatorname{Der}_{J_{\tiny{\mbox{cn}}}}$ & $\operatorname{Dim} \operatorname{Der}$ & \multicolumn{2}{c}{Lie Algebras  endowed abelian complex structures}\\
    \hdashline
     \multirow{1}{*}{$0$} & 4   &  $(\mathfrak{r}_2\times \mathfrak{r}_{2},J)$,  & $(\mathfrak{r}^{'}_{2},J_2)$\\
    \hdashline
    \multirow{1}{*}{$1$}  & 5   &                                                                            & $(\mathfrak{d}_{4,1},J)$\\
    \hdashline
    \multirow{2}{*}{$2$}  & 4   &                                              & $(\mathfrak{r}^{'}_{2},J_1(0,1))$  \\
                          & 8   &  $(\mathfrak{r}_2\times \mathbb{R}^{2},J)$   &                                     \\
    \hdashline
    \multirow{1}{*}{$4$}  & 10  &  $(\mathfrak{h}_{3}\times \mathbb{R},J)$   & \\
    \hdashline
    \multirow{1}{*}{$8$}  & 16  &  $(\mathfrak{a}_{4},J)$  & \\
    \end{tabular}
\end{center}
{\scriptsize Table 2: Dimension of  complex derivations.}
\end{table}

\begin{center}
\begin{tikzpicture}
\node[circle,fill,inner sep=0pt,minimum size=3pt,label=above:{\scriptsize$(\mathfrak{r}_2\times \mathfrak{r}_{2},J)$}] (r2r2) at (-1.5,0) {};
\node[circle,fill,inner sep=0pt,minimum size=3pt,label=above:{\scriptsize$(\mathfrak{r}^{'}_{2},J_2)$}] (r2p3) at (+1.5,0) {};
\node[ label=above:{\scriptsize$(\mathfrak{d}_{4,1},J)$}] (d41label) at (0,-0.9) {};
\node[circle,fill,inner sep=0pt,minimum size=3pt, label=above:{}] (d41) at (0,-1) {};
\node[circle,fill,inner sep=0pt,minimum size=3pt,label=left:{\scriptsize$(\mathfrak{r}_2\times \mathbb{R}^{2},J)$}] (r2z) at (-1.5,-2) {};
\node[circle,fill,inner sep=0pt,minimum size=3pt,label=right:{\scriptsize$(\mathfrak{r}^{'}_{2},J_1(0,1))$}] (r2p1) at (+1.5,-2) {};
\node[circle,fill,inner sep=0pt,minimum size=3pt,label=left:{\scriptsize$(\mathfrak{h}_{3}\times \mathbb{R},J)$ }] (h3) at (0,-3) {};
\node[circle,fill,inner sep=0pt,minimum size=3pt,label=below:{\scriptsize$(\mathfrak{a}_{4},J)$ }] (a4) at (0,-4) {};
\draw [->] (r2r2) edge (d41) (r2r2) edge (r2z) (r2p3) edge (d41) (r2p3) edge  (r2p1);
\draw [->] (d41) edge (h3) (r2z) edge (h3) (r2p1) edge (h3)   (h3) edge (a4);

\end{tikzpicture}
\end{center}

\subsection{$(\mathfrak{r}_2\times \mathfrak{r}_{2},J)$}

\begin{enumerate}

\item
$
(\mathfrak{r}_2\times \mathfrak{r}_{2},J)
\xrightarrow{\text{\,\,deg\,\,}}
(\mathfrak{d}_{4,1},J)
$\newline

Consider the curve $g_t = \left(
          \begin{array}{cc}
            A(t) & -B(t) \\
            B(t) & A(t) \\
          \end{array}
        \right)$ where $$A(t)=\left(
                                \begin{array}{cc}
                                  1 & 0 \\
                                  e^t & -e^t \\
                                \end{array}
                              \right)
                             \qquad
                             \mbox{and}
                             \qquad
                         B(t)=\left(
                                \begin{array}{cc}
                                  0 & 0 \\
                                  -e^{-t} & 0 \\
                                \end{array}
                              \right)
                         .$$
Note that $g_t \in \operatorname{GL}(\mathbb{R}^4, {J_{\tiny{\mbox{cn}}}})$ and its action on $\mu_2$ gives
$$
\left\{
\begin{array}{l}
[{\it e1},{\it e2}]={{\rm e}^{-2\,t}}{\it e4},
[{\it e1},{\it e3}]={{\rm e}^{-t}}{\it e2}+{\it e3}-{{\rm e}^{-3\,t}}{\it e4},
[{\it e1},{\it e4}]={\it e4},\\
{[{\it e2},{\it e3}]}={\it e4},
[{\it e2},{\it e4}]=-{{\rm e}^{-t}}{\it e4},
[{\it e3},{\it e4}]={{\rm e}^{-2\,t}}{\it e4}
\end{array}
\right.
$$
Taking the limit $g_t \cdot \mu_2$ as $t \to \infty$, we obtain $\mu_{10}(1)$.

\item
$
(\mathfrak{r}_2\times \mathfrak{r}_{2},J)
\centernot{ \xrightarrow{\text{\,\,deg\,\,}} }
(\mathfrak{r}^{'}_{2},J_1(0,1))
$\newline

Because $
\mathfrak{r}_2\times \mathfrak{r}_{2}
\centernot{ \xrightarrow{\text{\,\,deg\,\,}} }
\mathfrak{r}^{'}_{2}
$, since $\operatorname{Dim}\operatorname{Der}( \mathfrak{r}_2\times \mathfrak{r}_{2} )= 4 $ and $\operatorname{Dim}\operatorname{Der}(\mathfrak{r}^{'}_{2})= 4$,
and $\mathfrak{r}_2\times \mathfrak{r}_{2} \not\equiv \mathfrak{r}^{'}_{2}$.

\item
$
(\mathfrak{r}_2\times \mathfrak{r}_{2},J)
\xrightarrow{\text{\,\,deg\,\,}}
(\mathfrak{r}_2\times \mathbb{R}^{2},J)
$\newline

Consider the curve $g_t = \left(
          \begin{array}{cc}
            A(t) & -B(t) \\
            B(t) & A(t) \\
          \end{array}
        \right)$ where $$A(t)= \operatorname{diag}(1,e^t)
                             \qquad
                             \mbox{and}
                             \qquad
                         B(t)=0
                         .$$
Note that $g_t \in \operatorname{GL}(\mathbb{R}^4, {J_{\tiny{\mbox{cn}}}})$ and its action on $\mu_{2}$ gives
$$
\left\{
\begin{array}{l}
[{\it e1},{\it e3}]={\it e3},
[{\it e2},{\it e4}]={{\rm e}^{-t}}{\it e4}
\end{array}
\right.
$$
Taking the limit $g_t \cdot \mu_{2}$ as $t \to \infty$, we obtain $ \mu_{6}(0,0) $.

\end{enumerate}

\subsection{$(\mathfrak{r}^{'}_{2},J_2)$}

\begin{enumerate}

\item
$
(\mathfrak{r}^{'}_{2},J_2)
\xrightarrow{\text{\,\,deg\,\,}}
(\mathfrak{d}_{4,1},J)
$\newline

Consider the curve $g_t = \left(
          \begin{array}{cc}
            A(t) & -B(t) \\
            B(t) & A(t) \\
          \end{array}
        \right)$ where $$A(t)=\operatorname{diag}(1,e^t)
                             \qquad
                             \mbox{and}
                             \qquad
                         B(t)=0
                         .$$
Note that $g_t \in \operatorname{GL}(\mathbb{R}^4, {J_{\tiny{\mbox{cn}}}})$ and its action on $\mu_4$ gives
$$
\left\{
\begin{array}{l}
[{\it e1},{\it e3}]={\it e3},
[{\it e1},{\it e4}]={\it e4},
[{\it e2},{\it e3}]={\it e4},
[{\it e2},{\it e4}]=-{{\rm e}^{-2\,t}}{\it e3}
\end{array}
\right.
$$
Taking the limit $g_t \cdot \mu_{4}$ as $t \to \infty$, we obtain $\mu_{10}(1)$.

\item
$
(\mathfrak{r}^{'}_{2},J_2)
\xrightarrow{\text{\,\,deg\,\,}}
(\mathfrak{r}^{'}_{2},J_1(0,1))
$\newline

Consider the curve $g_t = \left(
          \begin{array}{cc}
            A(t) & -B(t) \\
            B(t) & A(t) \\
          \end{array}
        \right)$ where $$A(t)=\operatorname{diag}(1,e^t)
                             \qquad
                             \mbox{and}
                             \qquad
                         B(t)=\left(
                                \begin{array}{cc}
                                  0 & 0 \\
                                  e^t & 0 \\
                                \end{array}
                              \right)
                         .$$
Note that $g_t \in \operatorname{GL}(\mathbb{R}^4, {J_{\tiny{\mbox{cn}}}})$ and its action on $\mu_{4}$ gives
$$
\left\{
\begin{array}{l}
[{\it e1},{\it e2}]={\it e2}-{{\rm e}^{-t}}{\it e3},
[{\it e1},{\it e4}]={\it e4},
[{\it e2},{\it e3}]={\it e4},
[{\it e2},{\it e4}]={{\rm e}^{-t}}{\it e2}-{{\rm e}^{-2\,t}}{\it e3},\\
{[{\it e3},{\it e4}]}={\it e2}-{{\rm e}^{-t}}{\it e3}
\end{array}
\right.
$$
Taking the limit $g_t \cdot \mu_{4}$ as $t \to \infty$, we obtain $\mu_{3}(0,1)$.

\item
$
(\mathfrak{r}^{'}_{2},J_2)
\centernot{ \xrightarrow{\text{\,\,deg\,\,}} }
(\mathfrak{r}_2\times \mathbb{R}^{2},J)
$\newline

In fact, this holds because $
\mathfrak{r}^{'}_{2}
\centernot{ \xrightarrow{\text{\,\,deg\,\,}} }
\mathfrak{r}_2\times \mathbb{R}^{2}
$.  Let us now explore this in detail. Suppose for the sake of contradiction that $
\mathfrak{r}^{'}_{2}
{ \xrightarrow{\text{\,\,deg\,\,}} }
\mathfrak{r}_2\times \mathbb{R}^{2}
$. First, consider the \(\operatorname{GL}(\mathbb{R}^4)\)-equivariant continuous function \[ \kappa_{g}: C^{2}(\mathbb{R}^{4};\mathbb{R}^{4}) \rightarrow \Sigma^2(\mathbb{R}^{4}; \mathbb{R}). \]
By assumption, we would have $ \kappa_{g}( \mu_{4} ) {\xrightarrow{\text{\,\,deg\,\,}} } \kappa( \mu_{6}(0,0) ) $
where:
$$
\kappa_{g}(\mu_{4} ) = (2+4g)e_1^{\ast}\otimes e_{1}^{\ast} - 2e_2^{\ast}\otimes e_{2}^{\ast}
$$
and
$$
\kappa_{g}( \mu_{6}(0,0) ) = (1+g) e_1^{\ast}\otimes e_{1}^{\ast}
$$

Setting $g= -\frac{1}{2}$, we obtain a contradiction, since $\kappa_{-\frac{1}{2}}(\mu_{4} )$ is negative semi-definite bilinear form while $\kappa_{-\frac{1}{2}}(\mu_{6}(0,0)  ) $ is non-zero positive semi-definite bilinear form.

\end{enumerate}

\subsection{$(\mathfrak{d}_{4,1},J)$}

\begin{enumerate}

  \item
$
(\mathfrak{d}_{4,1},J)
\centernot{ \xrightarrow{\text{\,\,deg\,\,}} }
(\mathfrak{r}^{'}_{2},J_1(0,1))
$\newline

Because $
\mathfrak{d}_{4,1}
\centernot{ \xrightarrow{\text{\,\,deg\,\,}} }
\mathfrak{r}^{'}_{2}
$, since $\operatorname{Dim}\operatorname{Der}( \mathfrak{d}_{4,1} )=5$ and $\operatorname{Dim}\operatorname{Der}(\mathfrak{r}^{'}_{2})=4$.

\item
$
(\mathfrak{d}_{4,1},J)
\centernot{ \xrightarrow{\text{\,\,deg\,\,}} }
(\mathfrak{r}_2\times \mathbb{R}^{2},J)
$\newline

In fact, this holds because $
\mathfrak{d}_{4,1}
\centernot{ \xrightarrow{\text{\,\,deg\,\,}} }
\mathfrak{r}_2\times \mathbb{R}^{2}
$.  Let us now explore this in detail. Suppose for the sake of contradiction that $
\mathfrak{d}_{4,1}
{ \xrightarrow{\text{\,\,deg\,\,}} }
\mathfrak{r}_2\times \mathbb{R}^{2}
$. First, consider the \(\operatorname{GL}(\mathbb{R}^4)\)-equivariant continuous function \[ \kappa_{g}: C^{2}(\mathbb{R}^{4};\mathbb{R}^{4}) \rightarrow \Sigma^2(\mathbb{R}^{4}; \mathbb{R}). \]
By assumption, we would have $ \kappa_{g}( \mu_{10}(1) ) {\xrightarrow{\text{\,\,deg\,\,}} } \kappa( \mu_{6}(0,0) ) $
where:
$$
\kappa_{g}( \mu_{10}(1) ) = (2+4g)e_{1}^\ast \otimes  e_{1}^\ast
$$
and
$$
\kappa_{g}( \mu_{6}(0,0) ) = (1+g)e_{1}^\ast \otimes  e_{1}^\ast
$$

Setting $g= -\frac{1}{2}$, we obtain a contradiction, since $\kappa_{ -\frac{1}{2}}(  \mu_{10}(1) ) $ is the zero bilinear form while $\kappa_{-\frac{1}{2}}(  \mu_{6}(0,0) ) $ is non-zero positive semi-definite bilinear form.

\item
$
(\mathfrak{d}_{4,1},J)
{ \xrightarrow{\text{\,\,deg\,\,}} }
  (\mathfrak{h}_{3}\times \mathbb{R},J)
$\newline

Consider the curve $g_t = \left(
          \begin{array}{cc}
            A(t) & -B(t) \\
            B(t) & A(t) \\
          \end{array}
        \right)$ where $$A(t)= \operatorname{diag}(e^t,1)
                             \qquad
                             \mbox{and}
                             \qquad
                         B(t)=\left(
                                \begin{array}{cc}
                                  0  & 0 \\
                                  -e^{2t} & 0 \\
                                \end{array}
                              \right)
                         .$$
Note that $g_t \in \operatorname{GL}(\mathbb{R}^4, {J_{\tiny{\mbox{cn}}}})$ and its action on $\mu_{10}(1)$ gives
$$
\left\{
\begin{array}{l}
[{\it e1},{\it e3}]={\it e2}+{\frac {1}{{{\rm e}^{t}}}}{\it e3},
[{\it e1},{\it e4}]={\frac {1}{{{\rm e}^{t}}}}{\it e4},
[{\it e2},{\it e3}]={\frac {1}{{{\rm e}^{t}}}}{\it e4}
\end{array}
\right.
$$
Taking the limit $g_t \cdot \mu_{10}(1)$ as $t \to \infty$, we obtain $\mu_{1}$.

\end{enumerate}

\subsection{$(\mathfrak{r}^{'}_{2},J_1(0,1)) $}
\begin{enumerate}

  \item
$
(\mathfrak{r}^{'}_{2},J_1(0,1))
%\centernot
{ \xrightarrow{\text{\,\,deg\,\,}} }
(\mathfrak{h}_{3}\times \mathbb{R},J)
$\newline

Consider the curve $g_t = \left(
          \begin{array}{cc}
            A(t) & -B(t) \\
            B(t) & A(t) \\
          \end{array}
        \right)$ where $$A(t)= \operatorname{diag}(e^t,1)
                             \qquad
                             \mbox{and}
                             \qquad
                         B(t)=\left(
                                \begin{array}{cc}
                                  0 & 0 \\
                                  \frac{1}{2}e^{2t} & 0 \\
                                \end{array}
                              \right)
                         .$$
Note that $g_t \in \operatorname{GL}(\mathbb{R}^4, {J_{\tiny{\mbox{cn}}}})$ and its action on $\mu_{3}(0,1)$ gives
$$
\left\{
\begin{array}{l}
[{\it e1},{\it e2}]={{\rm e}^{-t}}{\it e2},
[{\it e1},{\it e3}]={\it e2},
[{\it e1},{\it e4}]={{\rm e}^{-t}}{\it e4},\\
{[{\it e2},{\it e3}]}={{\rm e}^{-t}}{\it e4},
[{\it e3},{\it e4}]={{\rm e}^{-t}}{\it e2}
\end{array}
\right.
$$
Taking the limit $g_t \cdot \mu_{3}(0,1)$ as $t \to \infty$, we obtain $\mu_{1}$.

\end{enumerate}

\subsection{$(\mathfrak{r}_2\times \mathbb{R}^{2},J) $}
\begin{enumerate}

  \item
$
(\mathfrak{r}_2\times \mathbb{R}^{2},J)
%\centernot
{ \xrightarrow{\text{\,\,deg\,\,}} }
(\mathfrak{h}_{3}\times \mathbb{R},J)
$\newline

Consider the curve $g_t = \left(
          \begin{array}{cc}
            A(t) & -B(t) \\
            B(t) & A(t) \\
          \end{array}
        \right)$ where $$A(t)=\operatorname{diag}(e^t,1)
                             \qquad
                             \mbox{and}
                             \qquad
                         B(t)=\left(
                                \begin{array}{cc}
                                  0 & 0 \\
                                  -e^{2t} & 0 \\
                                \end{array}
                              \right)
                         .$$
Note that $g_t \in \operatorname{GL}(\mathbb{R}^4, {J_{\tiny{\mbox{cn}}}})$ and its action on $ \mu_{6}(0,0) $ gives
$$
\left\{
\begin{array}{l}
[{\it e1},{\it e3}]={\it e2}+{{\rm e}^{-t}}{\it e3}
\end{array}
\right.
$$
Taking the limit $g_t \cdot \mu_{6}(0,0) $ as $t \to \infty$, we obtain $\mu_{1} $.

\end{enumerate}

\section{General Case}

\begin{landscape}

\begin{table}[ht]
\centering
\begin{center}
    \begin{tabular}{c : c: p{4.5cm} : p{3.5cm} : p{6cm} : p{4cm} }

      \multirow{2}{*}{ \begin{tabular}{l}$\operatorname{Dim}$ \\ $\operatorname{Der}_{J_{\tiny{\mbox{cn}}}}$ \end{tabular}} &
      \multirow{2}{*}{ \begin{tabular}{l}$\operatorname{Dim}$ \\ $\operatorname{Der}$ \end{tabular} } &
     \multicolumn{4}{c}{$4$-dimensional Lie Algebras endowed with complex strucuture}\\

     \multicolumn{2}{c:}{}&
     \multicolumn{2}{c:}{}&
     Unimodular &
     Abelian complex strcuture\\
    \hdashline
     \multirow{1}{*}{$0$} & 4   &   & &  & $(\mathfrak{r}_2\times \mathfrak{r}_{2},J)$,  $(\mathfrak{r}^{'}_{2},J_2)$  \\
    \hdashline
     \multirow{3}{*}{$1$} & 4  &  &  &$(\mathfrak{sl}(2,\mathbb{R})\times \mathbb{R}, J(t))$ , $(\mathfrak{so}(3,\mathbb{R})\times \mathbb{R}, J(t))$ &  $(\mathfrak{d}_{4,1},J)$ \\
                          & 5  &  $(\mathfrak{ d}^{'}_{4,\delta},J_{k_{4}})$, $(\mathfrak{ d}_{4,\lambda}, J_{k_{2}})$, $(\mathfrak{ h}_{ 4},J)$
                               &
                               & $(\mathfrak{ d}^{'}_{4,0},J_{k_{2}})$, $(\mathfrak{d}_4,J_2)$ &  \\
                          & 7  & $(\mathfrak{ d}_{4,1/2}, J_2)$,$(\mathfrak{ d}_{4,1/2}, J_3)$ &  & &  \\
   \hdashline
     \multirow{5}{*}{$2$} & 4  &
                               &$(\mathfrak{r}^{'}_{2 }, J_{1}(a,b))$, $(\mathfrak{r}^{'}_{2 }, J_{3})$
                               &
                               &$(\mathfrak{r}^{'}_{2 }, J_{1}(0,1))$ \\
                          & 5  &
                               &
                               &$(\mathfrak{d}_4,J_1)$
                               & \\
                          & 6  & $(\mathfrak{r}^{'}_{4,c,d},J_{k_{2}})$
                               & $(\mathfrak{rr}^{'}_{3,a}\times\mathbb{R},J_{k_{2}})$,
                               & $(\mathfrak{r}^{'}_{4,-1/2,d},J_{k_{2}})$, $(\mathfrak{rr}^{'}_{3,0}\times\mathbb{R},J)$
                               &  \\
                          & 7  &  $(\mathfrak{ d}_{4,1/2}, J_1)$ &  & & \\
                          & 8  &  $(\mathfrak{ r}_{4,\alpha,\alpha },J)$, $(\mathfrak{r}_{4,\alpha,1},J)$, $(\mathfrak{ r}_{4,1},J)$
                               &  $(\mathfrak{rr}_{3,1 }\times \mathbb{R},J)$
                               &  $(\mathfrak{ r}_{4,-1/2,-1/2},J)$
                               &  $(\mathfrak{r}_2\times \mathbb{R}^{2},J)$ \\
   \hdashline
     \multirow{3}{*}{$4$} & 4  &    & $(\mathfrak{r}^{'}_{2 }, J_{1}(0,-1))$  &  & \\
                          & 10  &  \multicolumn{4}{c}{$(\mathfrak{h}_{3}\times \mathbb{R},J)$}  \\
                          & 12  &  $(\mathfrak{ r}_{4,1,1},J)$ &  & & \\
      \hdashline
   $8$ & $16$ & \multicolumn{4}{c}{$(\mathfrak{a}_{4},J)$}\\
    \end{tabular}
\end{center}
{\scriptsize Table 3: Dimension of  complex derivations. \\ With $k_{j}=1,2,\ldots,j$, $(a,b) \neq (0,-1)$ and $\lambda \neq 1,1/2$}
\end{table}
\end{landscape}

\subsection{$(\mathfrak{ d}^{'}_{4,\delta},J_{k})$ with $\delta>0$ and ${k}=1,\ldots,4$}\label{delta4prima}
\begin{enumerate}

\item
$
(\mathfrak{ d}^{'}_{4,\delta},J_{k})
\centernot
{ \xrightarrow{\text{\,\,deg\,\,}} }
(\mathfrak{d}_4,J_1)
$\newline

Because $
\mathfrak{ d}^{'}_{4,\delta}
\centernot
{ \xrightarrow{\text{\,\,deg\,\,}} }
\mathfrak{d}_4
$, since $\operatorname{Dim}\operatorname{Der}( \mathfrak{ d}^{'}_{4,\delta} )= 5 $ and $\operatorname{Dim}\operatorname{Der}( \mathfrak{d}_4 )=5$, and
$ \mathfrak{ d}^{'}_{4,\delta} \not \equiv \mathfrak{d}_4 $.

 \item
$
(\mathfrak{ d}^{'}_{4,\delta},J_{k})
%\centernot
{ \xrightarrow{\text{\,\,deg\,\,}} }
(\mathfrak{r}^{'}_{4,c,d},J_{k_{2}})
$ with $(c,d) = (1/2,\delta)$ and  $ \lfloor k/2 \rfloor \, {\equiv} \, k_{2} (\operatorname{mod} 2)$ \newline

Consider the curve $g_t = \left(
          \begin{array}{cc}
            A(t) & -B(t) \\
            B(t) & A(t) \\
          \end{array}
        \right)$ where $$A(t)= \operatorname{diag}\left( (-1)^{\lfloor \frac{k+1}{2} \rfloor} \frac{1}{\delta} , e^t  \right)
                             \qquad
                             \mbox{and}
                             \qquad
                         B(t)=0
                         .$$
Note that $g_t \in \operatorname{GL}(\mathbb{R}^4, {J_{\tiny{\mbox{cn}}}})$ and its action on $\mu_{13, (-1)^k}((-1)^{ \lfloor \frac{k+1}{2} \rfloor} \frac{1}{2\delta} )$ gives
$$
\left\{
\begin{array}{l}
[{\it e1},{\it e2}]=\frac{1}{2}\,{\it e2}+ \left( -1 \right) ^{ \lfloor \frac{k}{2}\rfloor }\delta\,{\it e4},
[{\it e1},{\it e3}]={\it e3},
[{\it e1},{\it e4}]=- \left( -1 \right) ^{ \lfloor \frac{k}{2}\rfloor  }\delta\,{\it e2}+ \frac{1}{2}\,{\it e4},\\
{[{\it e2},{\it e4}]}={\rm e}^{-2\,t} { \frac {{} \left( -1 \right) ^{  \left\lfloor \frac{k+1}{2} \right\rfloor }}{\delta}}{\it e3}
\end{array}
\right.
$$
Taking the limit $g_t \cdot \mu_{13, (-1)^k}((-1)^{ \lfloor \frac{k+1}{2} \rfloor} \frac{1}{2\delta} )$ as $t \to \infty$, we obtain $\mu_{6}(\frac{1}{2}, (-1)^{ \lfloor k/2\rfloor} \delta) $.

  \item\label{delta4primar4prima}
$
(\mathfrak{ d}^{'}_{4,\delta},J_{k})
\centernot
{ \xrightarrow{\text{\,\,deg\,\,}} }
(\mathfrak{r}^{'}_{4,c,d},J_{k_{2}})
$ with $(c,d) \neq (1/2,\delta)$ or $ \lfloor k/2 \rfloor  \centernot{\equiv} k_{2} (\operatorname{mod} 2)$ \newline

Suppose for the sake of contradiction that $
(\mathfrak{ d}^{'}_{4,\delta},J_{k})
%\centernot
{ \xrightarrow{\text{\,\,deg\,\,}} }
(\mathfrak{r}^{'}_{4,c,d},J_{k_{2}})
$
. First, we consider the $\operatorname{GL}(\mathbb{R}^4, {J_{\tiny{\mbox{cn}}}})$-equivariant continuous function $\psi_{\alpha,\alpha}: C^{2}_{J_{\tiny{\mbox{cn}}}}(\mathbb{R}^{4};\mathbb{R}^{4}) \rightarrow C^{2}(\mathbb{R}^{4};\mathbb{R}^{4})$ with $\alpha = \frac{(-1)^{ \lfloor \frac{k}{2}\rfloor}}{2\delta}$. Therefore $$\mathfrak{A}:=(\mathbb{R}^{4}, \psi_{ \alpha, \alpha }( \mu_{13, (-1)^k}((-1)^{ \lfloor \frac{k+1}{2} \rfloor} \frac{1}{2\delta} )  ) { \xrightarrow{\text{\,\,deg\,\,}} } \mathfrak{B}:=(\mathbb{R}^{4}, \psi_{ \alpha ,\alpha}( \mu_{6}(c,d) ),$$
where:
$$
\mathfrak{A} =
\left\{
\begin{array}{l}
[{\it e1},{\it e2}]=-\frac{1}{2}{\frac { \left( -1 \right) ^{ \lfloor\frac{k+1}{2}\rfloor } }{\delta}}{\it e2}
+ \frac{1}{2}\,{\frac { \left( 2 \, \left( -1 \right) ^{k}{\delta}^{2}+ \left( -1 \right) ^{k} \right) }{{\delta}^{2}}}{\it e4},\\
{[{\it e1},{\it e3}]}=\frac{1}{2}\,{\frac { \left( -1 \right) ^{1+k}}{{\delta}^{2}}}{\it e1} + {\frac { \left( -1 \right) ^{ \lfloor\frac{k+1}{2}\rfloor } }{\delta}}{\it e3},\\
{[{\it e1},{\it e4}]}=-\frac{1}{2}\,{\frac { \left( 2\, \left( -1 \right) ^{k}{\delta}^{2}-
 \left( -1 \right) ^{1+k} \right) }{{\delta}^{2}}}{\it e2}  -\frac{1}{2}\,\frac { \left( -1 \right) ^{ \lfloor\frac{k+1}{2}\rfloor } }{\delta} { \it e4},\\
{[{\it e2},{\it e3}]}=\frac{1}{4}\,{\frac { \left( -1 \right) ^{k}}{{\delta}^{2}}}{\it e2}+ \frac{1}{2} \,
{\frac { \left( -1 \right)^{ \lfloor \frac{k}{2}\rfloor +k} }{\delta}}{\it e4},\\
{[{\it e2},{\it e4}]}=\frac{1}{2}\,{\frac { \left( -1 \right) ^{ 1 + \lfloor \frac{k}{2}\rfloor }}{\delta}}{\it e1}+{\it e3},\\
{[{\it e3},{\it e4}]}=\frac{1}{2}\,{\frac
{ \left( -1 \right) ^{ \lfloor \frac{k}{2}\rfloor +k}}{\delta}}{\it e2} + \frac{1}{4}\,{\frac { \left( -1 \right) ^{1+k}}{{\delta}^{2}}}{\it e4}
\end{array}
\right.
$$
and
$$
\mathfrak{B} =
\left\{
\begin{array}{l}
[{\it e1},{\it e2}]=-{\frac { \left( -c\delta+ \left( -1 \right)^{ \lfloor k/2 \rfloor} d \right) }{\delta}}{\it e2}
+{\frac { \left( d\delta+ \left( -1 \right)^{ \lfloor k/2 \rfloor}c \right) }{\delta}}{\it e4},\\
{[{\it e1},{\it e3}]}=\frac{1}{2}\,{ \frac { \left( -1 \right) ^{ 1 + \lfloor k/2\rfloor }}{\delta}}{\it e1}+{\it e3},\\
{[{\it e1},{\it e4}]}=-{\frac { \left( d\delta+  \left( -1 \right)^{ \lfloor k/2 \rfloor}c \right) }{\delta}}{\it e2}
-{\frac { \left( -c\delta + \left( -1 \right)^{ \lfloor k/2 \rfloor}d \right) }{\delta}}{\it e4},\\
{[{\it e2},{\it e3}]}=\frac{1}{2} \, {\frac { \left( -1 \right)^{ \lfloor k/2 \rfloor} c }{\delta}}{\it e2}
+ \frac{1}{2}\,{\frac { \left( -1 \right)^{ \lfloor k/2 \rfloor}d }{\delta}}{\it e4},\\
{[{\it e3},{\it e4}]} = \frac{1}{2} \,{\frac {\left( -1 \right)^{ \lfloor k/2 \rfloor}d }{\delta}}{\it e2}
+\frac{1}{2}\,\frac { \left( -1 \right) ^{1 + \lfloor k/2 \rfloor  }c }{\delta}{\it e4}
\end{array}
\right.
$$
Since $\mathfrak{A}$ is a unimodular algebra, it follows $\mathfrak{B}$ is also unimodular, and so:$$  \left\{ {\frac {\delta+2\,c\delta+2\,
\left( -1 \right) ^{ 1 +  \lfloor k/2 \rfloor}d}{\delta}},
\frac{1}{2}\,{\frac { \left( -1  \right) ^{1 + \lfloor k/2 \rfloor } \left( -1+2\,c  \right) }{\delta}} \right\} $$
Solving for $c$ and $d$, we have $c=\frac{1}{2}$ and $d=(-1)^{ \lfloor k/2\rfloor}\delta$. Thus $(\mathfrak{r}^{'}_{4,c,d}, J_{k_{2}} )$ is
$( \mathfrak{r}^{'}_{4,\frac{1}{2}, \delta} , J_{k_{2}} )$ with  $k_{2} \equiv \lfloor k/2\rfloor (\operatorname{mod} 2)$; which is a contradiction.

\item
$
(\mathfrak{ d}^{'}_{4,\delta},J_{k})
\centernot
{ \xrightarrow{\text{\,\,deg\,\,}} }
(\mathfrak{rr}^{'}_{3,a}\times\mathbb{R},J_{k_{2}})
$ with $a>0$ and $k_2=1,2$ \newline

In fact, this holds because $
\mathfrak{ d}^{'}_{4,\delta}
\centernot
{ \xrightarrow{\text{\,\,deg\,\,}} }
\mathfrak{rr}^{'}_{3,a}\times\mathbb{R}
$.  Let us now explore this in detail. Suppose for the sake of contradiction that $
\mathfrak{ d}^{'}_{4,\delta}
{ \xrightarrow{\text{\,\,deg\,\,}} }
\mathfrak{rr}^{'}_{3,a}\times\mathbb{R}
$. We have two cases:

\begin{itemize}
  \item If $a \neq \frac{2}{\sqrt{1+4\delta^2}}$

First, consider the \(\operatorname{GL}(\mathbb{R}^4)\)-equivariant continuous function \[ \kappa_{g}: C^{2}(\mathbb{R}^{4};\mathbb{R}^{4}) \rightarrow \Sigma^2(\mathbb{R}^{4}; \mathbb{R}). \]
By assumption, we would have $ \kappa_{g}( \mu_{13,+}(1/2\delta) ) {\xrightarrow{\text{\,\,deg\,\,}} } \kappa( \mu_{7}(a) ) $
where:
$$
\kappa_{g}(  \mu_{13,+}(1/2\delta) ) = -\frac{1}{2}\,{\frac { \left( 4\,{\delta}^{2}-3-8\,g \right) }{{
\delta}^{2}}}e_{1}^\ast \otimes e_{1}^\ast
$$
and
$$
\kappa_{g}( \mu_{7}(a) ) = (2\,{a}^{2}-2+4\,g{a}^{2}) e_{1}^\ast \otimes e_{1}^\ast
$$

Setting $g= \frac{1}{2}{\delta}^{2}-\frac{3}{8}$, we obtain a contradiction, since $\kappa_{g}(  \mu_{13,+}(1/2\delta) ) $ is zero bilinear form while $\kappa_{g}( \mu_{7}(a) ) = \frac{1}{2} \left( {a}^{2}-4+4\,{a}^{2}{\delta}^{2} \right) e_{1}^\ast \otimes e_{1}^\ast \neq 0$ since $a \neq \frac{2}{\sqrt{1+4\delta^2}}$.

  \item If $a = \frac{2}{\sqrt{1+4\delta^2}}$.

We consider the \(\operatorname{GL}(\mathbb{R}^4)\)-equivariant continuous function \[ \Psi_{f,g}: C^{2}(\mathbb{R}^{4};\mathbb{R}^{4}) \rightarrow \Sigma^{3}(\mathbb{R}^{4};\mathbb{R}). \]
Therefore $$ \Psi_{f,g}(  \mu_{13,+}(1/2\delta) ) \xrightarrow{\text{\,\,deg\,\,}} \Psi_{f,g}( \mu_{7}(  \frac{2}{\sqrt{1+4\delta^2}} )  ), $$
where:
$$
 \Psi_{f,g}( \mu_{13,+}(1/2\delta) ) = -\frac{1}{2}\,{\frac {  \left( -16\,g-18\,f+24\,f{\delta}^
{2}+36\,{\delta}^{2}-15 \right) }{{\delta}^{3}}} e_{1}^\ast \otimes e_{1}^\ast \otimes e_{1}^\ast
$$
and
$$
\Psi_{f,g}(  \mu_{7}(  \frac{2}{\sqrt{1+4\delta^2}} )= -8\,{\frac { \left( -8\,g-9\,f+12\,f{\delta}^{2}-
3+36\,{\delta}^{2} \right) }{ \left( 1+4\,{\delta}^{2} \right) ^{3/2}}
} e_{1}^\ast \otimes e_{1}^\ast \otimes e_{1}^\ast
$$
Taking $f=0$  and $g= \frac{9}{4} \delta^2 - \frac{15}{16}$, we obtain a contradiction, since $ \Psi_{0,\frac{9}{4} \delta^2 - \frac{15}{16}}(\mu_{13,+}(1/2\delta)) \equiv 0$
and $ \Psi_{0,\frac{9}{4} \delta^2 - \frac{15}{16}}( \mu_{7}(  \frac{2}{\sqrt{1+4\delta^2}}  ) = -\frac{36}{\sqrt{1+4 \delta^2}}  e_{1}^\ast \otimes e_{1}^\ast \otimes e_{1}^\ast $.

\end{itemize}

  \item
$
(\mathfrak{ d}^{'}_{4,\delta},J_{k})
\centernot
{ \xrightarrow{\text{\,\,deg\,\,}} }
(\mathfrak{r}^{'}_{4,-1/2,d},J_{k_{2}})
$\newline
ver Subsection \ref{delta4prima}, item \ref{delta4primar4prima}.

\item
$
(\mathfrak{ d}^{'}_{4,\delta},J_{k})
\centernot
{ \xrightarrow{\text{\,\,deg\,\,}} }
(\mathfrak{rr}^{'}_{3,0}\times\mathbb{R},J)
$\newline

In fact, this holds because $
\mathfrak{ d}^{'}_{4,\delta}
\centernot
{ \xrightarrow{\text{\,\,deg\,\,}} }
\mathfrak{rr}^{'}_{3,0}\times\mathbb{R}
$.  Let us now explore this in detail. Suppose for the sake of contradiction that $
\mathfrak{ d}^{'}_{4,\delta}
{ \xrightarrow{\text{\,\,deg\,\,}} }
\mathfrak{rr}^{'}_{3,0}\times\mathbb{R}
$. First, consider the \(\operatorname{GL}(\mathbb{R}^4)\)-equivariant continuous function \[ \kappa_{g}: C^{2}(\mathbb{R}^{4};\mathbb{R}^{4}) \rightarrow \Sigma^2(\mathbb{R}^{4}; \mathbb{R}). \]
By assumption, we would have $ \kappa_{g}( \mu_{13,+}(\frac{1}{2\delta}) ) {\xrightarrow{\text{\,\,deg\,\,}} } \kappa( \mu_{7}(0) ) $
where:
$$
\kappa_{g}(\mu_{13,+}(\frac{1}{2\delta})) = -\frac{1}{2}\,{\frac {  \left( 4\,{\delta}^{2}-3-8\,g \right) }{{
\delta}^{2}}} e_{1}^\ast \otimes e_{1}^\ast
$$
and
$$
\kappa_{g}( \mu_{7}(0) ) = -2e_{1}^\ast \otimes e_{1}^\ast
$$

Setting $g=\frac{1}{2}{\delta}^{2}- \frac{3}{8}$, we obtain a contradiction, since $\kappa_{g}( \mu_{13,+}(\frac{1}{2\delta} ) ) $ is zero bilinear form while $\kappa_{g}(  \mu_{7}(0) ) $ is non-zero negative semi-definite bilinear form.

\item
$
(\mathfrak{ d}^{'}_{4,\delta},J_{k})
\centernot
{ \xrightarrow{\text{\,\,deg\,\,}} }
(\mathfrak{ d}_{4,1/2}, J_1)
$\newline

In fact, this holds because $
\mathfrak{ d}^{'}_{4,\delta}
\centernot
{ \xrightarrow{\text{\,\,deg\,\,}} }
\mathfrak{ d}_{4,1/2}
$.  Let us now explore this in detail. Suppose for the sake of contradiction that $
\mathfrak{ d}^{'}_{4,\delta}
{ \xrightarrow{\text{\,\,deg\,\,}} }
\mathfrak{ d}_{4,1/2}
$. First, consider the \(\operatorname{GL}(\mathbb{R}^4)\)-equivariant continuous function \[ \kappa_{g}: C^{2}(\mathbb{R}^{4};\mathbb{R}^{4}) \rightarrow \Sigma^2(\mathbb{R}^{4}; \mathbb{R}). \]
By assumption, we would have $ \kappa_{g}(  \mu_{13,+}(\frac{1}{2\delta})  ) {\xrightarrow{\text{\,\,deg\,\,}} } \kappa( \mu_{10}(2) ) $
where:
$$
\kappa_{g}(\mu_{13,+}(\frac{1}{2\delta})) = -\frac{1}{2}\,{\frac { \left( 4\,{\delta}^{2}-3-8\,g \right) }{{
\delta}^{2}}}
 e_{1}^\ast \otimes e_{1}^\ast
$$
and
$$
\kappa_{g}( \mu_{10}(2) ) = (6+16g) e_{1}^\ast \otimes e_{1}^\ast
$$

Setting $g= \frac{1}{2}\,{\delta}^{2}- \frac{3}{8} $, we obtain a contradiction, since $\kappa_{g}( \mu_{13,+}(\frac{1}{2\delta}) ) $ is zero bilinear form while
$\kappa_{g}( \mu_{10}(2) ) = 8\delta^2  e_{1}^\ast \otimes e_{1}^\ast$ is non-zero positive semi-definite bilinear form.

\item\label{delta4primar4a}
$
(\mathfrak{ d}^{'}_{4,\delta},J_{k})
\centernot
{ \xrightarrow{\text{\,\,deg\,\,}} }
(\mathfrak{ r}_{4,a,a },J)
$ with $ a \neq 0$  \newline

In fact, this holds because $
\mathfrak{ d}^{'}_{4,\delta}
\centernot
{ \xrightarrow{\text{\,\,deg\,\,}} }
\mathfrak{ r}_{4,a,a }
$.  Let us now explore this in detail. Suppose for the sake of contradiction that $
\mathfrak{ d}^{'}_{4,\delta}
%\centernot
{ \xrightarrow{\text{\,\,deg\,\,}} }
\mathfrak{ r}_{4,a,a }
$. We have two cases:

\begin{itemize}

  \item If $1+2\,{a}^{2}+ \frac{1}{8} \left( 4\,{\delta}^{2}-3 \right)  \left( 2\,a+1  \right) ^{2} \neq 0$.

First, consider the \(\operatorname{GL}(\mathbb{R}^4)\)-equivariant continuous function \[ \kappa_{g}: C^{2}(\mathbb{R}^{4};\mathbb{R}^{4}) \rightarrow \Sigma^2(\mathbb{R}^{4}; \mathbb{R}). \]
By assumption, we would have $ \kappa_{g}(  \mu_{13,+}(\frac{1}{2\delta})   ) {\xrightarrow{\text{\,\,deg\,\,}} }
\kappa( \mu_{6}(a,0)  ) $
where:
$$
\kappa_{g}( \mu_{13,+}(\frac{1}{2\delta}) ) = -\frac{1}{2}\,{\frac {   \left( 4\,{\delta}^{2}-3-8\,g \right) }{{\delta}^{2}}} e_1^{\ast} \otimes e_1^{\ast}
$$
and
$$
\kappa_{g}( \mu_{6}(a,0)   ) =  (1+2\,{a}^{2}+g \left( 2\,a+1 \right) ^{2}) e_1^{\ast} \otimes e_1^{\ast}
$$

Setting $g=  \frac{1}{2} {\delta}^{2}- \frac{3}{8} $, we obtain a contradiction, since $\kappa_{g}( \mu_{13,+}(\frac{1}{2\delta}) ) $ is the zero bilinear form while $\kappa_{g}( \mu_{6}(a,0)  ) = 1+2\,{a}^{2}+ \frac{1}{8}\, \left( 4\,{\delta}^{2}-3 \right)  \left( 2\,a+1 \right) ^{2} e_1^{\ast} \otimes e_1^{\ast} $ is non-zero bilinear form.

  \item If $1+2\,{a}^{2}+ \frac{1}{8} \left( 4\,{\delta}^{2}-3 \right)  \left( 2\,a+1  \right) ^{2} = \frac{1}{2}\, \left( 2\,a+1 \right) ^{2}{\delta}^{2}+ \frac{1}{8}\, \left( 2\,a-1
 \right)  \left( 2\,a-5 \right) =0$, and so
 \begin{eqnarray}\label{eq1}
\frac{9}{4}\, \left( 2\,a+1 \right) ^{3}{\delta}^{2}+{\frac {9}{16}}\, \left(
2\,a-1 \right)  \left( 4\,{a}^{2}-8\,a-5 \right)
&=&0.
 \end{eqnarray}

We consider the \(\operatorname{GL}(\mathbb{R}^4)\)-equivariant continuous function \[ \Psi_{f,g}: C^{2}(\mathbb{R}^{4};\mathbb{R}^{4}) \rightarrow \Sigma^{3}(\mathbb{R}^{4};\mathbb{R}). \]
Therefore $$ \Psi_{f,g}(  \mu_{13,+}(\frac{1}{2\delta} ) ) \xrightarrow{\text{\,\,deg\,\,}} \Psi_{f,g}(  \mu_{6}(a,0)   ), $$
where:
$$
 \Psi_{f,g}( \mu_{13,+}(\frac{1}{2\delta}) ) =  -\frac{1}{2}\,{\frac {  \left( -16\,g-18\,f+24\,f{\delta}^{2}+36\,{\delta}^{2}-15 \right) }{{\delta}^{3}}} e_1^{\ast} \otimes e_1^{\ast} \otimes e_1^{\ast}
$$
and
$$
\Psi_{f,g}(  \mu_{6}(a,0)  )= (\left( 8\,g+12\,f+12 \right) {a}^{3}+ \left( 12\,g+6\,f \right) {a}^{
2}+ \left( 6\,g+6\,f \right) a+6+g+3\,f) e_1^{\ast} \otimes e_1^{\ast} \otimes e_1^{\ast}
$$
Taking $f=0$  and $g= -{\frac {15}{16}}+ \frac{9}{4}\,{\delta}^{2} $, we obtain a contradiction, since $ \Psi_{f,g}(\mu_{13,+}(1/2\delta)) \equiv 0$
and $ \Psi_{f,g}(  \mu_{6}(a,0)  ) = \frac{9}{4}\, \left( 2\,a+1 \right) ^{3}{\delta}^{2}+{\frac {9}{16}}\, \left(
2\,a-1 \right)  \left( 4\,{a}^{2}-8\,a-9 \right)
 e_{1}^\ast \otimes e_{1}^\ast \otimes e_{1}^\ast \not\equiv 0$ because of Equation \ref{eq1} and $\delta \neq 0$.

\end{itemize}

\item
$
(\mathfrak{ d}^{'}_{4,\delta},J_{k})
\centernot
{ \xrightarrow{\text{\,\,deg\,\,}} }
(\mathfrak{r}_{4,\alpha,1},J)
$ with $\alpha \neq0$ and $-1 < \alpha \leq 1$  \newline
See Subsection \ref{delta4prima} items \ref{delta4primar4a} and using $\mathfrak{r}_{4,\alpha,1} \equiv \mathfrak{r}_{4,\frac{1}{\alpha},\frac{1}{\alpha}}$.

\item
$
(\mathfrak{ d}^{'}_{4,\delta},J_{k})
\centernot
{ \xrightarrow{\text{\,\,deg\,\,}} }
(\mathfrak{ r}_{4,1},J)
$\newline

In fact, this holds because $
\mathfrak{ d}^{'}_{4,\delta}
\centernot
{ \xrightarrow{\text{\,\,deg\,\,}} }
\mathfrak{ r}_{4,1}
$.  Let us now explore this in detail. Suppose for the sake of contradiction that $
\mathfrak{ d}^{'}_{4,\delta}
{ \xrightarrow{\text{\,\,deg\,\,}} }
\mathfrak{ r}_{4,1}
$. We have two cases:

\begin{itemize}

\item If $\delta \neq \frac{1}{6}\,\sqrt {3}$

Consider the \(\operatorname{GL}(\mathbb{R}^4)\)-equivariant continuous function \[ \kappa_{g}: C^{2}(\mathbb{R}^{4};\mathbb{R}^{4}) \rightarrow \Sigma^2(\mathbb{R}^{4}; \mathbb{R}). \]
By assumption, we would have $ \kappa_{g}( \mu_{13,+}(\frac{1}{2\delta})  ) {\xrightarrow{\text{\,\,deg\,\,}} } \kappa( \mu_{9} ) $
where:
$$
\kappa_{g}( \mu_{13,+}(\frac{1}{2\delta}) ) = -\frac{1}{2}\, \frac{ \left( 4\,{\delta}^{2}-3-8\,g \right) } {{\delta}^{2}} e_1^{\ast} \otimes e_1^{\ast}
$$
and
$$
\kappa_{g}( \mu_9 ) = 3+9\,g e_1^{\ast} \otimes e_1^{\ast}
$$

Setting $g= \frac{1}{2}\,{\delta}^{2}-\frac{3}{8}$, we obtain a contradiction, since $\kappa_{g}( \mu_{13,+}(\frac{1}{2\delta}) ) $ is the zero bilinear form while $\kappa_{g}( \mu_9  ) = -\frac{3}{8}+\frac{9}{2}\,{\delta}^{2} e_1^{\ast} \otimes e_1^{\ast}$ is non-zero bilinear form.

 \item If $\delta = \frac{1}{6}\,\sqrt {3}$

We consider the \(\operatorname{GL}(\mathbb{R}^4)\)-equivariant continuous function \[ \Psi_{f,g}: C^{2}(\mathbb{R}^{4};\mathbb{R}^{4}) \rightarrow \Sigma^{3}(\mathbb{R}^{4};\mathbb{R}). \]
Therefore $$ \Psi_{f,g}(  \mu_{13,+}(\frac{1}{2\delta} ) ) \xrightarrow{\text{\,\,deg\,\,}} \Psi_{f,g}(  \mu_{9}   ), $$
where:
$$
 \Psi_{f,g}( \mu_{13,+}(\frac{1}{2\delta}) ) =   -\frac{1}{2}\,{\frac {  \left( -16\,g-18\,f+24\,f{\delta}^{2}+36\,{\delta}^{2}-15 \right) }{{\delta}^{3}}} e_1^{\ast} \otimes e_1^{\ast} \otimes e_1^{\ast}
$$
and
$$
\Psi_{f,g}(  \mu_{9}  )= (27 g+27f+18) e_1^{\ast} \otimes e_1^{\ast} \otimes e_1^{\ast}
$$
Taking $f= 0 $  and $g= -\frac{3}{4} $, we obtain a contradiction, since $ \Psi_{f,g}( \mu_{13,+}(\frac{1}{2\delta}) ) \equiv 0$
and $ \Psi_{f,g}(  \mu_{9}  ) = -\frac{9}{4}  e_1^{\ast} \otimes e_1^{\ast} \otimes e_1^{\ast}$.

\end{itemize}

\item
$
(\mathfrak{ d}^{'}_{4,\delta},J_{k})
\centernot
{ \xrightarrow{\text{\,\,deg\,\,}} }
(\mathfrak{rr}_{3,1 }\times \mathbb{R},J)
$\newline

In fact, this holds because $
\mathfrak{ d}^{'}_{4,\delta}
\centernot
{ \xrightarrow{\text{\,\,deg\,\,}} }
\mathfrak{rr}_{3,1 }\times \mathbb{R}
$.  Let us now explore this in detail. Suppose for the sake of contradiction that $
\mathfrak{ d}^{'}_{4,\delta}
{ \xrightarrow{\text{\,\,deg\,\,}} }
\mathfrak{rr}_{3,1 }\times \mathbb{R}
$. First, consider the \(\operatorname{GL}(\mathbb{R}^4)\)-equivariant continuous function \[ \kappa_{g}: C^{2}(\mathbb{R}^{4};\mathbb{R}^{4}) \rightarrow \Sigma^2(\mathbb{R}^{4}; \mathbb{R}). \]
By assumption, we would have $ \kappa_{g}( \mu_{13,+}(\frac{1}{2\delta}) ) {\xrightarrow{\text{\,\,deg\,\,}} } \kappa( \mu_{3}(0,0) ) $
where:
$$
\kappa_{g}( \mu_{13,+}(\frac{1}{2\delta}) ) = -\frac{1}{2}\, \frac{ \left( 4\,{\delta}^{2}-3-8\,g \right) } {{\delta}^{2}} e_1^{\ast} \otimes e_1^{\ast}
$$
and
$$
\kappa_{g}( \mu_{3}(0,0) ) = (2+4g) e_1^{\ast} \otimes e_1^{\ast}
$$

Setting $g= \frac{1}{2}\,{\delta}^{2}-\frac{3}{8}$, we obtain a contradiction, since $\kappa_{g}( \mu_{13,+}(\frac{1}{2\delta})  ) $ is the zero bilinear form while $\kappa_{g}(  \mu_{3}(0,0) ) = (\frac{1}{2}+2\,{\delta}^{2}) e_1^{\ast} \otimes e_1^{\ast}$ is non-zero positive semi-definite bilinear form.

\item
$
(\mathfrak{ d}^{'}_{4,\delta},J_{k})
\centernot
{ \xrightarrow{\text{\,\,deg\,\,}} }
(\mathfrak{ r}_{4,-1/2,-1/2},J)
$\newline
See Subsection \ref{delta4prima} items \ref{delta4primar4a}.

\item
$
(\mathfrak{ d}^{'}_{4,\delta},J_{k})
\centernot
{ \xrightarrow{\text{\,\,deg\,\,}} }
(\mathfrak{r}_2\times \mathbb{R}^{2},J)
$\newline

In fact, this holds because $
\mathfrak{ d}^{'}_{4,\delta}
\centernot
{ \xrightarrow{\text{\,\,deg\,\,}} }
\mathfrak{r}_2\times \mathbb{R}^{2}
$.  Let us now explore this in detail. Suppose for the sake of contradiction that $
\mathfrak{ d}^{'}_{4,\delta}
{ \xrightarrow{\text{\,\,deg\,\,}} }
\mathfrak{r}_2\times \mathbb{R}^{2}
$. First, consider the \(\operatorname{GL}(\mathbb{R}^4)\)-equivariant continuous function \[ \kappa_{g}: C^{2}(\mathbb{R}^{4};\mathbb{R}^{4}) \rightarrow \Sigma^2(\mathbb{R}^{4}; \mathbb{R}). \]
By assumption, we would have $ \kappa_{g}( \mu_{13,+}(\frac{1}{2\delta}) ) {\xrightarrow{\text{\,\,deg\,\,}} } \kappa( \mu_{6}(0,0) ) $
where:
$$
\kappa_{g}( \mu_{13,+}(\frac{1}{2\delta}) ) = -\frac{1}{2}\, \frac{ \left( 4\,{\delta}^{2}-3-8\,g \right) } {{\delta}^{2}} e_1^{\ast} \otimes e_1^{\ast}
$$
and
$$
\kappa_{g}( \mu_{6}(0,0) ) = (1+g) e_1^{\ast} \otimes e_1^{\ast}
$$

Setting $g= \frac{1}{2}\,{\delta}^{2}-\frac{3}{8}$, we obtain a contradiction, since $\kappa_{g}( \mu_{13,+}(\frac{1}{2\delta}) ) $ is the zero bilinear form while $\kappa_{g}( \mu_{6}(0,0) ) = (\frac{5}{8}+\frac{1}{2}\,{\delta}^{2})e_1^{\ast} \otimes e_1^{\ast} $ is non-zero positive semi-definite bilinear form.

  \item
$
(\mathfrak{ d}^{'}_{4,\delta},J_{k})
%\centernot
{ \xrightarrow{\text{\,\,deg\,\,}} }
(\mathfrak{h}_{3}\times \mathbb{R},J)
$\newline

Consider the curve $g_t = \left(
          \begin{array}{cc}
            A(t) & -B(t) \\
            B(t) & A(t) \\
          \end{array}
        \right)$ where $$A(t)= 0
                             \qquad
                             \mbox{and}
                             \qquad
                         B(t)= \left(
                                 \begin{array}{cc}
                                   0 & -e^{\frac{1}{1}t} \\
                                   e^t & 0 \\
                                 \end{array}
                               \right)
                         .$$
Note that $g_t \in \operatorname{GL}(\mathbb{R}^4, {J_{\tiny{\mbox{cn}}}})$ and its action on $\mu_{13, (-1)^k}((-1)^{ \lfloor \frac{k+1}{2} \rfloor} \frac{1}{2\delta} )$ gives
$$
\left\{
\begin{array}{l}
[{\it e1},{\it e3}]={\it e2},
[{\it e1},{\it e4}]=\frac{1}{2}\,{\frac { \left( -1 \right) ^{  \lfloor \frac{k+1}{2} \rfloor }{{\rm e}^{-t}}}{\delta}}{\it e1}+ \left( -1 \right) ^{k}{{\rm e}^{-t}}{\it e3},\\
{[{\it e2},{\it e4}]}={\frac { \left( -1 \right)^{ \lfloor \frac{k+1}{2} \rfloor }{{\rm e}^{-t}}}{\delta}}{\it e2},
[{\it e3},{\it e4}]= - \left( -1 \right) ^{k}{{\rm e}^{-t}}{\it e1}+ \frac{1}{2}\,{\frac { \left( -1 \right) ^{ \lfloor \frac{k+1}{2}\rfloor }{{\rm e}^{-t}} }{\delta}}{\it e3}
\end{array}
\right.
$$
Taking the limit $g_t \cdot \mu_{13, (-1)^k}((-1)^{ \lfloor \frac{k+1}{2} \rfloor} \frac{1}{2\delta} )$ as $t \to \infty$, we obtain $\mu_{1}$.

  \item
$
(\mathfrak{ d}^{'}_{4,\delta},J_{k})
\centernot
{ \xrightarrow{\text{\,\,deg\,\,}} }
(\mathfrak{ r}_{4,1,1},J)
$\newline
See Subsection \ref{delta4prima} items \ref{delta4primar4a}.

\end{enumerate}

\subsection{$(\mathfrak{ d}_{4,\lambda}, J_{k_{2}})$ with $\lambda>1/2$, $\lambda \neq 1$ and ${k_{2}}=1,2$}\label{delta4normal}
\begin{enumerate}

\item
$
(\mathfrak{ d}_{4,\lambda}, J_{k_{2}})
\centernot
{ \xrightarrow{\text{\,\,deg\,\,}} }
(\mathfrak{d}_4,J_2)
$\newline

Because $
\mathfrak{ d}_{4,\lambda}
\centernot
{ \xrightarrow{\text{\,\,deg\,\,}} }
\mathfrak{d}_4
$, since $\operatorname{Dim}\operatorname{Der}( \mathfrak{ d}_{4,\lambda} )=5$ and $\operatorname{Dim}\operatorname{Der}(\mathfrak{d}_4)=5$, and $\mathfrak{ d}_{4,\lambda} \not \equiv \mathfrak{d}_4$.

  \item\label{delta4normalr4prima}
$
(\mathfrak{ d}_{4,\lambda}, J_{k_{2}})
\centernot
{ \xrightarrow{\text{\,\,deg\,\,}} }
(\mathfrak{r}^{'}_{4,c,d},J_{k_{2}})
$ $d>0$ \newline

In fact, this holds because $
\mathfrak{ d}_{4,\lambda}
\centernot
{ \xrightarrow{\text{\,\,deg\,\,}} }
\mathfrak{r}^{'}_{4,c,d}
$.  Let us now explore this in detail. Suppose for the sake of contradiction that $
\mathfrak{ d}_{4,\lambda}
{ \xrightarrow{\text{\,\,deg\,\,}} }
\mathfrak{r}^{'}_{4,c,d}
$. We have two cases:

\begin{itemize}
  \item If $-\frac{1}{2}\, \left( 2\,c+1 \right) ^{2}{\lambda}^{2} + \frac{1}{2}\, \left( 2\,c+1  \right) ^{2}\lambda+\frac{1}{2}-2\,{d}^{2}-2\,c \neq 0$

First, consider the \(\operatorname{GL}(\mathbb{R}^4)\)-equivariant continuous function \[ \kappa_{g}: C^{2}(\mathbb{R}^{4};\mathbb{R}^{4}) \rightarrow \Sigma^2(\mathbb{R}^{4}; \mathbb{R}). \]
By assumption, we would have $ \kappa_{g}( \mu_{10}(\frac{1}{\lambda}) ) {\xrightarrow{\text{\,\,deg\,\,}} } \kappa( \mu_{6}(c,d) ) $
where:
$$
\kappa_{g}( \mu_{10}(\frac{1}{\lambda}) ) =2\,{\frac {{\lambda}^{2}+1-\lambda+2\,g}{{\lambda}^{2}}} e_1^{\ast} \otimes e_1^{\ast}
$$
and
$$
\kappa_{g}( \mu_{6}(c,d) ) = (1-2\,{d}^{2}+2\,{c}^{2}+g+4\,gc+4\,g{c}^{2} )e_1^{\ast} \otimes e_1^{\ast}
$$
Setting $g= -\frac{1}{2}\,{\lambda}^{2}-\frac{1}{2}+\frac{1}{2}\,\lambda $, we obtain a contradiction, since $\kappa_{g}(\mu_{10}(\frac{1}{\lambda})   ) $ is zero bilinear form while $\kappa_{g}(\mu_{6}(c,d)  ) = (-\frac{1}{2}\, \left( 2\,c+1 \right) ^{2}{\lambda}^{2} + \frac{1}{2}\, \left( 2\,c+1  \right) ^{2}\lambda+\frac{1}{2}-2\,{d}^{2}-2\,c) e_1^{\ast} \otimes e_1^{\ast}$ is non-zero semi-definite bilinear form.

\item If $-\frac{1}{2}\, \left( 2\,c+1 \right) ^{2}{\lambda}^{2} + \frac{1}{2}\, \left( 2\,c+1  \right) ^{2}\lambda+\frac{1}{2}-2\,{d}^{2}-2\,c =  0$ and so
\begin{eqnarray}\label{eq2}
  -\frac{9}{4} \, \left( 2\,c+1 \right) ^{3}{\lambda}^{2}+\frac{9}{4}\, \left( 2\,c+1
 \right) ^{3}\lambda +\frac{9}{4} -18\,{c}^{2}-\frac{9}{2}\,c-18\,c{d}^{2}-9\,{d}^{2}  = 0
\end{eqnarray}

We consider the \(\operatorname{GL}(\mathbb{R}^4)\)-equivariant continuous function \[ \Psi_{f,g}: C^{2}(\mathbb{R}^{4};\mathbb{R}^{4}) \rightarrow \Sigma^{3}(\mathbb{R}^{4};\mathbb{R}). \]
Therefore $$ \Psi_{0,g}( \mu_{10}(\frac{1}{\lambda})  ) \xrightarrow{\text{\,\,deg\,\,}} \Psi_{0,g}(  \mu_{6}(c,d)  ), $$
where:
$$
 \Psi_{0,g}( \mu_{10}(\frac{1}{\lambda}) ) = 2\,{\frac {  \left( 4\,g+6-9\,\lambda+9\,{\lambda} ^{2} \right) }{{\lambda}^{3}}} e_1^{\ast} \otimes e_1^{\ast}\otimes e_1^{\ast}
$$
and
$$
\Psi_{0,g}( \mu_{6}(c,d)  )= (8\,g{c}^{3}+12\,g{c}^{2}+6\,gc+g+6+12\,{c}^{3}-36\,c{d}^{2}) e_1^{\ast} \otimes e_1^{\ast}\otimes e_1^{\ast}
$$
Taking  $g= -\frac{3}{2}+\frac{9}{4}\,\lambda- \frac{9}{4}\,{\lambda}^{2} $ we obtain  $ \Psi_{0,g}( \mu_{10}(\frac{1}{\lambda}) ) \equiv 0$
and $ \Psi_{0,g}(   \mu_{6}(c,d ) ) = P(\lambda,c,d)  e_1^{\ast} \otimes e_1^{\ast}\otimes e_1^{\ast} $ with
\begin{eqnarray}\label{eq3}
 P(\lambda,c,d)& = & -\frac{9}{4}\, \left( 2\,c+1 \right) ^{3}{\lambda}^{2}+ \frac{9}{4}\, \left( 2\,c+1  \right) ^{3}\lambda+ \frac{9}{2}-18\,{c}^{2}-9\,c -36\,c{d}^{2}.
\end{eqnarray}
Note that $P(\lambda,c,d) \neq 0$ because, on the contrary, by equating Equations \ref{eq2} and \ref{eq3},  we would obtain $c= \frac{1}{2}$, and therefore $\frac{1}{2}\, \left( 2\,\lambda-1 \right) ^{2}+2\,{d}^{2}=0 $ which can not be since $\lambda \neq \frac{1}{2}$.

Thus, we obtain a contradiction in this case as well.

\end{itemize}

\item
$
(\mathfrak{ d}_{4,\lambda}, J_{k_{2}})
\centernot
{ \xrightarrow{\text{\,\,deg\,\,}} }
(\mathfrak{rr}^{'}_{3,a}\times\mathbb{R},J_{k_{2}})
$ with $a>0$ and $k_2=1,2$ \newline

In fact, this holds because $
\mathfrak{ d}_{4,\lambda}
\centernot
{ \xrightarrow{\text{\,\,deg\,\,}} }
\mathfrak{rr}^{'}_{3,a}\times\mathbb{R}
$.  Let us now explore this in detail. Suppose for the sake of contradiction that $
\mathfrak{ d}_{4,\lambda}
{ \xrightarrow{\text{\,\,deg\,\,}} }
\mathfrak{rr}^{'}_{3,a}\times\mathbb{R}
$. We have two cases:

\begin{itemize}
\item  If $-2-2\,{a}^{2}{\lambda}^{2}+2\,{a}^{2}\lambda \neq 0$

First, consider the \(\operatorname{GL}(\mathbb{R}^4)\)-equivariant continuous function \[ \kappa_{g}: C^{2}(\mathbb{R}^{4};\mathbb{R}^{4}) \rightarrow \Sigma^2(\mathbb{R}^{4}; \mathbb{R}). \]
By assumption, we would have $ \kappa_{g}( \mu_{10}(\frac{1}{\lambda}) ) {\xrightarrow{\text{\,\,deg\,\,}} } \kappa( \mu_{7}(a) ) $
where:
$$
\kappa_{g}(  \mu_{10}(\frac{1}{\lambda})  ) = 2\,{\frac {{\lambda}^{2}+1-\lambda+2\,g}{{\lambda}^{2}}} e_{1}^{\ast} \otimes e_{1}^{\ast}
$$
and
$$
\kappa_{g}( \mu_{7}(a) ) = (2\,{a}^{2}-2+4\,g{a}^{2}) e_{1}^{\ast} \otimes e_{1}^{\ast}
$$

Setting $g= -\frac{1}{2}\,{\lambda}^{2}-\frac{1}{2}+\frac{1}{2}\,\lambda $, we obtain a contradiction, since $\kappa_{g}( \mu_{10}(\frac{1}{\lambda}) ) $ is zero bilinear form while $\kappa_{g}( \mu_{7}(a) ) = (-2-2\,{a}^{2}{\lambda}^{2}+2\,{a}^{2}\lambda)e_{1}^{\ast} \otimes e_{1}^{\ast} $ is non-zero bilinear form.

  \item If $-2-2\,{a}^{2}{\lambda}^{2}+2\,{a}^{2}\lambda = 0$, and so
\begin{eqnarray}\label{eq4}
18\,{a}^{3}\lambda  - 18\,{a}^{3}{\lambda}^{2} -18\,a   &=& 0
\end{eqnarray}

We consider the \(\operatorname{GL}(\mathbb{R}^4)\)-equivariant continuous function \[ \Psi_{f,g}: C^{2}(\mathbb{R}^{4};\mathbb{R}^{4}) \rightarrow \Sigma^{3}(\mathbb{R}^{4};\mathbb{R}). \]
Therefore $$ \Psi_{0,g}( \mu_{10}(\frac{1}{\lambda})  ) \xrightarrow{\text{\,\,deg\,\,}} \Psi_{0,g}(  \mu_{7}(a)  ), $$
where:
$$
 \Psi_{0,g}( \mu_{10}(\frac{1}{\lambda}) ) = 2\,{\frac {4\,g+6-9\,\lambda+9\,{\lambda}^{2}}{{\lambda}^{3}}} e_{1}^{\ast} \otimes e_{1}^{\ast}\otimes e_{1}^{\ast}
$$
and
$$
\Psi_{0,g}( \mu_{7}(a)  )= (8\,{a}^{3}g-36\,a+12\,{a}^{3}) e_{1}^{\ast} \otimes e_{1}^{\ast}\otimes e_{1}^{\ast}
$$
Setting $g= -\frac{3}{2} + \frac{9}{4}\,\lambda-\frac{9}{4}\,{\lambda}^{2}  $, we obtain a contradiction, since $ \Psi_{0,g}( \mu_{10}(\frac{1}{\lambda}) ) \equiv 0$
and $ \Psi_{0,g}(  \mu_{7}(a) ) = P(\lambda,a) e_{1}^{\ast} \otimes e_{1}^{\ast}\otimes e_{1}^{\ast} $ with
\begin{eqnarray*}
P(\lambda,a) & = &  18\,{a}^{3}\lambda-18\,{a}^{3}{\lambda}^{2}-36\,a\\
             & = &  -18a
\end{eqnarray*}
since Equation (\ref{eq4}).

\end{itemize}

  \item
$
(\mathfrak{ d}_{4,\lambda}, J_{k_{2}})
\centernot
{ \xrightarrow{\text{\,\,deg\,\,}} }
(\mathfrak{r}^{'}_{4,-1/2,d},J_{k_{2}})
$\newline
See Subsection \ref{delta4normal} item \ref{delta4normalr4prima}

\item
$
(\mathfrak{ d}_{4,\lambda}, J_{k_{2}})
\centernot
{ \xrightarrow{\text{\,\,deg\,\,}} }
(\mathfrak{rr}^{'}_{3,0}\times\mathbb{R},J)
$\newline

In fact, this holds because $
\mathfrak{ d}_{4,\lambda}
\centernot
{ \xrightarrow{\text{\,\,deg\,\,}} }
\mathfrak{rr}^{'}_{3,0}\times\mathbb{R}
$.  Let us now explore this in detail. Suppose for the sake of contradiction that $
\mathfrak{ d}_{4,\lambda}
{ \xrightarrow{\text{\,\,deg\,\,}} }
\mathfrak{rr}^{'}_{3,0}\times\mathbb{R}
$. First, consider the \(\operatorname{GL}(\mathbb{R}^4)\)-equivariant continuous function \[ \kappa_{g}: C^{2}(\mathbb{R}^{4};\mathbb{R}^{4}) \rightarrow \Sigma^2(\mathbb{R}^{4}; \mathbb{R}). \]
By assumption, we would have $ \kappa_{g}( \mu_{10}(\frac{1}{\lambda})  ) {\xrightarrow{\text{\,\,deg\,\,}} } \kappa( \mu_{7}(0) ) $
where:
$$
\kappa_{g}( \mu_{10}(\frac{1}{\lambda}) ) = 2\,{\frac {{\lambda}^{2}+1-\lambda+2\,g}{{\lambda}^{2}}} e_{1}^{\ast} \otimes e_{1}^{\ast}
$$
and
$$
\kappa_{g}( \mu_{7}(0) ) = -2 e_{1}^{\ast} \otimes e_{1}^{\ast}
$$

Setting $g= -\frac{1}{2}\,{\lambda}^{2}-\frac{1}{2}+\frac{1}{2}\,\lambda $, we obtain a contradiction, since $\kappa_{g}(  \mu_{10}(\frac{1}{\lambda})  ) $ is zero bilinear form while $\kappa_{g}(  \mu_{7}(0)  ) $ is non-zero negative semi-definite bilinear form.

\item
$
(\mathfrak{ d}_{4,\lambda}, J_{k_{2}})
\centernot
{ \xrightarrow{\text{\,\,deg\,\,}} }
(\mathfrak{ d}_{4,1/2}, J_1)
$\newline

In fact, this holds because $
\mathfrak{ d}_{4,\lambda}
\centernot
{ \xrightarrow{\text{\,\,deg\,\,}} }
\mathfrak{ d}_{4,1/2}
$.  Let us now explore this in detail. Suppose for the sake of contradiction that $
\mathfrak{ d}_{4,\lambda}
{ \xrightarrow{\text{\,\,deg\,\,}} }
\mathfrak{ d}_{4,1/2}
$. First, consider the \(\operatorname{GL}(\mathbb{R}^4)\)-equivariant continuous function \[ \kappa_{g}: C^{2}(\mathbb{R}^{4};\mathbb{R}^{4}) \rightarrow \Sigma^2(\mathbb{R}^{4}; \mathbb{R}). \]
By assumption, we would have $ \kappa_{g}( \mu_{10}(\frac{1}{\lambda}) ) {\xrightarrow{\text{\,\,deg\,\,}} } \kappa( \mu_{10}(2) ) $
where:
$$
\kappa_{g}( \mu_{10}(\frac{1}{\lambda}) ) = 2\,{\frac {{\lambda}^{2}+1-\lambda+2\,g}{{\lambda}^{2}}} e_{1}^{\ast} \otimes e_{1}^{\ast}
$$
and
$$
\kappa_{g}( \mu_{10}(2) ) = (6 + 16g) e_{1}^{\ast} \otimes e_{1}^{\ast}
$$

Setting $g= -\frac{1}{2}\,{\lambda}^{2}-\frac{1}{2}+\frac{1}{2}\,\lambda $, we obtain a contradiction, since $\kappa_{g}( \mu_{10}(\frac{1}{\lambda}) ) $ is zero bilinear form while $\kappa_{g}(  \mu_{10}(2) ) = -2\, \left( 2\,\lambda-1 \right) ^{2} e_{1}^{\ast} \otimes e_{1}^{\ast} $ is non-zero negative semi-definite bilinear form (since $\lambda \neq \frac{1}{2}$).

\item\label{delta4normalr4normal}
$
(\mathfrak{ d}_{4,\lambda}, J_{k_{2}})
\centernot
{ \xrightarrow{\text{\,\,deg\,\,}} }
(\mathfrak{ r}_{4,\alpha,\alpha },J)
$ with $a\neq 0$. \newline

In fact, this holds because $
\mathfrak{ d}_{4,\lambda}
\centernot
{ \xrightarrow{\text{\,\,deg\,\,}} }
\mathfrak{ r}_{4,\alpha,\alpha }
$.  Let us now explore this in detail. Suppose for the sake of contradiction that $
\mathfrak{ d}_{4,\lambda}
{ \xrightarrow{\text{\,\,deg\,\,}} }
\mathfrak{ r}_{4,\alpha,\alpha }
$. We have two cases:

\begin{itemize}

\item If $-\frac{1}{2}\, \left( 2\,a+1 \right) ^{2}{\lambda}^{2} + \frac{1}{2}\, \left( 2\,a+1  \right) ^{2}\lambda + \frac{1}{2}-2\,a \neq 0$

First, consider the \(\operatorname{GL}(\mathbb{R}^4)\)-equivariant continuous function \[ \kappa_{g}: C^{2}(\mathbb{R}^{4};\mathbb{R}^{4}) \rightarrow \Sigma^2(\mathbb{R}^{4}; \mathbb{R}). \]
By assumption, we would have $ \kappa_{g}( \mu_{10}(\frac{1}{\lambda}) ) {\xrightarrow{\text{\,\,deg\,\,}} } \kappa( \mu_{6}(a,0) ) $
where:
$$
\kappa_{g}( \mu_{10}(\frac{1}{\lambda}) ) =  2\,{\frac {{\lambda}^{2}+1-\lambda+2\,g}{{\lambda}^{2}}} e_{1}^{\ast} \otimes e_{1}^{\ast}
$$
and
$$
\kappa_{g}( \mu_{6}(a,0)  ) = (1+2\,{a}^{2}+g+4\,ga+4\,g{a}^{2}) e_{1}^{\ast} \otimes e_{1}^{\ast}
$$

Setting $g= -\frac{1}{2}\,{\lambda}^{2}-\frac{1}{2}+\frac{1}{2}\,\lambda $, we obtain a contradiction, since $\kappa_{g}( \mu_{10}(\frac{1}{\lambda}) ) $ is zero bilinear form while $\kappa_{g}( \mu_{6}(a,0)  ) =
(-\frac{1}{2}\, \left( 2\,a+1 \right) ^{2}{\lambda}^{2} + \frac{1}{2}\, \left( 2\,a+1  \right) ^{2}\lambda + \frac{1}{2}-2\,a )e_{1}^{\ast} \otimes e_{1}^{\ast}$ is non-zero bilinear form.

\item If $-\frac{1}{2}\, \left( 2\,a+1 \right) ^{2}{\lambda}^{2} + \frac{1}{2}\, \left( 2\,a+1  \right) ^{2}\lambda + \frac{1}{2}-2\,a = 0$, then we have
\begin{eqnarray}\label{eq5}
 -\frac{9}{4}\, \left( 2\,a+1 \right) ^{3}{\lambda}^{2}+\frac{9}{4}\, \left( 2\,a+1 \right) ^{3}\lambda-\frac{9}{2}\,a-18\,{a}^{2}+\frac{9}{4} = 0
\end{eqnarray}

We consider the \(\operatorname{GL}(\mathbb{R}^4)\)-equivariant continuous function \[ \Psi_{f,g}: C^{2}(\mathbb{R}^{4};\mathbb{R}^{4}) \rightarrow \Sigma^{3}(\mathbb{R}^{4};\mathbb{R}). \]
Therefore $$ \Psi_{0,g}( \mu_{10}(\frac{1}{\lambda})  ) \xrightarrow{\text{\,\,deg\,\,}} \Psi_{0,g}( \mu_{6}(a,0)  ), $$
where:
$$
 \Psi_{0,g}(  \mu_{10}(\frac{1}{\lambda}) ) = 2\,{\frac {4\,g+6-9\,\lambda+9\,{\lambda}^{2}}{{\lambda}^{3}}} e_{1}^{\ast} \otimes e_{1}^{\ast}\otimes e_{1}^{\ast}
$$
and
$$
\Psi_{0,g}(  \mu_{6}(a,0) )=  (8\,g{a}^{3}+12\,g{a}^{2}+6\,ga+g+6+12\,{a}^{3}) e_{1}^{\ast} \otimes e_{1}^{\ast}\otimes e_{1}^{\ast}
$$
Setting $g= -\frac{3}{2} + \frac{9}{4}\,\lambda-\frac{9}{4}\,{\lambda}^{2}  $ we obtain $ \Psi_{0,g}(  \mu_{10}(\frac{1}{\lambda}) ) \equiv 0$
and $ \Psi_{0,g}(  \mu_{6}(a,0)  ) = P(\lambda,a) e_{1}^{\ast} \otimes e_{1}^{\ast}\otimes e_{1}^{\ast} $ where
\begin{eqnarray*}
  P(\lambda,a) & = &( -\frac{9}{4}\, \left( 2\,a+1 \right) ^{3}{\lambda}^{2} + \frac{9}{4}\, \left( 2\,a+1  \right) ^{3}\lambda + \frac{9}{2}-9\,a-18\,{a}^{2} ) \\
               & = & -\frac{9}{2}\,a+\frac{9}{4}
\end{eqnarray*}
as shown in Equation Equation \ref{eq5}. Note that $P(\lambda,a) \neq 0$, since otherwise we would have $\lambda = \frac{1}{2}$ from Equation (\ref{eq5}) again.

\end{itemize}

\item
$
(\mathfrak{ d}_{4,\lambda}, J_{k_{2}})
\centernot
{ \xrightarrow{\text{\,\,deg\,\,}} }
(\mathfrak{r}_{4,\alpha,1},J)
$\newline
See Subsection \ref{delta4normal} items \ref{delta4normalr4normal} and using $\mathfrak{r}_{4,\alpha,1} \equiv \mathfrak{r}_{4,\frac{1}{\alpha},\frac{1}{\alpha}}$.

\item
$
(\mathfrak{ d}_{4,\lambda}, J_{k_{2}})
\centernot
{ \xrightarrow{\text{\,\,deg\,\,}} }
(\mathfrak{ r}_{4,1},J)
$\newline

In fact, this holds because $
\mathfrak{ d}_{4,\lambda}
\centernot
{ \xrightarrow{\text{\,\,deg\,\,}} }
\mathfrak{ r}_{4,1}
$.  Let us now explore this in detail. Suppose for the sake of contradiction that $
\mathfrak{ d}_{4,\lambda}
{ \xrightarrow{\text{\,\,deg\,\,}} }
\mathfrak{ r}_{4,1}
$. First, consider the \(\operatorname{GL}(\mathbb{R}^4)\)-equivariant continuous function \[ \kappa_{g}: C^{2}(\mathbb{R}^{4};\mathbb{R}^{4}) \rightarrow \Sigma^2(\mathbb{R}^{4}; \mathbb{R}). \]
By assumption, we would have $ \kappa_{g}( \mu_{10}(\frac{1}{\lambda}) ) {\xrightarrow{\text{\,\,deg\,\,}} } \kappa( \mu_{9} ) $
where:
$$
\kappa_{g}( \mu_{10}(\frac{1}{\lambda}) ) = 2\,{\frac {{\lambda}^{2}+1-\lambda+2\,g}{{\lambda}^{2}}} e_{1}^{\ast} \otimes e_{1}^{\ast}
$$
and
$$
\kappa_{g}( \mu_{9} ) =  (3+9g) e_{1}^{\ast} \otimes e_{1}^{\ast}
$$

Setting $g= -\frac{1}{2}\,{\lambda}^{2}-\frac{1}{2}+\frac{1}{2}\,\lambda $, we obtain a contradiction, since $\kappa_{g}(  \mu_{10}(\frac{1}{\lambda}) ) $ is zero bilinear form while $\kappa_{g}( \mu_{9} ) = -\frac{9}{2}( (\lambda - \frac{1}{2})^2 + \frac{1}{12})  $ is non-zero negative semi-definite bilinear form.

\item
$
(\mathfrak{ d}_{4,\lambda}, J_{k_{2}})
\centernot
{ \xrightarrow{\text{\,\,deg\,\,}} }
(\mathfrak{rr}_{3,1 }\times \mathbb{R},J)
$\newline

In fact, this holds because $
\mathfrak{ d}_{4,\lambda}
\centernot
{ \xrightarrow{\text{\,\,deg\,\,}} }
\mathfrak{rr}_{3,1 }\times \mathbb{R}
$.  Let us now explore this in detail. Suppose for the sake of contradiction that $
\mathfrak{ d}_{4,\lambda}
{ \xrightarrow{\text{\,\,deg\,\,}} }
\mathfrak{rr}_{3,1 }\times \mathbb{R}
$. First, consider the \(\operatorname{GL}(\mathbb{R}^4)\)-equivariant continuous function \[ \kappa_{g}: C^{2}(\mathbb{R}^{4};\mathbb{R}^{4}) \rightarrow \Sigma^2(\mathbb{R}^{4}; \mathbb{R}). \]
By assumption, we would have $ \kappa_{g}( \mu_{10}(\frac{1}{\lambda}) ) {\xrightarrow{\text{\,\,deg\,\,}} } \kappa( \mu_{3}(0,0) ) $
where:
$$
\kappa_{g}( \mu_{10}(\frac{1}{\lambda}) ) = 2\,{\frac {{\lambda}^{2}+1-\lambda+2\,g}{{\lambda}^{2}}} e_{1}^{\ast} \otimes e_{1}^{\ast}
$$
and
$$
\kappa_{g}( \mu_{3}(0,0) ) = (2 + 4g) e_{1}^{\ast} \otimes e_{1}^{\ast}
$$

Setting $g= -\frac{1}{2}\,{\lambda}^{2}-\frac{1}{2}+\frac{1}{2}\,\lambda $, we obtain a contradiction, since $\kappa_{g}(  \mu_{10}(\frac{1}{\lambda})  ) $ is zero bilinear form while $\kappa_{g}( \mu_{3}(0,0)  ) = -2\,\lambda\, \left( \lambda-1 \right) e_{1}^{\ast} \otimes e_{1}^{\ast} $ is non-zero  bilinear form (since $\lambda \neq 1$).

\item
$
(\mathfrak{ d}_{4,\lambda}, J_{k_{2}})
\centernot
{ \xrightarrow{\text{\,\,deg\,\,}} }
(\mathfrak{ r}_{4,-1/2,-1/2},J)
$\newline
See Subsection \ref{delta4normal} items \ref{delta4normalr4normal}

\item
$
(\mathfrak{ d}_{4,\lambda}, J_{k_{2}})
\centernot
{ \xrightarrow{\text{\,\,deg\,\,}} }
(\mathfrak{r}_2\times \mathbb{R}^{2},J)
$\newline

In fact, this holds because
$
\mathfrak{ d}_{4,\lambda}
\centernot
{ \xrightarrow{\text{\,\,deg\,\,}} }
\mathfrak{r}_2\times \mathbb{R}^{2}
$
.  Let us now explore this in detail. Suppose for the sake of contradiction that $
\mathfrak{ d}_{4,\lambda}
{ \xrightarrow{\text{\,\,deg\,\,}} }
\mathfrak{r}_2\times \mathbb{R}^{2}
$. We have two cases:

\begin{itemize}

  \item If ${\lambda}^{2}-\lambda -1 \neq 0$

First, consider the \(\operatorname{GL}(\mathbb{R}^4)\)-equivariant continuous function \[ \kappa_{g}: C^{2}(\mathbb{R}^{4};\mathbb{R}^{4}) \rightarrow \Sigma^2(\mathbb{R}^{4}; \mathbb{R}). \]
By assumption, we would have $ \kappa_{g}( \mu_{10}(\frac{1}{\lambda}) ) {\xrightarrow{\text{\,\,deg\,\,}} } \kappa( \mu_{6}(0,0) ) $
where:
$$
\kappa_{g}( \mu_{10}(\frac{1}{\lambda}) ) = 2\,{\frac {{\lambda}^{2}+1-\lambda+2\,g}{{\lambda}^{2}}} e_{1}^{\ast} \otimes e_{1}^{\ast}
$$
and
$$
\kappa_{g}( \mu_{6}(0,0) ) = (1+g) e_{1}^{\ast} \otimes e_{1}^{\ast}
$$

Setting $g= -\frac{1}{2}\,{\lambda}^{2}-\frac{1}{2}+\frac{1}{2}\,\lambda $, we obtain a contradiction, since $\kappa_{g}( \mu_{10}(\frac{1}{\lambda})  ) $ is zero bilinear form while $\kappa_{g}( \mu_{6}(0,0) ) = -\frac{1}{2}( {\lambda}^{2}-\lambda -1)$ is non-zero bilinear form.

  \item If ${\lambda}^{2}-\lambda -1 = 0$

We consider the \(\operatorname{GL}(\mathbb{R}^4)\)-equivariant continuous function \[ \Psi_{f,g}: C^{2}(\mathbb{R}^{4};\mathbb{R}^{4}) \rightarrow \Sigma^{3}(\mathbb{R}^{4};\mathbb{R}). \]
Therefore $$ \Psi_{0,g}( \mu_{10}(\frac{1}{\lambda})  ) \xrightarrow{\text{\,\,deg\,\,}} \Psi_{0,g}(  \mu_{6}(0,0)  ), $$
where:
$$
 \Psi_{0,g}( \mu_{10}(\frac{1}{\lambda}) ) = 2\,{\frac {4\,g+6-9\,\lambda+9\,{\lambda}^{2}}{{\lambda}^{3}}} e_{1}^{\ast} \otimes e_{1}^{\ast}\otimes e_{1}^{\ast}
$$
and
$$
\Psi_{0,g}(  \mu_{6}(0,0)  )= (6+g) e_{1}^{\ast} \otimes e_{1}^{\ast}\otimes e_{1}^{\ast}
$$
Setting $g= -\frac{3}{2} + \frac{9}{4}\,\lambda-\frac{9}{4}\,{\lambda}^{2}  $, we obtain a contradiction, since $ \Psi_{0,g}( \mu_{10}(\frac{1}{\lambda})  ) \equiv 0$
and $ \Psi_{0,g}(  \mu_{6}(0,0) ) = -\frac{9}{4}({\lambda}^{2}-\lambda-2)e_{1}^{\ast} \otimes e_{1}^{\ast}\otimes e_{1}^{\ast} \neq 0 $ (by hypothesis).

\end{itemize}

  \item
$
(\mathfrak{ d}_{4,\lambda}, J_{k_{2}})
%\centernot
{ \xrightarrow{\text{\,\,deg\,\,}} }
(\mathfrak{h}_{3}\times \mathbb{R},J)
$\newline

Let $a= \frac{1}{\lambda}$ if $k=1$, and let $a= \frac{1}{1-\lambda}$ if $k=2$. Consider the curve $g_t = \left(
          \begin{array}{cc}
            A(t) & -B(t) \\
            B(t) & A(t) \\
          \end{array}
        \right)$ where $$A(t)= \operatorname{diag}(e^t,1)
                             \qquad
                             \mbox{and}
                             \qquad
                         B(t)= \left(
                                 \begin{array}{cc}
                                   0 & 0 \\
                                   \frac{e^{2t}}{a-2} & 0 \\
                                 \end{array}
                               \right)
                         .$$
Note that $g_t$ is well-defined since $\lambda\neq \frac{1}{2}$ and $g_t $ is in $\operatorname{GL}(\mathbb{R}^4, {J_{\tiny{\mbox{cn}}}})$. The action of $g_t $ on $\mu_{10}(a)$ gives
$$
\left\{
\begin{array}{l}
[{\it e1},{\it e2}]={{\rm e}^{-t}} \left( a-1 \right) {\it e2},
[{\it e1},{\it e3}]={\it e2}+{{\rm e}^{-t}}{\it e3},
[{\it e1},{\it e4}]={{\rm e}^{-t}}a{\it e4},[{\it e2},{\it e3}]={{\rm e}^{-t}}{\it e4}
\end{array}
\right.
$$
Taking the limit $g_t \cdot \mu_{10}(a)$ as $t \to \infty$, we obtain $\mu_{1}$.

  \item
$
(\mathfrak{ d}_{4,\lambda}, J_{k_{2}})
\centernot
{ \xrightarrow{\text{\,\,deg\,\,}} }
(\mathfrak{ r}_{4,1,1},J)
$\newline
See Subsection \ref{delta4normal} items \ref{delta4normalr4normal}.

\end{enumerate}

\subsection{$(\mathfrak{ h}_{ 4},J)$}\label{h4}
\begin{enumerate}

\item
$
(\mathfrak{ h}_{ 4},J)
\centernot
{ \xrightarrow{\text{\,\,deg\,\,}} }
(\mathfrak{d}_4,J_2)
$\newline

Because $
\mathfrak{ h}_{ 4}
\centernot
{ \xrightarrow{\text{\,\,deg\,\,}} }
\mathfrak{d}_4
$, since $\operatorname{Dim}\operatorname{Der}( \mathfrak{ h}_{ 4} )=5$ and $\operatorname{Dim}\operatorname{Der}(\mathfrak{d}_4)=5$.

  \item\label{h4r4prima}
$
(\mathfrak{ h}_{ 4},J)
\centernot
{ \xrightarrow{\text{\,\,deg\,\,}} }
(\mathfrak{r}^{'}_{4,c,d},J_{k_{2}})
$ with $d>0$ \newline

In fact, this holds because $
\mathfrak{ h}_{ 4}
\centernot
{ \xrightarrow{\text{\,\,deg\,\,}} }
\mathfrak{r}^{'}_{4,c,d}
$ .  Let us now explore this in detail. Suppose for the sake of contradiction that $
\mathfrak{ h}_{ 4}
{ \xrightarrow{\text{\,\,deg\,\,}} }
\mathfrak{r}^{'}_{4,c,d}
$. We have two cases:

\begin{itemize}
\item If $4\,{c}^{2}-12\,c+5-16\,{d}^{2} \neq 0$

First, consider the \(\operatorname{GL}(\mathbb{R}^4)\)-equivariant continuous function \[ \kappa_{g}: C^{2}(\mathbb{R}^{4};\mathbb{R}^{4}) \rightarrow \Sigma^2(\mathbb{R}^{4}; \mathbb{R}). \]
By assumption, we would have $ \kappa_{g}( \mu_{14} ) {\xrightarrow{\text{\,\,deg\,\,}} } \kappa( \mu_{6}(c,d) ) $
where:
$$
\kappa_{g}( \mu_{14} ) = (6+16g) e_{1}^{\ast} \otimes e_{1}^{\ast}
$$
and
$$
\kappa_{g}( \mu_{6}(c,d) ) = (1-2\,{d}^{2}+2\,{c}^{2}+g+4\,gc+4\,g{c}^{2}) e_{1}^{\ast} \otimes e_{1}^{\ast}
$$

Setting $g= -\frac{3}{8}$, we obtain a contradiction, since $\kappa_{g}( \mu_{14} ) $ is zero bilinear form while $\kappa_{g}(  \mu_{6}(c,d)  )  = \frac{1}{8}(4\,{c}^{2}-12\,c+5-16\,{d}^{2})$ is non-zero bilinear form.

\item If $4\,{c}^{2}-12\,c+5-16\,{d}^{2} = 0$, and so

\begin{eqnarray}\label{eq6}
  9\,{c}^{3}-27\,{c}^{2}+{\frac {45}{4}}\,c-36\,c{d}^{2}=0
\end{eqnarray}

We consider the \(\operatorname{GL}(\mathbb{R}^4)\)-equivariant continuous function \[ \Psi_{f,g}: C^{2}(\mathbb{R}^{4};\mathbb{R}^{4}) \rightarrow \Sigma^{3}(\mathbb{R}^{4};\mathbb{R}). \]
Therefore $$ \Psi_{0,g}( \mu_{14} ) \xrightarrow{\text{\,\,deg\,\,}} \Psi_{0,g}( \mu_{6}(c,d) ), $$
where:
$$
 \Psi_{0,g}( \mu_{14}  ) = (64 g+60)  e_{1}^{\ast} \otimes e_{1}^{\ast} \otimes e_{1}^{\ast}
$$
and
$$
\Psi_{0,g}(  \mu_{6}(c,d)  )= (8\,g{c}^{3}+12\,g{c}^{2}+6\,gc+g+6+12\,{c}^{3}-36\,c{d}^{2})  e_{1}^{\ast} \otimes e_{1}^{\ast} \otimes e_{1}^{\ast}
$$
Setting $g= -\frac{15}{16} $ we obtain  $ \Psi_{0,g}( \mu_{14}) \equiv 0$
and $ \Psi_{0,g}(  \mu_{6}(c,d) ) = P(c,d)  e_{1}^{\ast} \otimes e_{1}^{\ast} \otimes e_{1}^{\ast}$ with
\begin{eqnarray*}
  P(c,d) & = &  \frac{9}{2}\,{c}^{3}-{\frac {45}{4}}\,{c}^{2}-{\frac {45}{8}}\,c+{\frac {81}{16}}-36\,c{d}^{2}\\
         & = &  (\frac{9}{2}\,{c}^{3}-{\frac {45}{4}}\,{c}^{2}-{\frac {45}{8}}\,c+{\frac {81}{16}})-(9\,{c}^{3}-27\,{c}^{2}+{\frac {45}{4}}\,c)\\
         & = &-{\frac {9}{16}}\, \left( 2\,c-1 \right)  \left( -3+2\,c \right) ^{2}.
\end{eqnarray*}
as shown in Equation \ref{eq6}. Note that $P(c,d) \neq 0$, since otherwise we would have  $ c= \frac{1}{2} $ or $c= \frac{3}{2}$ which would imply, from Equation (\ref{eq5}), that $d=0$ or $d= \pm \frac{1}{2} \sqrt{-1}$ respectively; both of which are impossible.

\end{itemize}

\item
$
(\mathfrak{ h}_{ 4},J)
\centernot
{ \xrightarrow{\text{\,\,deg\,\,}} }
(\mathfrak{rr}^{'}_{3,a}\times\mathbb{R},J_{k_{2}})
$ with $a>0$ and $k_2=1,2$ \newline

In fact, this holds because $
\mathfrak{ h}_{ 4}
\centernot
{ \xrightarrow{\text{\,\,deg\,\,}} }
\mathfrak{rr}^{'}_{3,a}\times\mathbb{R}$ .  Let us now explore this in detail. Suppose for the sake of contradiction that $
\mathfrak{ h}_{ 4}
{ \xrightarrow{\text{\,\,deg\,\,}} }
\mathfrak{rr}^{'}_{3,a}\times\mathbb{R}$. We have two cases:

\begin{itemize}
  \item If $a \neq 2$.

First, consider the \(\operatorname{GL}(\mathbb{R}^4)\)-equivariant continuous function \[ \kappa_{g}: C^{2}(\mathbb{R}^{4};\mathbb{R}^{4}) \rightarrow \Sigma^2(\mathbb{R}^{4}; \mathbb{R}). \]
By assumption, we would have $ \kappa_{g}( \mu_{14}  ) {\xrightarrow{\text{\,\,deg\,\,}} } \kappa( \mu_{7}(a) ) $
where:
$$
\kappa_{g}(  \mu_{14}  ) = (6+16g) e_{1}^{\ast} \otimes  e_{1}^{\ast}
$$
and
$$
\kappa_{g}( \mu_{7}(a) ) = (2\,{a}^{2}-2+4\,g{a}^{2}) e_{1}^{\ast} \otimes  e_{1}^{\ast}
$$

Setting $g= -\frac{3}{8}$, we obtain a contradiction, since $\kappa_{g}(   \mu_{14}  ) $ is zero bilinear form while $\kappa_{g}( \mu_{7}(a) ) = (\frac{1}{2}\,{a}^{2}-2)e_{1}^{\ast} \otimes  e_{1}^{\ast} $ is non-zero bilinear form.

  \item If $a=2$

We consider the \(\operatorname{GL}(\mathbb{R}^4)\)-equivariant continuous function \[ \Psi_{f,g}: C^{2}(\mathbb{R}^{4};\mathbb{R}^{4}) \rightarrow \Sigma^{3}(\mathbb{R}^{4};\mathbb{R}). \]
Therefore $$ \Psi_{0,g}(  \mu_{14}  ) \xrightarrow{\text{\,\,deg\,\,}} \Psi_{0,g}(   \mu_{7}(a)  ), $$
where:
$$
 \Psi_{0,g}(  \mu_{14}  ) = (64g+60)e_{1}^{\ast} \otimes  e_{1}^{\ast} \otimes  e_{1}^{\ast}
$$
and
$$
\Psi_{0,g}(   \mu_{7}(a)  ) = (8\,{a}^{3}g-36\,a+12\,{a}^{3})e_{1}^{\ast} \otimes  e_{1}^{\ast} \otimes  e_{1}^{\ast} = (64\,g+24) e_{1}^{\ast} \otimes  e_{1}^{\ast} \otimes  e_{1}^{\ast}
$$
Setting $g= -\frac{15}{16} $, we obtain a contradiction, since $ \Psi_{0,g}( \mu_{14} ) \equiv 0$
and $ \Psi_{0,g}( \mu_{7}(a) ) = -36 e_{1}^{\ast} \otimes  e_{1}^{\ast} \otimes  e_{1}^{\ast} \neq 0$.

\end{itemize}

  \item
$
(\mathfrak{ h}_{ 4},J)
\centernot
{ \xrightarrow{\text{\,\,deg\,\,}} }
(\mathfrak{r}^{'}_{4,-1/2,d},J_{k_{2}})
$\newline
See Subsection \ref{h4} items \ref{h4r4prima}

\item
$
(\mathfrak{ h}_{ 4},J)
\centernot
{ \xrightarrow{\text{\,\,deg\,\,}} }
(\mathfrak{rr}^{'}_{3,0}\times\mathbb{R},J)
$\newline

In fact, this holds because $
\mathfrak{ h}_{ 4}
\centernot
{ \xrightarrow{\text{\,\,deg\,\,}} }
\mathfrak{rr}^{'}_{3,0}\times\mathbb{R}
$.  Let us now explore this in detail. Suppose for the sake of contradiction that $
\mathfrak{ h}_{ 4}
{ \xrightarrow{\text{\,\,deg\,\,}} }
\mathfrak{rr}^{'}_{3,0}\times\mathbb{R}
$. First, consider the \(\operatorname{GL}(\mathbb{R}^4)\)-equivariant continuous function \[ \kappa_{g}: C^{2}(\mathbb{R}^{4};\mathbb{R}^{4}) \rightarrow \Sigma^2(\mathbb{R}^{4}; \mathbb{R}). \]
By assumption, we would have $ \kappa_{g}(  \mu_{14}  ) {\xrightarrow{\text{\,\,deg\,\,}} } \kappa( \mu_{7}(0) ) $
where:
$$
\kappa_{g}( \mu_{14} ) = \kappa_{g}(  \mu_{14}  ) = (6+16g) e_{1}^{\ast} \otimes  e_{1}^{\ast}
$$
and
$$
\kappa_{g}( \mu_{7}(0) ) = -2 e_{1}^{\ast} \otimes  e_{1}^{\ast}
$$

Setting $g= -\frac{3}{8}$, we obtain a contradiction, since $\kappa_{g}( \mu_{14}  ) $ is zero bilinear form while $\kappa_{g}( \mu_{7}(0) ) $ is non-zero negative semi-definite bilinear form.

\item
$
(\mathfrak{ h}_{ 4},J)
%\centernot
{ \xrightarrow{\text{\,\,deg\,\,}} }
(\mathfrak{ d}_{4,1/2}, J_1)
$\newline

Consider the curve $g_t = \left(
          \begin{array}{cc}
            A(t) & -B(t) \\
            B(t) & A(t) \\
          \end{array}
        \right)$ where $$A(t)= \operatorname{diag}(1, e^{-t})
                             \qquad
                             \mbox{and}
                             \qquad
                         B(t)= 0
                         .$$
Note that $g_t \in \operatorname{GL}(\mathbb{R}^4, {J_{\tiny{\mbox{cn}}}})$ and its action on $\mu_{14}$ gives
$$
\left\{
\begin{array}{l}
[{\it e1},{\it e2}]={\it e2},
[{\it e1},{\it e3}]={{\rm e}^{-t}}{\it e2}+{\it e3},
[{\it e1},{\it e4}]=2\,{\it e4},[{\it e2},{\it e3}]={\it e4}
\end{array}
\right.
$$
Taking the limit $g_t \cdot \mu_{14}$ as $t \to \infty$, we obtain $\mu_{10}(2)$.

\item\label{h4r4normal}
$
(\mathfrak{ h}_{ 4},J)
\centernot
{ \xrightarrow{\text{\,\,deg\,\,}} }
(\mathfrak{ r}_{4, a, a },J)
$ with $a \neq 0$ and $ a \neq \frac{1}{2}$ \newline

In fact, this holds because $
\mathfrak{ h}_{ 4}
\centernot
{ \xrightarrow{\text{\,\,deg\,\,}} }
\mathfrak{ r}_{4, a, a }$ .  Let us now explore this in detail. Suppose for the sake of contradiction that $
\mathfrak{ h}_{ 4}
{ \xrightarrow{\text{\,\,deg\,\,}} }
\mathfrak{ r}_{4, a, a }$. We have two cases:

\begin{itemize}
  \item If $ 2a-5 \neq 0$

  First, consider the \(\operatorname{GL}(\mathbb{R}^4)\)-equivariant continuous function \[ \kappa_{g}: C^{2}(\mathbb{R}^{4};\mathbb{R}^{4}) \rightarrow \Sigma^2(\mathbb{R}^{4}; \mathbb{R}). \]
By assumption, we would have $ \kappa_{g}( \mu_{14} ) {\xrightarrow{\text{\,\,deg\,\,}} } \kappa( \mu_{6}(a,0) ) $
where:
$$
\kappa_{g}( \mu_{14}  ) = (6+16g)e_{1}^{\ast} \otimes e_{1}^{\ast}
$$
and
$$
\kappa_{g}( \mu_{6}(a,0) ) = (1+2\,{a}^{2}+g+4\,ga+4\,g{a}^{2}) e_{1}^{\ast} \otimes e_{1}^{\ast}
$$

Setting $g= -\frac{3}{8}$, we obtain a contradiction, since $\kappa_{g}( \mu_{14}  ) $ is the zero bilinear form while $\kappa_{g}( \mu_{6}(a,0) ) = \frac{1}{8}\, \left( 2\,a-1 \right)  \left( 2\,a-5 \right) e_{1}^{\ast} \otimes e_{1}^{\ast}$ is non-zero bilinear form (since $a \neq \frac{1}{2}$).

\item If $ 2a-5 = 0$

We consider the \(\operatorname{GL}(\mathbb{R}^4)\)-equivariant continuous function \[ \Psi_{f,g}: C^{2}(\mathbb{R}^{4};\mathbb{R}^{4}) \rightarrow \Sigma^{3}(\mathbb{R}^{4};\mathbb{R}). \]
Therefore $$ \Psi_{0,g}( \mu_{14}   ) \xrightarrow{\text{\,\,deg\,\,}} \Psi_{0,g}(  \mu_{6}(a,0)  ), $$
where:
$$
 \Psi_{0,g}( \mu_{14} ) = (64g+60) e_{1}^{\ast} \otimes e_{1}^{\ast} \otimes e_{1}^{\ast}
$$
and
$$
\Psi_{0,g}( \mu_{6}(a,0)   ) = (216g + \frac{387}{2})e_{1}^{\ast} \otimes e_{1}^{\ast} \otimes e_{1}^{\ast}
$$
Setting $g= -\frac{15}{16} $, we obtain a contradiction, since  $ \Psi_{0,g}( \mu_{14}  ) \equiv 0$ and $ \Psi_{0,g}(  \mu_{6}(a,0) ) = -9  e_{1}^{\ast} \otimes e_{1}^{\ast} \otimes e_{1}^{\ast}$ .

\end{itemize}

\item\label{h4r4normalunmedio}
$
(\mathfrak{ h}_{ 4},J)
\centernot
{ \xrightarrow{\text{\,\,deg\,\,}} }
(\mathfrak{ r}_{4, \frac{1}{2}, \frac{1}{2} },J)
$  \newline

First, we consider the $\operatorname{GL}(\mathbb{R}^4, {J_{\tiny{\mbox{cn}}}})$-equivariant continuous function $\psi_{1,0}: C^{2}_{J_{\tiny{\mbox{cn}}}}(\mathbb{R}^{4};\mathbb{R}^{4}) \rightarrow C^{2}(\mathbb{R}^{4};\mathbb{R}^{4})$. Therefore $$\mathfrak{A}:=(\mathbb{R}^{4}, \psi_{1,0}( \mu_{14} )) { \xrightarrow{\text{\,\,deg\,\,}} } \mathfrak{B}:=(\mathbb{R}^{4}, \psi_{1,0}( \mu_{6}(\frac{1}{2},0) )),$$
where:
$$
\mathfrak{A} =
\left\{
\begin{array}{l}
[{\it e1},{\it e2}]={\it e2}+{\it e4},
[{\it e1},{\it e3}]=-{\it e1}+{\it e2}+{\it e3}+{\it e4},
[{\it e1},{\it e4}]=-2\,{\it e2}+2\,{\it e4},\\
{[{\it e2},{\it e3}]}=-{\it e2}+{\it e4}
\end{array}
\right.
$$
and
$$
\mathfrak{B} =
\left\{
\begin{array}{l}
[{\it e1},{\it e2}]=\frac{1}{2}\,{\it e2}+\frac{1}{2}\,{\it e4},
[{\it e1},{\it e3}]=-{\it e1}+{\it e3},
[{\it e1},{\it e4}]=-\frac{1}{2}\,{\it e2}+\frac{1}{2}\,{\it e4}
\end{array}
\right.
$$

Consider the $\operatorname{GL}(\mathbb{R}^4)$-equivariant continuous function $$\kappa_g: C^{2}(\mathbb{R}^{4};\mathbb{R}^{4}) \rightarrow \Sigma^2(\mathbb{R}^{4}; \mathbb{R}).$$
By assumption, we would have $ \kappa_{g}(\mathfrak{A} ) {\xrightarrow{\text{\,\,deg\,\,}} } \kappa( \mathfrak{B} ) $.
The matrix of $\kappa_{g}(\mathfrak{A} )$ with respect to the canonical basis is:$$
\left[ \begin {array}{cccc} 2+16\,g&0&4+8\,g&0\\ \noalign{\medskip}0&0
&0&0\\ \noalign{\medskip}4+8\,g&0&2+4\,g&0\\ \noalign{\medskip}0&0&0&0
\end {array} \right]
$$
while the matrix of $\kappa_{g}( \mathfrak{B} )$ with respect to the canonical basis is:$$
 \left[ \begin {array}{cccc} 1+4\,g&0&1+2\,g&0\\ \noalign{\medskip}0&0
&0&0\\ \noalign{\medskip}1+2\,g&0&1+g&0\\ \noalign{\medskip}0&0&0&0
\end {array} \right]
$$
Setting $g= -\frac{1}{2}$, we obtain a contradiction, since the dimension of the radical of $\kappa_{-\frac{1}{2}}(\mathfrak{A})$ is $3$  and  the dimension of the radical of $\kappa_{-\frac{1}{2}}(\mathfrak{B})$ is $2$.

\item
$
(\mathfrak{ h}_{ 4},J)
\centernot
{ \xrightarrow{\text{\,\,deg\,\,}} }
(\mathfrak{r}_{4,\alpha,1},J)
$ with $-1 < \alpha \leq 1$, $\alpha \neq 0$ \newline
See Subsection \ref{h4} items \ref{h4r4prima} and using $\mathfrak{r}_{4,\alpha,1} \equiv \mathfrak{r}_{4,\frac{1}{\alpha},\frac{1}{\alpha}}$.

\item
$
(\mathfrak{ h}_{ 4},J)
\centernot
{ \xrightarrow{\text{\,\,deg\,\,}} }
(\mathfrak{ r}_{4,1},J)
$\newline

In fact, this holds because $
\mathfrak{ h}_{ 4}
\centernot
{ \xrightarrow{\text{\,\,deg\,\,}} }
\mathfrak{ r}_{4,1}
$.  Let us now explore this in detail. Suppose for the sake of contradiction that $
\mathfrak{ h}_{ 4}
{ \xrightarrow{\text{\,\,deg\,\,}} }
\mathfrak{ r}_{4,1}
$. First, we consider the $\operatorname{GL}(\mathbb{R}^4)$-equivariant continuous function $$\kappa_g: C^{2}(\mathbb{R}^{4};\mathbb{R}^{4}) \rightarrow \Sigma^2(\mathbb{R}^{4}; \mathbb{R}).$$
By assumption, we would have $ \kappa_{g}( \mu_{14} ) {\xrightarrow{\text{\,\,deg\,\,}} } \kappa( \mu_{9} ) $
where:
$$
\kappa_{g}( \mu_{14} ) = (6+16g) e_1^{\ast} \otimes e_1^{\ast}
$$
and
$$
\kappa_{g}( \mu_{9} ) = (3+9g)e_1^{\ast} \otimes e_1^{\ast}
$$

Setting $g= -\frac{3}{8}$, we obtain a contradiction, since $\kappa_{g}( \mu_{14} ) $ is zero bilinear form while $\kappa_{g}( \mu_{9} ) $ is non-zero negative semi-definite bilinear form.

\item
$
(\mathfrak{ h}_{ 4},J)
\centernot
{ \xrightarrow{\text{\,\,deg\,\,}} }
(\mathfrak{rr}_{3,1 }\times \mathbb{R},J)
$\newline

In fact, this holds because $
\mathfrak{ h}_{ 4}
\centernot
{ \xrightarrow{\text{\,\,deg\,\,}} }
\mathfrak{rr}_{3,1 }\times \mathbb{R}
$.  Let us now explore this in detail. Suppose for the sake of contradiction that $
\mathfrak{ h}_{ 4}
{ \xrightarrow{\text{\,\,deg\,\,}} }
\mathfrak{rr}_{3,1 }\times \mathbb{R}
$. First, consider the \(\operatorname{GL}(\mathbb{R}^4)\)-equivariant continuous function \[ \kappa_{g}: C^{2}(\mathbb{R}^{4};\mathbb{R}^{4}) \rightarrow \Sigma^2(\mathbb{R}^{4}; \mathbb{R}). \]
By assumption, we would have $ \kappa_{g}( \mu_{14} ) {\xrightarrow{\text{\,\,deg\,\,}} } \kappa( \mu_{3}(0,0) ) $
where:
$$
\kappa_{g}( \mu_{14} ) = (6+16g) e_1^{\ast} \otimes e_1^{\ast}
$$
and
$$
\kappa_{g}( \mu_{9} ) = ( 2+4g)e_1^{\ast} \otimes e_1^{\ast}
$$

Setting $g= -\frac{3}{8}$, we obtain a contradiction, since $\kappa_{g}( \mu_{14} ) $ is zero bilinear form while $\kappa_{g}( \mu_{3}(0,0) ) $ is non-zero positive semi-definite bilinear form.

\item
$
(\mathfrak{ h}_{ 4},J)
\centernot
{ \xrightarrow{\text{\,\,deg\,\,}} }
(\mathfrak{ r}_{4,-1/2,-1/2},J)
$\newline
See Subsection \ref{h4} items \ref{h4r4normal} and \ref{h4r4normalunmedio}.

\item
$
(\mathfrak{ h}_{ 4},J)
\centernot
{ \xrightarrow{\text{\,\,deg\,\,}} }
(\mathfrak{r}_2\times \mathbb{R}^{2},J)
$\newline

In fact, this holds because $
\mathfrak{ h}_{ 4}
\centernot
{ \xrightarrow{\text{\,\,deg\,\,}} }
\mathfrak{r}_2\times \mathbb{R}^{2}
$.  Let us now explore this in detail. Suppose for the sake of contradiction that $
\mathfrak{ h}_{ 4}
{ \xrightarrow{\text{\,\,deg\,\,}} }
\mathfrak{r}_2\times \mathbb{R}^{2}
$. First, consider the \(\operatorname{GL}(\mathbb{R}^4)\)-equivariant continuous function \[ \kappa_{g}: C^{2}(\mathbb{R}^{4};\mathbb{R}^{4}) \rightarrow \Sigma^2(\mathbb{R}^{4}; \mathbb{R}). \]
By assumption, we would have $ \kappa_{g}( \mu_{14} ) {\xrightarrow{\text{\,\,deg\,\,}} } \kappa( \mu_{6}(0,0) ) $
where:
$$
\kappa_{g}( \mu_{14} ) = (6+16g) e_1^{\ast} \otimes e_1^{\ast}
$$
and
$$
\kappa_{g}( \mu_{6}(0,0) ) = (1+g  )e_1^{\ast} \otimes e_1^{\ast}
$$

Setting $g= -\frac{3}{8}$, we obtain a contradiction, since $\kappa_{g}( \mu_{14} ) $ is zero bilinear form while $\kappa_{g}( \mu_{6}(0,0)) $ is non-zero positive semi-definite bilinear form.

  \item
$
(\mathfrak{ h}_{ 4},J)
%\centernot
{ \xrightarrow{\text{\,\,deg\,\,}} }
(\mathfrak{h}_{3}\times \mathbb{R},J)
$\newline

Consider the curve $g_t = \left(
          \begin{array}{cc}
            A(t) & -B(t) \\
            B(t) & A(t) \\
          \end{array}
        \right)$ where $$A(t)= \operatorname{diag}(e^t,e^{2t})
                             \qquad
                             \mbox{and}
                             \qquad
                         B(t)= 0
                         .$$
Note that $g_t \in \operatorname{GL}(\mathbb{R}^4, {J_{\tiny{\mbox{cn}}}})$ and its action on $\mu_{14}$ gives
$$
\left\{
\begin{array}{l}
[{\it e1},{\it e2}]={\frac {{1}}{{{\rm e}^{t}}}}\it e2,
[{\it e1},{\it e3}]={\it e2}+{\frac {1}{{{\rm e}^{t}}}}{\it e3},
[{\it e1},{\it e4}]={\frac {2}{{{\rm e}^{t}}}}{\it e4},
[{\it e2},{\it e3}]={\frac {1}{{{\rm e}^{t}}}}{\it e4}
\end{array}
\right.
$$
Taking the limit $g_t \cdot \mu_{14}$ as $t \to \infty$, we obtain $\mu_{1}$.

  \item
$
(\mathfrak{ h}_{ 4},J)
\centernot
{ \xrightarrow{\text{\,\,deg\,\,}} }
(\mathfrak{ r}_{4,1,1},J)
$\newline

See Subsection \ref{h4} items \ref{h4r4normal} and \ref{h4r4normalunmedio}.

\end{enumerate}

\subsection{$(\mathfrak{ d}_{4,1/2}, J_2)$}\label{delta4unmedio}
\begin{enumerate}

\item
$
(\mathfrak{ d}_{4,1/2}, J_2)
\centernot
{ \xrightarrow{\text{\,\,deg\,\,}} }
(\mathfrak{ d}_{4,1/2}, J_1)
$\newline

Suppose for the sake of contradiction that $
(\mathfrak{ d}_{4,1/2}, J_2)
%\centernot
{ \xrightarrow{\text{\,\,deg\,\,}} }
(\mathfrak{ d}_{4,1/2}, J_1)
$. First, we consider the $\operatorname{GL}(\mathbb{R}^4, {J_{\tiny{\mbox{cn}}}})$-equivariant continuous function $\psi_{1,0}: C^{2}_{J_{\tiny{\mbox{cn}}}}(\mathbb{R}^{4};\mathbb{R}^{4}) \rightarrow C^{2}(\mathbb{R}^{4};\mathbb{R}^{4})$. Therefore $$\mathfrak{A}:=(\mathbb{R}^{4}, \psi_{1,0}(\mu_{11}(+))) { \xrightarrow{\text{\,\,deg\,\,}} } \mathfrak{B}:=(\mathbb{R}^{4}, \psi_{1,0}( \mu_{10}(2) ),$$
where:
$$
\mathfrak{A} =
\left\{
\begin{array}{l}
[{\it e1},{\it e2}]=\frac{1}{2}{\it e2}+ \frac{1}{2}{\it e4},
[{\it e1},{\it e3}]=-{\it e1}+{\it e3},\\
{[{\it e1},{\it e4}]}=-\frac{1}{2}{\it e2} + \frac{1}{2}{\it e4},
[{\it e2},{\it e4}]=-{\it e1}+{\it e3}
\end{array}
\right.
$$
and
$$
\mathfrak{B} =
\left\{
\begin{array}{l}
[{\it e1},{\it e2}]={\it e2}+{\it e4},
[{\it e1},{\it e3}]=-{\it e1}+{\it e3},\\
{[{\it e1},{\it e4}]}=-2\,{\it e2}+2\,{\it e4},
[{\it e2},{\it e3}]=-{\it e2}+{\it e4}
\end{array}
\right.
$$
Consider the $\operatorname{GL}(\mathbb{R}^4)$-equivariant continuous function $\kappa: C^{2}(\mathbb{R}^{4};\mathbb{R}^{4}) \rightarrow \Sigma^2(\mathbb{R}^{4}; \mathbb{R})$ defined by
$$
\kappa(\lambda)(X,Y):= \operatorname{Trace}(\lambda(X, \lambda(Y , \square)));
$$
the symmetric bilinear form $\kappa(\lambda)$ is just the Killing form of $(\mathbb{R}^4, \lambda)$. Therefore $\kappa(\mathfrak{A}) {\xrightarrow{\text{\,\,deg\,\,}} } \kappa(\mathfrak{B}) $. The matrix of $\kappa(\mathfrak{A})$ with respect to the canonical basis is:
$$
\left[ \begin {array}{cccc} 1&0&1&0\\ \noalign{\medskip}0&1&0&0
\\ \noalign{\medskip}1&0&1&0\\ \noalign{\medskip}0&0&0&1\end {array}
 \right],
$$
while the matrix of $\kappa(\mathfrak{B})$ with respect to the canonical basis is:
$$
 \left[ \begin {array}{cccc} 2&0&4&0\\ \noalign{\medskip}0&0&0&0
\\ \noalign{\medskip}4&0&2&0\\ \noalign{\medskip}0&0&0&0\end {array}
 \right]
$$
which is a contradiction because $\kappa(\mathfrak{A})$ is positive semi-definite bilinear form and $\kappa(\mathfrak{B})$ is a bilinear form with signature $(-,0,0,+)$.

\item\label{delta4unmedior4normal}
$
(\mathfrak{ d}_{4,1/2}, J_2)
\centernot
{ \xrightarrow{\text{\,\,deg\,\,}} }
(\mathfrak{ r}_{4,\alpha,\alpha },J)
$ with $\alpha \neq 0, \frac{1}{2}$. \newline

In fact, this holds because $
\mathfrak{ d}_{4,1/2}
\centernot
{ \xrightarrow{\text{\,\,deg\,\,}} }
\mathfrak{ r}_{4,\alpha,\alpha }
$.  Let us now explore this in detail. Suppose for the sake of contradiction that $
\mathfrak{ d}_{4,1/2}
{ \xrightarrow{\text{\,\,deg\,\,}} }
\mathfrak{ r}_{4,\alpha,\alpha }
$. We have two cases:

\begin{itemize}
\item If $a \neq \frac{5}{2}$.

First, consider the \(\operatorname{GL}(\mathbb{R}^4)\)-equivariant continuous function \[ \kappa_{g}: C^{2}(\mathbb{R}^{4};\mathbb{R}^{4}) \rightarrow \Sigma^2(\mathbb{R}^{4}; \mathbb{R}). \]
By assumption, we would have $ \kappa_{g}( \mu_{11}(+) ) {\xrightarrow{\text{\,\,deg\,\,}} } \kappa( \mu_{6}(a,0) ) $
where:
$$
\kappa_{g}( \mu_{11}(+) ) = (\frac{3}{2}+4g) e_{1}^{\ast} \otimes e_{1}^{\ast}
$$
and
$$
\kappa_{g}( \mu_{6}(a,0) ) = (1+2\,{a}^{2}+g+4\,ga+4\,g{a}^{2}) e_{1}^{\ast} \otimes e_{1}^{\ast}
$$

Setting $g= -\frac{3}{8}$, we obtain a contradiction, since $\kappa_{g}( \mu_{11}(+) ) $ is the zero  bilinear form while $\kappa_{g}( \mu_{6}(a,0) ) = \frac{1}{8}\, \left( 2\,a-1 \right)  \left( 2\,a-5 \right) e_{1}^{\ast} \otimes e_{1}^{\ast} $ is non-zero bilinear form.

\item If $a = \frac{5}{2}$.

We consider the \(\operatorname{GL}(\mathbb{R}^4)\)-equivariant continuous function \[ \Psi_{f,g}: C^{2}(\mathbb{R}^{4};\mathbb{R}^{4}) \rightarrow \Sigma^{3}(\mathbb{R}^{4};\mathbb{R}). \]
Therefore $$ \Psi_{0,g}( \mu_{11}(+)  ) \xrightarrow{\text{\,\,deg\,\,}} \Psi_{0,g}(  \mu_{6}(a,0)  ), $$
where:
$$
 \Psi_{0,g}( \mu_{11}(+) ) = (8g + \frac{15}{2} ) e_{1}^{\ast} \otimes e_{1}^{\ast} \otimes e_{1}^{\ast}
$$
and
\begin{eqnarray*}
\Psi_{0,g}(  \mu_{6}(a,0)  ) &=& (8\,g{a}^{3}+12\,g{a}^{2}+6\,ga+g+6+12\,{a}^{3})e_{1}^{\ast} \otimes e_{1}^{\ast} \otimes e_{1}^{\ast}\\
                             &=& (216\,g+{\frac {387}{2}}) e_{1}^{\ast} \otimes e_{1}^{\ast} \otimes e_{1}^{\ast}
\end{eqnarray*}
Setting $g= -\frac{15}{16}$, we obtain a contradiction, since $ \Psi_{0,g}( \mu_{11}(+) ) \equiv 0$
and $ \Psi_{0,g}(   \mu_{6}(a,0) ) = -9 e_{1}^{\ast} \otimes e_{1}^{\ast} \otimes e_{1}^{\ast} $.

\end{itemize}

\item
$
(\mathfrak{ d}_{4,1/2}, J_2)
%\centernot
{ \xrightarrow{\text{\,\,deg\,\,}} }
(\mathfrak{ r}_{4,1/2,1/2 },J)
$\newline

Consider the curve $g_t = \left(
          \begin{array}{cc}
            A(t) & -B(t) \\
            B(t) & A(t) \\
          \end{array}
        \right)$ where $$A(t)= \operatorname{diag}(1,e^t)
                             \qquad
                             \mbox{and}
                             \qquad
                         B(t)= 0
                         .$$
Note that $g_t \in \operatorname{GL}(\mathbb{R}^4, {J_{\tiny{\mbox{cn}}}})$ and its action on $\mu_{11}(+)$ gives
$$
\left\{
\begin{array}{l}
[{\it e1},{\it e2}]=\frac{1}{2}{\it e2},
[{\it e1},{\it e3}]={\it e3},
[{\it e1},{\it e4}]=\frac{1}{2}{\it e4},
[{\it e2},{\it e4}]={{\rm e}^{-2\,t}}{\it e3}
\end{array}
\right.
$$
Taking the limit $g_t \cdot \mu_{11}(+)$ as $t \to \infty$, we obtain $\mu_{6}(\frac{1}{2},0)$.

\item
$
(\mathfrak{ d}_{4,1/2}, J_2)
\centernot
{ \xrightarrow{\text{\,\,deg\,\,}} }
(\mathfrak{r}_{4,\alpha,1},J)
$ with $-1 < \alpha \leq 1$ and $\alpha \neq 0$\newline
See Subsection \ref{delta4unmedio} item \ref{delta4unmedior4normal}

\item\label{delta4unmedior4uno}
$
(\mathfrak{ d}_{4,1/2}, J_2)
\centernot
{ \xrightarrow{\text{\,\,deg\,\,}} }
(\mathfrak{ r}_{4,1},J)
$\newline

In fact, this holds because $
\mathfrak{ d}_{4,1/2}
\centernot
{ \xrightarrow{\text{\,\,deg\,\,}} }
\mathfrak{ r}_{4,1}
$.  Let us now explore this in detail. Suppose for the sake of contradiction that $
\mathfrak{ d}_{4,1/2}
{ \xrightarrow{\text{\,\,deg\,\,}} }
\mathfrak{ r}_{4,1}
$. First, consider the \(\operatorname{GL}(\mathbb{R}^4)\)-equivariant continuous function \[ \kappa_{g}: C^{2}(\mathbb{R}^{4};\mathbb{R}^{4}) \rightarrow \Sigma^2(\mathbb{R}^{4}; \mathbb{R}). \]
By assumption, we would have $ \kappa_{g}( \mu_{11}(+) ) {\xrightarrow{\text{\,\,deg\,\,}} } \kappa( \mu_{9} ) $
where:
$$
\kappa_{g}( \mu_{11}(+) ) = (\frac{3}{2}+4g) e_{1}^{\ast} \otimes e_{1}^{\ast}
$$
and
$$
\kappa_{g}( \mu_{9} ) = (3 + 9 g) e_{1}^{\ast} \otimes e_{1}^{\ast}
$$

Setting $g= -\frac{3}{8}$, we obtain a contradiction, since $\kappa_{g}( \mu_{11}(+) )$ is the zero bilinear form while $\kappa_{g}( \mu_{9} )= -\frac{3}{8}e_{1}^{\ast} \otimes e_{1}^{\ast} $ is non-zero negative semi-definite bilinear form.

\item\label{delta4unmedior3uno}
$
(\mathfrak{ d}_{4,1/2}, J_2)
\centernot
{ \xrightarrow{\text{\,\,deg\,\,}} }
(\mathfrak{rr}_{3,1 }\times \mathbb{R},J)
$\newline

In fact, this holds because $
\mathfrak{ d}_{4,1/2}
\centernot
{ \xrightarrow{\text{\,\,deg\,\,}} }
\mathfrak{rr}_{3,1 }\times \mathbb{R}
$.  Let us now explore this in detail. Suppose for the sake of contradiction that $
\mathfrak{ d}_{4,1/2}
{ \xrightarrow{\text{\,\,deg\,\,}} }
\mathfrak{rr}_{3,1 }\times \mathbb{R}
$. First, consider the \(\operatorname{GL}(\mathbb{R}^4)\)-equivariant continuous function \[ \kappa_{g}: C^{2}(\mathbb{R}^{4};\mathbb{R}^{4}) \rightarrow \Sigma^2(\mathbb{R}^{4}; \mathbb{R}). \]
By assumption, we would have $ \kappa_{g}( \mu_{11}(+) ) {\xrightarrow{\text{\,\,deg\,\,}} } \kappa( \mu_{3}(0,0) ) $
where:
$$
\kappa_{g}( \mu_{11}(+) ) = (\frac{3}{2}+4g) e_{1}^{\ast} \otimes e_{1}^{\ast}
$$
and
$$
\kappa_{g}( \mu_{3}(0,0) ) = (2+4g ) e_{1}^{\ast} \otimes e_{1}^{\ast}
$$

Setting $g= -\frac{3}{8}$, we obtain a contradiction, since $\kappa_{g}( \mu_{11}(+) )$ is the zero bilinear form while $\kappa_{g}( \mu_{3}(0,0) )=  \frac{1}{2} e_{1}^{\ast} \otimes e_{1}^{\ast} $ is non-zero positive semi-definite bilinear form.

\item
$
(\mathfrak{ d}_{4,1/2}, J_2)
\centernot
{ \xrightarrow{\text{\,\,deg\,\,}} }
(\mathfrak{ r}_{4,-1/2,-1/2},J)
$\newline
See Subsection \ref{delta4unmedio} item \ref{delta4unmedior4normal}

\item\label{delta4unmedior2}
$
(\mathfrak{ d}_{4,1/2}, J_2)
\centernot
{ \xrightarrow{\text{\,\,deg\,\,}} }
(\mathfrak{r}_2\times \mathbb{R}^{2},J)
$\newline

In fact, this holds because $
\mathfrak{ d}_{4,1/2}
\centernot
{ \xrightarrow{\text{\,\,deg\,\,}} }
\mathfrak{r}_2\times \mathbb{R}^{2}
$.  Let us now explore this in detail. Suppose for the sake of contradiction that $
\mathfrak{ d}_{4,1/2}
{ \xrightarrow{\text{\,\,deg\,\,}} }
\mathfrak{r}_2\times \mathbb{R}^{2}
$. First, consider the \(\operatorname{GL}(\mathbb{R}^4)\)-equivariant continuous function \[ \kappa_{g}: C^{2}(\mathbb{R}^{4};\mathbb{R}^{4}) \rightarrow \Sigma^2(\mathbb{R}^{4}; \mathbb{R}). \]
By assumption, we would have $ \kappa_{g}( \mu_{11}(+) ) {\xrightarrow{\text{\,\,deg\,\,}} } \kappa( \mu_{6}(0,0) ) $
where:
$$
\kappa_{g}( \mu_{11}(+) ) = (\frac{3}{2}+4g) e_{1}^{\ast} \otimes e_{1}^{\ast}
$$
and
$$
\kappa_{g}( \mu_{6}(0,0) ) = (1+g ) e_{1}^{\ast} \otimes e_{1}^{\ast}
$$

Setting $g= -\frac{3}{8}$, we obtain a contradiction, since $\kappa_{g}( \mu_{11}(+) )$ is the zero bilinear form while $\kappa_{g}( \mu_{6}(0,0) )=  \frac{5}{8} e_{1}^{\ast} \otimes e_{1}^{\ast} $ is non-zero positive semi-definite bilinear form.

  \item
$
(\mathfrak{ d}_{4,1/2}, J_2)
%\centernot
{ \xrightarrow{\text{\,\,deg\,\,}} }
(\mathfrak{h}_{3}\times \mathbb{R},J)
$\newline

Consider the curve $g_t = \left(
          \begin{array}{cc}
            A(t) & -B(t) \\
            B(t) & A(t) \\
          \end{array}
        \right)$ where $$A(t)= 0
                             \qquad
                             \mbox{and}
                             \qquad
                         B(t)= \left(
                                 \begin{array}{cc}
                                   0 & e^t \\
                                   -e^{2t} & 0 \\
                                 \end{array}
                               \right)
                         .$$
Note that $g_t \in \operatorname{GL}(\mathbb{R}^4, {J_{\tiny{\mbox{cn}}}})$ and its action on $\mu_{11}(+) $ gives
$$
\left\{
\begin{array}{l}
[{\it e1},{\it e3}]={\it e2},
[{\it e1},{\it e4}]=\frac{1}{2}\,{{\rm e}^{-2\,t}}{\it e1},
[{\it e2},{\it e4}]={{\rm e}^{-2\,t}}{\it e2},
[{\it e3},{\it e4}]=\frac{1}{2}\,{{\rm e}^{-2\,t}}{\it e3}
\end{array}
\right.
$$
Taking the limit $g_t \cdot \mu_{11}(+) $ as $t \to \infty$, we obtain $\mu_{1}$.

  \item
$
(\mathfrak{ d}_{4,1/2}, J_2)
\centernot
{ \xrightarrow{\text{\,\,deg\,\,}} }
(\mathfrak{ r}_{4,1,1},J)
$\newline
See Subsection \ref{delta4unmedio} item \ref{delta4unmedior4normal}
\end{enumerate}

\subsection{$(\mathfrak{ d}_{4,1/2}, J_3)$}
\begin{enumerate}

\item
$
(\mathfrak{ d}_{4,1/2}, J_3)
%\centernot
{ \xrightarrow{\text{\,\,deg\,\,}} }
(\mathfrak{ d}_{4,1/2}, J_1)
$\newline

Consider the curve $g_t = \left(
          \begin{array}{cc}
            A(t) & -B(t) \\
            B(t) & A(t) \\
          \end{array}
        \right)$ where $$A(t)=\left(
                               \begin{array}{cc}
                                 -\frac{1}{2} & 0 \\
                                 \frac{\sqrt{2}}{2} e^t & e^t \\
                               \end{array}
                             \right)
                             \qquad
                             \mbox{and}
                             \qquad
                         B(t)=0.$$
Note that $g_t \in \operatorname{GL}(\mathbb{R}^4, {J_{\tiny{\mbox{cn}}}})$ and its action on $g_t \cdot \mu_{11}(-)$ gives
$$
\left\{
\begin{array}{l}
[{\it e1},{\it e2}]={\it e2},
[{\it e1},{\it e3}]={\it e3},
[{\it e1},{\it e4}]=-1/2\,\sqrt {2}{{\rm e}^{-t}}{\it e3}+2\,{\it e4},\\
{[{\it e2},{\it e3}]}=-1/2\,\sqrt {2}{{\rm e}^{-t}}{\it e3}+{\it e4},
[{\it e2},{\it e4}]=-1/2\,{{\rm e}^{-2\,t}}{\it e3}+1/2\,\sqrt {2}{{\rm e}^{-t}}{\it e4}.
\end{array}
\right.
$$
Taking the limit $g_t \cdot \mu_{11}(-)$ as $t \to \infty$, we obtain $\mu_{10}(2)$.

\item
$
(\mathfrak{ d}_{4,1/2}, J_3)
\centernot
{ \xrightarrow{\text{\,\,deg\,\,}} }
(\mathfrak{ r}_{4,\alpha,\alpha },J)
$\newline
See Subsection \ref{delta4unmedio} item \ref{delta4unmedior4normal}

\item
$
(\mathfrak{ d}_{4,1/2}, J_3)
%\centernot
{ \xrightarrow{\text{\,\,deg\,\,}} }
(\mathfrak{ r}_{4,1/2,1/2 },J)
$\newline

Consider the curve $g_t = \left(
          \begin{array}{cc}
            A(t) & -B(t) \\
            B(t) & A(t) \\
          \end{array}
        \right)$ where $$A(t)= \operatorname{diag}(-1,e^t)
                             \qquad
                             \mbox{and}
                             \qquad
                         B(t)= 0
                         .$$
Note that $g_t \in \operatorname{GL}(\mathbb{R}^4, {J_{\tiny{\mbox{cn}}}})$ and its action on $\mu_{11}(-) $ gives
$$
\left\{
\begin{array}{l}
[{\it e1},{\it e2}]=\frac{1}{2}\,{\it e2},
[{\it e1},{\it e3}]={\it e3},
[{\it e1},{\it e4}]=\frac{1}{2}\,{\it e4},
[{\it e2},{\it e4}]=-{{\rm e}^{-2\,t}}{\it e3}
\end{array}
\right.
$$
Taking the limit $g_t \cdot \mu_{11}(-) $ as $t \to \infty$, we obtain $\mu_{6}(\frac{1}{2},0)$.
\item
$
(\mathfrak{ d}_{4,1/2}, J_3)
\centernot
{ \xrightarrow{\text{\,\,deg\,\,}} }
(\mathfrak{r}_{4,\alpha,1},J)
$ with $-1 < \alpha \leq -1$ and $\alpha \neq 0$ \newline
See Subsection \ref{delta4unmedio} item \ref{delta4unmedior4normal}

\item
$
(\mathfrak{ d}_{4,1/2}, J_3)
\centernot
{ \xrightarrow{\text{\,\,deg\,\,}} }
(\mathfrak{ r}_{4,1},J)
$\newline
See Subsection \ref{delta4unmedio} item \ref{delta4unmedior4uno}

\item
$
(\mathfrak{ d}_{4,1/2}, J_3)
\centernot
{ \xrightarrow{\text{\,\,deg\,\,}} }
(\mathfrak{rr}_{3,1 }\times \mathbb{R},J)
$\newline
See Subsection \ref{delta4unmedio} item \ref{delta4unmedior3uno}

\item
$
(\mathfrak{ d}_{4,1/2}, J_3)
\centernot
{ \xrightarrow{\text{\,\,deg\,\,}} }
(\mathfrak{ r}_{4,-1/2,-1/2},J)
$\newline
See Subsection \ref{delta4unmedio} item \ref{delta4unmedior4normal}

\item
$
(\mathfrak{ d}_{4,1/2}, J_3)
\centernot
{ \xrightarrow{\text{\,\,deg\,\,}} }
(\mathfrak{r}_2\times \mathbb{R}^{2},J)
$\newline
See Subsection \ref{delta4unmedio} item \ref{delta4unmedior2}

  \item
$
(\mathfrak{ d}_{4,1/2}, J_3)
%\centernot
{ \xrightarrow{\text{\,\,deg\,\,}} }
(\mathfrak{h}_{3}\times \mathbb{R},J)
$\newline

Consider the curve $g_t = \left(
          \begin{array}{cc}
            A(t) & -B(t) \\
            B(t) & A(t) \\
          \end{array}
        \right)$ where $$A(t)= 0
                             \qquad
                             \mbox{and}
                             \qquad
                         B(t)= \left(
                                 \begin{array}{cc}
                                   0 & e^t \\
                                   -e^{2t} & 0 \\
                                 \end{array}
                               \right)
                         .$$
Note that $g_t \in \operatorname{GL}(\mathbb{R}^4, {J_{\tiny{\mbox{cn}}}})$ and its action on $\mu_{11}(-)$ gives
$$
\left\{
\begin{array}{l}
[{\it e1},{\it e3}]={\it e2},
[{\it e1},{\it e4}]=-\frac{1}{2}{{\rm e}^{-2\,t}}{\it e1},
[{\it e2},{\it e4}]=-{{\rm e}^{-2\,t}}{\it e2},
[{\it e3},{\it e4}]=-\frac{1}{2}\,{{\rm e}^{-2\,t}}{\it e3}
\end{array}
\right.
$$
Taking the limit $g_t \cdot  \mu_{11}(-)$ as $t \to \infty$, we obtain $\mu_{1}$.

  \item
$
(\mathfrak{ d}_{4,1/2}, J_3)
\centernot
{ \xrightarrow{\text{\,\,deg\,\,}} }
(\mathfrak{ r}_{4,1,1},J)
$\newline
See Subsection \ref{delta4unmedio} item \ref{delta4unmedior4normal}
\end{enumerate}

\subsection{$(\mathfrak{r}^{'}_{2 }, J_{1}(a,b))$ with $b\neq 0$, $(a,b)\neq (0,\pm 1)$}
\begin{enumerate}

  \item
$
(\mathfrak{r}^{'}_{2 }, J_{1}(a,b))
\centernot
{ \xrightarrow{\text{\,\,deg\,\,}} }
(\mathfrak{r}^{'}_{2 }, J_{1}(0,-1))
$\newline

Suppose for the sake of contradiction that $(\mathfrak{r}^{'}_{2 }, J_{1}(a,b))
\centernot
{ \xrightarrow{\text{\,\,deg\,\,}} }
(\mathfrak{r}^{'}_{2 }, J_{1}(0,-1))$. We have two cases:

\begin{itemize}

\item If $a\neq 0$.

First, we consider the $\operatorname{GL}(\mathbb{R}^4, {J_{\tiny{\mbox{cn}}}})$-equivariant continuous function $\psi_{\alpha,\beta}: C^{2}_{J_{\tiny{\mbox{cn}}}}(\mathbb{R}^{4};\mathbb{R}^{4}) \rightarrow C^{2}(\mathbb{R}^{4};\mathbb{R}^{4})$ with $ \alpha={\frac {{a}^{2}+1-b}{ab}}$ and $\beta=\frac{1}{a} $. Therefore $$\mathfrak{A}:=(\mathbb{R}^{4}, \psi_{\alpha,\beta}( \mu_{3}(a,b) ) { \xrightarrow{\text{\,\,deg\,\,}} } \mathfrak{B}:=(\mathbb{R}^{4}, \psi_{\alpha,\beta}( \mu_{3}(0,-1) ),$$
where:
$$
\mathfrak{A} =
\left\{
\begin{array}{l}
[{\it e1},{\it e2}]={\frac { \left( {a}^{2}+1-{b}^{2} \right) {\it e4 }}{ab}},
[{\it e1},{\it e4}]=-{\frac { \left( {a}^{2}+1-{b}^{2} \right) {\it e2}}{ab}},\\
{[{\it e2},{\it e3}]}={\frac { \left( {b}^{2}+{a }^{2}+1 \right) {\it e4}}{b}},
[{\it e3},{\it e4}]={\frac { \left( {b}^{2}+{a}^{2}+1 \right) {\it e2}}{b}}
\end{array}
\right.
$$
and
$$
\mathfrak{B} =
\left\{
\begin{array}{l}
[{\it e1},{\it e2}]={\it e2}+{\frac { \left( {a}^{2}+1+b \right) {\it e4}}{ab}},
[{\it e1},{\it e4}]=-{\frac { \left( {a}^{2}+1+b \right) {\it e2}}{ab}}+{\it e4},\\
{[{\it e2},{\it e3}]}={\frac { \left( {a}^{2}+1+b \right) {\it e2}}{ab}}-{\it e4},
[{\it e3},{\it e4}]=-{\it e2}-{\frac { \left( {a}^{2}+1+b \right) {\it e4}}{ab}}
\end{array}
\right.
$$
 We obtain a contradiction, since $\mathfrak{A}$ is a unimodular algebra but $\mathfrak{B}$ is not unimodular.

\item If $a=0$ and $b \neq \pm 1$.

We consider the $\operatorname{GL}(\mathbb{R}^4, {J_{\tiny{\mbox{cn}}}})$-equivariant continuous function $\theta_{\beta}: C^{2}_{J_{\tiny{\mbox{cn}}}}(\mathbb{R}^{4};\mathbb{R}^{4}) \rightarrow C^{2}(\mathbb{R}^{4};\mathbb{R}^{4})$ with $\beta = \frac{b}{1-b} $, which is well defined because $b\neq 1$. Therefore $$\mathfrak{A}:=(\mathbb{R}^{4}, \theta_{\beta}( \mu_{3}(0,b) ) { \xrightarrow{\text{\,\,deg\,\,}} } \mathfrak{B}:=(\mathbb{R}^{4}, \theta_{\beta}( \mu_{3}(0,-1) ),$$
where:
$$
\mathfrak{A} =
\left\{
\begin{array}{l}
[{\it e1},{\it e2}]= \left( 1+b \right) {\it e4},
[{\it e1},{\it e4}]= \left( -b-1 \right) {\it e2}
\end{array}
\right.
$$
and
$$
\mathfrak{B} =
\left\{
\begin{array}{l}
[{\it e1},{\it e2}]= {\frac { \left( 1+b \right) {\it e4}}{1-b}},
[{\it e1},{\it e4}]={\frac { \left( 1+b \right) {\it e2}}{b-1}},\\
{[{\it e2},{\it e3}]}={\frac { \left( 1+b \right) {\it e2}}{1-b}},
[{\it e3},{\it e4}]={\frac { \left( 1+b \right) {\it e4}}{b-1}}.
\end{array}
\right.
$$
 We obtain a contradiction, since $\mathfrak{A}$ is a unimodular algebra but $\mathfrak{B}$ is not unimodular (since $b\neq -1$)

\end{itemize}

  \item
$
(\mathfrak{r}^{'}_{2 }, J_{1}(a,b))
%\centernot
{ \xrightarrow{\text{\,\,deg\,\,}} }
(\mathfrak{h}_{3}\times \mathbb{R},J)
$ with $a\neq 0$ or $b \neq -1$ \newline

Consider the curve $g_t = \left(
          \begin{array}{cc}
            A(t) & -B(t) \\
            B(t) & A(t) \\
          \end{array}
        \right)$ where $$A(t)= \left(
                                 \begin{array}{cc}
                                   e^t & 0 \\
                                   -\frac{a e^{2t}}{a^2+(b+1)^2 } & 1 \\
                                 \end{array}
                               \right)
                             \qquad
                             \mbox{and}
                             \qquad
                         B(t)= \left(
                                 \begin{array}{cc}
                                   0 & 0 \\
                                   \frac{(b+1)e^{2t}}{a^2+(b+1)^2} & 0 \\
                                 \end{array}
                               \right)
                         .$$
Note that $g_t \in \operatorname{GL}(\mathbb{R}^4, {J_{\tiny{\mbox{cn}}}})$ and its action on $\mu_{3}(a,b) $ gives
$$
\left\{
\begin{array}{l}
[{\it e1},{\it e2}]={{\rm e}^{-t}}{\it e2},
[{\it e1},{\it e3}]={\it e2},
{[{\it e1},{\it e4}]}={{\rm e}^{-t}}{\it e4},\\
{[{\it e2},{\it e3}]}={{\rm e}^{-t}}a{\it e2}+{{\rm e}^{-t}}b{\it e4},
[{\it e3},{\it e4}]={{\rm e}^{-t}}b{\it e2}-{{\rm e}^{-t}}a{\it e4}
\end{array}
\right.
$$
Taking the limit $g_t \cdot \mu_{3}(a,b)$ as $t \to \infty$, we obtain $\mu_{1}$.

 \item
$
(\mathfrak{r}^{'}_{2 }, J_{1}(a,b))
\centernot
{ \xrightarrow{\text{\,\,deg\,\,}} }
(\mathfrak{ r}_{4,1,1},J)
$\newline

In fact, this holds because $
\mathfrak{r}^{'}_{2 }
\centernot
{ \xrightarrow{\text{\,\,deg\,\,}} }
\mathfrak{ r}_{4,1,1}
$ since $\operatorname{Dim}([\mathfrak{r}^{'}_{2 },\mathfrak{r}^{'}_{2 }])=2$ and $\operatorname{Dim}([\mathfrak{ r}_{4,1,1} ,\mathfrak{ r}_{4,1,1} ] )=3$.

\end{enumerate}

\subsection{$(\mathfrak{r}^{'}_{2 }, J_{3})$}
\begin{enumerate}

  \item
$
(\mathfrak{r}^{'}_{2 }, J_{3})
{ \xrightarrow{\text{\,\,deg\,\,}} }
(\mathfrak{r}^{'}_{2 }, J_{1}(0,-1))
$\newline

Consider the curve $g_t = \left(
          \begin{array}{cc}
            A(t) & -B(t) \\
            B(t) & A(t) \\
          \end{array}
        \right)$ where $$A(t)= \operatorname{diag}(1,e^{-t})
                             \qquad
                             \mbox{and}
                             \qquad
                         B(t)=0
                         .$$
Note that $g_t \in \operatorname{GL}(\mathbb{R}^4, {J_{\tiny{\mbox{cn}}}})$ and its action on $ \mu_{5} $ gives
$$
\left\{
\begin{array}{l}
[{\it e1},{\it e2}]={\it e2},
[{\it e1},{\it e3}]={{\rm e}^{-t}}{\it e4},
[{\it e1},{\it e4}]={\it e4},
[{\it e2},{\it e3}]=-{\it e4},
[{\it e3},{\it e4}]=-{\it e2}
\end{array}
\right.
$$
Taking the limit $g_t \cdot \mu_{5} $ as $t \to \infty$, we obtain $ \mu_{3}(0,-1) $.

  \item
$
(\mathfrak{r}^{'}_{2 }, J_{3})
%\centernot
{ \xrightarrow{\text{\,\,deg\,\,}} }
(\mathfrak{h}_{3}\times \mathbb{R},J)
$\newline

Consider the curve $g_t = \left(
          \begin{array}{cc}
            A(t) & -B(t) \\
            B(t) & A(t) \\
          \end{array}
        \right)$ where $$A(t)= \operatorname{diag}(e^{t},0)
                             \qquad
                             \mbox{and}
                             \qquad
                         B(t)=\operatorname{diag}(0, -e^{2t})
                         .$$
Note that $g_t \in \operatorname{GL}(\mathbb{R}^4, {J_{\tiny{\mbox{cn}}}})$ and its action on $\mu_{5}$ gives
$$
\left\{
\begin{array}{l}
[{\it e1},{\it e2}]={{\rm e}^{-t}}{\it e2},
[{\it e1},{\it e3}]={\it e2},
[{\it e1},{\it e4}]={{\rm e}^{-t}}{\it e4},\\
{[{\it e2},{\it e3}]}=-{{\rm e}^{-t}}{\it e4},
[{\it e3},{\it e4}]=-{{\rm e}^{-t}}{\it e2}
\end{array}
\right.
$$
Taking the limit $g_t \cdot \mu_{5}$ as $t \to \infty$, we obtain $ \mu_{1}$.

 \item
$
(\mathfrak{r}^{'}_{2 }, J_{3})
\centernot
{ \xrightarrow{\text{\,\,deg\,\,}} }
(\mathfrak{ r}_{4,1,1},J)
$\newline

In fact, this holds because $
\mathfrak{r}^{'}_{2 }
\centernot
{ \xrightarrow{\text{\,\,deg\,\,}} }
\mathfrak{ r}_{4,1,1}
$ since $\operatorname{Dim}([\mathfrak{r}^{'}_{2 },\mathfrak{r}^{'}_{2 }])=2$ and $\operatorname{Dim}([\mathfrak{ r}_{4,1,1} ,\mathfrak{ r}_{4,1,1} ] )=3$.

\end{enumerate}

\subsection{$(\mathfrak{r}^{'}_{4,c,d},J_{k})$ with $d>0$ and $k=1,2$}
\begin{enumerate}

  \item
$
(\mathfrak{r}^{'}_{4,c,d},J_{k})
%\centernot
{ \xrightarrow{\text{\,\,deg\,\,}} }
(\mathfrak{h}_{3}\times \mathbb{R},J)
$ \newline

Consider the curve $g_t = \left(
          \begin{array}{cc}
            A(t) & -B(t) \\
            B(t) & A(t) \\
          \end{array}
        \right)$ where $$A(t)=\left(
                                \begin{array}{cc}
                                  e^t & 0 \\
                                  \frac{(-1)^{k+1}e^{2t}d}{d^2+(c-1)^2} & 1 \\
                                \end{array}
                              \right)
                             \qquad
                             \mbox{and}
                             \qquad
                         B(t)= \left(
                                 \begin{array}{cc}
                                   0 & 0 \\
                                   \frac{e^{2t}(c-1)}{d^2+(c-1)^2} & 0 \\
                                 \end{array}
                               \right)
                         .$$
Note that $g_t \in \operatorname{GL}(\mathbb{R}^4, {J_{\tiny{\mbox{cn}}}})$ and its action on $\mu_{6}(c,d)$ gives
$$
\left\{
\begin{array}{l}
[{\it e1},{\it e2}]={{\rm e}^{-t}}c{\it e2}-{{\rm e}^{-t}} \left( -1 \right) ^{k}d{\it e4},
[{\it e1},{\it e3}]={\it e2}+{{\rm e}^{-t}}{\it e3},\\
{[{\it e1},{\it e4}]}={{\rm e}^{-t}} \left( -1 \right) ^{k}d{\it e2}+{{\rm e}^{-t}}c{\it e4}
\end{array}
\right.
$$
Taking the limit $g_t \cdot \mu_{6}(c,d)$ as $t \to \infty$, we obtain $\mu_{1}$.

 \item
$
(\mathfrak{r}^{'}_{4,c,d},J_{k})
\centernot
{ \xrightarrow{\text{\,\,deg\,\,}} }
(\mathfrak{ r}_{4,1,1},J)
$ with $k=1,2$, $d>0$ and $c\neq -\frac{1}{2}$.\newline

In fact, this holds because $
\mathfrak{r}^{'}_{4,c,d}
\centernot
{ \xrightarrow{\text{\,\,deg\,\,}} }
\mathfrak{ r}_{4,1,1}
$ .  Let us now explore this in detail. Suppose for the sake of contradiction that $
\mathfrak{r}^{'}_{4,c,d}
{ \xrightarrow{\text{\,\,deg\,\,}} }
\mathfrak{ r}_{4,1,1}
$. We have two cases:

\begin{itemize}
  \item If ${1-2\,c+{c}^{2}-3\,{d}^{2}} \neq 0$

First, consider the \(\operatorname{GL}(\mathbb{R}^4)\)-equivariant continuous function \[ \kappa_{g}: C^{2}(\mathbb{R}^{4};\mathbb{R}^{4}) \rightarrow \Sigma^2(\mathbb{R}^{4}; \mathbb{R}). \]
By assumption, we would have $ \kappa_{g}(  \mu_{6}(c,d) ) {\xrightarrow{\text{\,\,deg\,\,}} } \kappa( \mu_{6}(1,0) ) $
where:
$$
\kappa_{g}( \mu_{6}(c,d) ) = (1-2\,{d}^{2}+2\,{c}^{2}+g+4\,gc+4\,g{c}^{2})e_1^{\ast} \otimes e_1^{\ast}
$$
and
$$
\kappa_{g}( \mu_{6}(1,0) ) = (3+9g)e_1^{\ast} \otimes e_1^{\ast}
$$

Setting $g= -{\frac {1-2\,{d}^{2}+2\,{c}^{2}}{ \left( 1+2\,c \right) ^{2}}}$, we obtain a contradiction, since $\kappa_{g}( \mu_{6}(c,d) ) $ is zero bilinear form while $\kappa_{g}(  \mu_{6}(1,0)  ) = -6\,{\frac {1-2\,c+{c}^{2}-3\,{d}^{2}}{ \left( 1+2\,c \right) ^{2}}} e_1^{\ast} \otimes e_1^{\ast}$ is non-zero bilinear form.

\item If ${1-2\,c+{c}^{2}-3\,{d}^{2}} = 0$, and so
\begin{eqnarray}\label{eq7}
  9\,c-18\,{c}^{2}+9\,{c}^{3}-27\,c{d}^{2} & = &0
\end{eqnarray}

We consider the \(\operatorname{GL}(\mathbb{R}^4)\)-equivariant continuous function \[ \Psi_{f,g}: C^{2}(\mathbb{R}^{4};\mathbb{R}^{4}) \rightarrow \Sigma^{3}(\mathbb{R}^{4};\mathbb{R}). \]
Therefore $$ \Psi_{0,g}(  \mu_{6}(c,d)  ) \xrightarrow{\text{\,\,deg\,\,}} \Psi_{0,g}(  \mu_{6}(1,0)  ), $$
where:
$$
 \Psi_{0,g}( \mu_{6}(c,d) ) = (8\,g{c}^{3}+12\,g{c}^{2}+6\,gc+g+6+12\,{c}^{3}-36\,c{d}^{2}) e_1^{\ast} \otimes e_1^{\ast} \otimes e_1^{\ast}
$$
and
$$
\Psi_{0,g}( \mu_{6}(1,0)  )= (27 g+18)  e_1^{\ast} \otimes e_1^{\ast} \otimes e_1^{\ast}
$$
Setting $g=-6\,{\frac {1+2\,{c}^{3}-6\,c{d}^{2}}{ \left( 2\,c+1 \right) ^{3}}}$ we obtain  $ \Psi_{0,g}( \mu_{6}(c,d) ) \equiv 0$
and $ \Psi_{0,g}(  \mu_{6}(1,0) ) = P(c,d) e_1^{\ast} \otimes e_1^{\ast} \otimes e_1^{\ast}$ with
\begin{eqnarray*}
P(c,d) & = &-36\,{\frac {4+5\,{c}^{3}-27\,c{d}^{2}-6\,{c}^{2}-3\,c}{ \left( 2\,c+1  \right) ^{3}}} \\
       & = & -36\,{\frac {4+5\,{c}^{3}-( 9\,c-18\,{c}^{2}+9\,{c}^{3} )-6\,{c}^{2}-3\,c}{ \left( 2\,c+1 \right) ^{3}}}\\
       & = &144\,{\frac { \left( c-1 \right) ^{3}}{ \left( 2\,c+1 \right) ^{3}}}
\end{eqnarray*}

as shown in Equation (\ref{eq7}). Note that $P( c,d ) \neq 0$, since otherwise we would have $ c=1 $, and thus $d=0$ from Equation (\ref{eq7}) again.
This is a contradiction.

\end{itemize}

 \item
$
(\mathfrak{r}^{'}_{4,-\frac{1}{2},d},J_{k})
\centernot
{ \xrightarrow{\text{\,\,deg\,\,}} }
(\mathfrak{ r}_{4,1,1},J)
$ with $k=1,2$, $d>0$.\newline
In fact, this holds because $
\mathfrak{r}^{'}_{4,-\frac{1}{2},d}
\centernot
{ \xrightarrow{\text{\,\,deg\,\,}} }
\mathfrak{ r}_{4,1,1}
$: the Lie algebra $\mathfrak{r}^{'}_{4,-\frac{1}{2},d}$ is a unimodular, while $\mathfrak{ r}_{4,1,1}$ is not unimodular.

\end{enumerate}

\subsection{$(\mathfrak{rr}^{'}_{3,a}\times\mathbb{R},J_{k})$ with $a>0$ and $k=1,2$}
\begin{enumerate}

  \item
$
(\mathfrak{rr}^{'}_{3,a}\times\mathbb{R},J_{k})
%\centernot
{ \xrightarrow{\text{\,\,deg\,\,}} }
(\mathfrak{h}_{3}\times \mathbb{R},J)
$\newline

Consider the curve $g_t = \left(
          \begin{array}{cc}
            A(t) & -B(t) \\
            B(t) & A(t) \\
          \end{array}
        \right)$ where $$A(t)=\left(
                                \begin{array}{cc}
                                  e^t & 0 \\
                                  \frac{e^{2t}}{a^2+1} & 1 \\
                                \end{array}
                              \right)
                             \qquad
                             \mbox{and}
                             \qquad
                         B(t)= \left(
                                 \begin{array}{cc}
                                   0 & 0 \\
                                   \frac{(-1)^{k+1} a e^{2t}}{a^2+1} & 0 \\
                                 \end{array}
                               \right)
                         .$$
Note that $g_t \in \operatorname{GL}(\mathbb{R}^4, {J_{\tiny{\mbox{cn}}}})$ and its action on $\mu_{7}(a)$ gives
$$
\left\{
\begin{array}{l}
[{\it e1},{\it e2}]={{\rm e}^{-t}} \left( -1 \right) ^{k+1}a{\it e2}+{{\rm e}^{-t}}{\it e4},
[{\it e1},{\it e3}]={\it e2},\\
{[{\it e1},{\it e4}]}=-{{\rm e}^{-t}}{\it e2} + {{\rm e}^{-t}} \left( -1 \right) ^{k+1}a{\it e4}
\end{array}
\right.
$$
Taking the limit $g_t \cdot \mu_{7}(a)$ as $t \to \infty$, we obtain $\mu_{1} $.

 \item
$
(\mathfrak{rr}^{'}_{3,a}\times\mathbb{R},J_{k})
\centernot
{ \xrightarrow{\text{\,\,deg\,\,}} }
(\mathfrak{ r}_{4,1,1},J)
$ with $k=1,2$ and $a>0$ \newline

In fact, this holds because $
\mathfrak{rr}^{'}_{3,a}\times\mathbb{R}
\centernot
{ \xrightarrow{\text{\,\,deg\,\,}} }
\mathfrak{ r}_{4,1,1}
$,  since $\operatorname{Dim}([\mathfrak{rr}^{'}_{3,a}\times\mathbb{R},\mathfrak{rr}^{'}_{3,a}\times\mathbb{R}])=2$ and $\operatorname{Dim}([\mathfrak{ r}_{4,1,1} ,\mathfrak{ r}_{4,1,1} ] )=3$.

\end{enumerate}

\subsection{$(\mathfrak{ d}_{4,1/2}, J_1)$}
\begin{enumerate}

  \item
$
(\mathfrak{ d}_{4,1/2}, J_1)
%\centernot
{ \xrightarrow{\text{\,\,deg\,\,}} }
(\mathfrak{h}_{3}\times \mathbb{R},J)
$\newline

Consider the curve $g_t = \left(
          \begin{array}{cc}
            A(t) & -B(t) \\
            B(t) & A(t) \\
          \end{array}
        \right)$ where $$A(t)= \operatorname{diag}(e^{t},0)
                             \qquad
                             \mbox{and}
                             \qquad
                         B(t)= \left(
                                 \begin{array}{cc}
                                   0 & 0 \\
                                   \frac{1}{2}e^{2t} & -1 \\
                                 \end{array}
                               \right)
                         .$$
Note that $g_t \in \operatorname{GL}(\mathbb{R}^4, {J_{\tiny{\mbox{cn}}}})$ and its action on $\mu_{10}(2) $ gives
$$
\left\{
\begin{array}{l}
[{\it e1},{\it e2}]=2\,{{\rm e}^{-t}}{\it e2},
[{\it e1},{\it e3}]={\it e2}+{{\rm e}^{-t}}{\it e3},
[{\it e1},{\it e4}]={{\rm e}^{-t}}{\it e4},
[{\it e3},{\it e4}]={{\rm e}^{-t}}{\it e2}
\end{array}
\right.
$$
Taking the limit $g_t \cdot \mu_{10}(2)$ as $t \to \infty$, we obtain $\mu_{1} $.

 \item
$
(\mathfrak{ d}_{4,1/2}, J_1)
\centernot
{ \xrightarrow{\text{\,\,deg\,\,}} }
(\mathfrak{ r}_{4,1,1},J)
$\newline

See Subsection \ref{delta4unmedio} item \ref{delta4unmedior4normal}.

\end{enumerate}

\subsection{$(\mathfrak{ r}_{4,\alpha,\alpha },J)$ with $ \alpha \neq 0,1$}\label{r4aa}
\begin{enumerate}

  \item\label{r4aah3}
$
(\mathfrak{ r}_{4,\alpha,\alpha },J)
%\centernot
{ \xrightarrow{\text{\,\,deg\,\,}} }
(\mathfrak{h}_{3}\times \mathbb{R},J)
$\newline

Consider the curve $g_t = \left(
          \begin{array}{cc}
            A(t) & -B(t) \\
            B(t) & A(t) \\
          \end{array}
        \right)$ where $$A(t)= \operatorname{diag}(e^t,1)
                             \qquad
                             \mbox{and}
                             \qquad
                         B(t)= \left(
                                 \begin{array}{cc}
                                   0 & 0 \\
                                   \frac{e^{2t}}{\alpha-1} & 0 \\
                                 \end{array}
                               \right)
                         .$$
Note that $g_t \in \operatorname{GL}(\mathbb{R}^4, {J_{\tiny{\mbox{cn}}}})$ and its action on $\mu_{6}(\alpha,0)$ gives
$$
\left\{
\begin{array}{l}
[{\it e1},{\it e2}]={{\rm e}^{-t}}\alpha\,{\it e2},
[{\it e1},{\it e3}]={\it e2}+{{\rm e}^{-t}}{\it e3},
[{\it e1},{\it e4}]={{\rm e}^{-t}}\alpha\,{\it e4}
\end{array}
\right.
$$
Taking the limit $g_t \cdot \mu_{6}(\alpha,0)$ as $t \to \infty$, we obtain $\mu_{1}$.

 \item\label{r4ar411}.
$
(\mathfrak{ r}_{4, a, a },J)
\centernot
{ \xrightarrow{\text{\,\,deg\,\,}} }
(\mathfrak{ r}_{4,1,1},J)
$ with $a\neq 0,1, -\frac{1}{2}$, \newline

In fact, this holds because $
\mathfrak{ r}_{4, a, a }
\centernot
{ \xrightarrow{\text{\,\,deg\,\,}} }
\mathfrak{ r}_{4,1,1}
$.  Let us now explore this in detail. Suppose for the sake of contradiction that $
\mathfrak{ r}_{4, a, a }
{ \xrightarrow{\text{\,\,deg\,\,}} }
\mathfrak{ r}_{4,1,1}
$. First, consider the \(\operatorname{GL}(\mathbb{R}^4)\)-equivariant continuous function \[ \kappa_{g}: C^{2}(\mathbb{R}^{4};\mathbb{R}^{4}) \rightarrow \Sigma^2(\mathbb{R}^{4}; \mathbb{R}). \]
By assumption, we would have $ \kappa_{g}( \mu_{6}(a,0) ) {\xrightarrow{\text{\,\,deg\,\,}} } \kappa( \mu_{6}(1,0) ) $
where:
$$
\kappa_{g}( \mu_{6}(a,0) ) = (1+2\,{a}^{2}+g+4\,ga+4\,g{a}^{2}) e_{1}^{\ast} \otimes e_{1}^{\ast}
$$
and
$$
\kappa_{g}( \mu_{6}(1,0) ) = (3+9g ) e_{1}^{\ast} \otimes e_{1}^{\ast}
$$

Setting $g= -{\frac {1+2\,{a}^{2}}{ \left( 1+2\,a \right) ^{2}}}$, we obtain a contradiction, since $\kappa_{g}( \mu_{6}(a,0) ) $ is zero bilinear form while $\kappa_{g}( \mu_{6}(1,0) ) = -6\,{\frac { \left( a-1 \right) ^{2}}{ \left( 1+2\,a \right) ^{2}}} e_{1}^{\ast} \otimes e_{1}^{\ast} $ is non-zero negative semi-definite bilinear form.

\item
$
(\mathfrak{ r}_{4, -\frac{1}{2}, -\frac{1}{2} },J)
\centernot
{ \xrightarrow{\text{\,\,deg\,\,}} }
(\mathfrak{ r}_{4,1,1},J)
$\newline
In fact, this holds because $
\mathfrak{ r}_{4, -\frac{1}{2}, -\frac{1}{2} }
\centernot
{ \xrightarrow{\text{\,\,deg\,\,}} }
\mathfrak{ r}_{4,1,1}
$: the Lie algebra $\mathfrak{ r}_{4, -\frac{1}{2}, -\frac{1}{2} }$ is a unimodular, while $\mathfrak{ r}_{4,1,1}$ is not unimodular.

\end{enumerate}

\subsection{$(\mathfrak{r}_{4,\alpha,1},J)$ with $ -1< \alpha \leq 1$, $\alpha \neq 0$. }
\begin{enumerate}

  \item
$
(\mathfrak{r}_{4,\alpha,1},J)
%\centernot
{ \xrightarrow{\text{\,\,deg\,\,}} }
(\mathfrak{h}_{3}\times \mathbb{R},J)
$\newline
See subsection \ref{r4aa} item \ref{r4aah3}.

 \item
$
(\mathfrak{r}_{4,\alpha,1},J)
\centernot
{ \xrightarrow{\text{\,\,deg\,\,}} }
(\mathfrak{ r}_{4,1,1},J)
$\newline
See subsection \ref{r4aa} item \ref{r4ar411}, and using $\mathfrak{r}_{4,\alpha,1} \equiv \mathfrak{r}_{4,\frac{1}{\alpha},\frac{1}{\alpha}}$.

\end{enumerate}

\subsection{$(\mathfrak{ r}_{4,1},J)$}
\begin{enumerate}

  \item
$
(\mathfrak{ r}_{4,1},J)
%\centernot
{ \xrightarrow{\text{\,\,deg\,\,}} }
(\mathfrak{h}_{3}\times \mathbb{R},J)
$\newline

Consider the curve $g_t = \left(
          \begin{array}{cc}
            A(t) & -B(t) \\
            B(t) & A(t) \\
          \end{array}
        \right)$ where $$A(t)= \operatorname{diag}(e^t,e^{2t})
                             \qquad
                             \mbox{and}
                             \qquad
                         B(t)= 0
                         .$$
Note that $g_t \in \operatorname{GL}(\mathbb{R}^4, {J_{\tiny{\mbox{cn}}}})$ and its action on $ \mu_{9}$ gives
$$
\left\{
\begin{array}{l}
[{\it e1},{\it e2}]={{\rm e}^{-t}}{\it e2},
[{\it e1},{\it e3}]={\it e2}+{{\rm e}^{-t}}{\it e3},
[{\it e1},{\it e4}]={{\rm e}^{-t}}{\it e4}.
\end{array}
\right.
$$
Taking the limit $g_t \cdot \mu_{9}$ as $t \to \infty$, we obtain $\mu_{1}$.

 \item
$
(\mathfrak{ r}_{4,1},J)
%\centernot
{ \xrightarrow{\text{\,\,deg\,\,}} }
(\mathfrak{ r}_{4,1,1},J)
$\newline

Consider the curve $g_t = \left(
          \begin{array}{cc}
            A(t) & -B(t) \\
            B(t) & A(t) \\
          \end{array}
        \right)$ where $$A(t)= \operatorname{diag}(1,e^{-t})
                             \qquad
                             \mbox{and}
                             \qquad
                         B(t)=0
                         .$$
Note that $g_t \in \operatorname{GL}(\mathbb{R}^4, {J_{\tiny{\mbox{cn}}}})$ and its action on $\mu_{9}$ gives
$$
\left\{
\begin{array}{l}
[{\it e1},{\it e2}]={\it e2},
[{\it e1},{\it e3}]={\it e3}+{{\rm e}^{-t}}{\it e4},
[{\it e1},{\it e4}]={\it e4}
\end{array}
\right.
$$
Taking the limit $g_t \cdot \mu_{9}$ as $t \to \infty$, we obtain $\mu_{6}(1,0)$.

\end{enumerate}

\subsection{$(\mathfrak{rr}_{3,1 }\times \mathbb{R},J)$}
\begin{enumerate}

  \item
$
(\mathfrak{rr}_{3,1 }\times \mathbb{R},J)
%\centernot
{ \xrightarrow{\text{\,\,deg\,\,}} }
(\mathfrak{h}_{3}\times \mathbb{R},J)
$\newline

Consider the curve $g_t = \left(
          \begin{array}{cc}
            A(t) & -B(t) \\
            B(t) & A(t) \\
          \end{array}
        \right)$ where $$A(t)= \operatorname{diag}(e^t,1)
                             \qquad
                             \mbox{and}
                             \qquad
                         B(t)= \left(
                                 \begin{array}{cc}
                                   0 & 0 \\
                                   e^{2t} & 0 \\
                                 \end{array}
                               \right)
                         .$$
Note that $g_t \in \operatorname{GL}(\mathbb{R}^4, {J_{\tiny{\mbox{cn}}}})$ and its action on $ \mu_{3}(0,0) $ gives
$$
\left\{
\begin{array}{l}
[{\it e1},{\it e2}]={{\rm e}^{-t}}{\it e2},
[{\it e1},{\it e3}]={\it e2},
[{\it e1},{\it e4}]={{\rm e}^{-t}}{\it e4}.
\end{array}
\right.
$$
Taking the limit $g_t \cdot \mu_{3}(0,0)$ as $t \to \infty$, we obtain $\mu_{1}$.

 \item
$
(\mathfrak{rr}_{3,1 }\times \mathbb{R},J)
\centernot
{ \xrightarrow{\text{\,\,deg\,\,}} }
(\mathfrak{ r}_{4,1,1},J)
$\newline

In fact, this holds because $
\mathfrak{rr}_{3,1 }\times \mathbb{R}
\centernot
{ \xrightarrow{\text{\,\,deg\,\,}} }
\mathfrak{ r}_{4,1,1}
$ since $\operatorname{Dim}([\mathfrak{rr}_{3,1 }\times \mathbb{R},\mathfrak{rr}_{3,1 }\times \mathbb{R}])=2$ and $\operatorname{Dim}([\mathfrak{ r}_{4,1,1} ,\mathfrak{ r}_{4,1,1} ] )=3$.

\end{enumerate}

\subsection{$(\mathfrak{r}^{'}_{2 }, J_{1}(0,-1)) $}
\begin{enumerate}

  \item
$
(\mathfrak{r}^{'}_{2 }, J_{1}(0,-1))
%\centernot
{ \xrightarrow{\text{\,\,deg\,\,}} }
(\mathfrak{a}_{4},J)
$\newline

Consider the curve $g_t = \left(
          \begin{array}{cc}
            A(t) & -B(t) \\
            B(t) & A(t) \\
          \end{array}
        \right)$ where $$A(t)= \operatorname{diag}(e^t,1)
                             \qquad
                             \mbox{and}
                             \qquad
                         B(t)= 0
                         .$$
Note that $g_t \in \operatorname{GL}(\mathbb{R}^4, {J_{\tiny{\mbox{cn}}}})$ and its action on $\mu_{3}(0,-1)$ gives
$$
\left\{
\begin{array}{l}
[{\it e1},{\it e2}]={{\rm e}^{-t}}{\it e2},
[{\it e1},{\it e4}]={{\rm e}^{-t}}{\it e4},
[{\it e2},{\it e3}]=-{{\rm e}^{-t}}{\it e4},
[{\it e3},{\it e4}]=-{{\rm e}^{-t}}{\it e2}
\end{array}
\right.
$$
Taking the limit $g_t \cdot \mu_{3}(0,-1)$ as $t \to \infty$, we obtain $\mu_{0}$.

\end{enumerate}

\subsection{$(\mathfrak{h}_{3}\times \mathbb{R},J) $}
\begin{enumerate}

\item
$
{(\mathfrak{rh}_{3}\times \mathbb{R},J)}
\xrightarrow{\text{\,\,deg\,\,}}
{(\mathfrak{a}_{4},J)}
$\newline

Consider the curve $g_t = \left(
          \begin{array}{cc}
            A(t) & -B(t) \\
            B(t) & A(t) \\
          \end{array}
        \right)$ where $$A(t)= \operatorname{diag}(1,e^{-t})
                             \qquad
                             \mbox{and}
                             \qquad
                         B(t)=0
                         .$$
Note that $g_t \in \operatorname{GL}(\mathbb{R}^4, {J_{\tiny{\mbox{cn}}}})$ and its action on $\mu_{1}$ gives
$$
\left\{
\begin{array}{l}
[{\it e1},{\it e3}]={\frac {1}{{{\rm e}^{t}}}}{\it e2}
\end{array}
\right.
$$
Taking the limit $g_t \cdot \mu_{1}$ as $t \to \infty$, we obtain $\mu_{0}$.

\end{enumerate}

\subsection{$ (\mathfrak{ r}_{4,1,1},J)$}
\begin{enumerate}

  \item
$
(\mathfrak{ r}_{4,1,1},J)
%\centernot
{ \xrightarrow{\text{\,\,deg\,\,}} }
(\mathfrak{a}_{4},J)
$\newline

Consider the curve $g_t = \left(
          \begin{array}{cc}
            A(t) & -B(t) \\
            B(t) & A(t) \\
          \end{array}
        \right)$ where $$A(t)= \operatorname{diag}(e^t,1)
                             \qquad
                             \mbox{and}
                             \qquad
                         B(t)= 0
                         .$$
Note that $g_t \in \operatorname{GL}(\mathbb{R}^4, {J_{\tiny{\mbox{cn}}}})$ and its action on $\mu_{6}(1,0)$ gives
$$
\left\{
\begin{array}{l}
[{\it e1},{\it e2}]={{\rm e}^{-t}}{\it e2},
[{\it e1},{\it e3}]={{\rm e}^{-t}}{\it e3},
[{\it e1},{\it e4}]={{\rm e}^{-t}}{\it e4}
\end{array}
\right.
$$
Taking the limit $g_t \cdot \mu_{6}(1,0)$ as $t \to \infty$, we obtain $\mu_{0}$.

\end{enumerate}

\newpage

\begin{center}
\begin{tikzpicture}

\node[circle,fill,inner sep=0pt,minimum size=3pt,label=left:{$(\mathfrak{r}_{2}',J_3)$}] (a1) at (0,+0) {};
\node[circle,fill,inner sep=0pt,minimum size=3pt,label=above:{\begin{tabular}{c} \scriptsize$(\mathfrak{r}_{2}',J_1(a,b))$ \\ \scriptsize$ (a,b)\neq (0,\pm1) $ \end{tabular}}] (a2) at (2,+0) {};
\node[circle,fill,inner sep=0pt,minimum size=3pt,label=above:{$(\mathfrak{r}_{2}',J_2)$}] (a3) at (4,+0) {};
\node[circle,fill,inner sep=0pt,minimum size=3pt,label=above:{\begin{tabular}{c} \scriptsize$(\mathfrak{rr}'_{3,a}\times \mathbb{R},J_k)$ \\ \scriptsize$a>0,\, k=1,2$ \end{tabular}}] (a4) at (6,+0) {};
\node[circle,fill,inner sep=0pt,minimum size=3pt,label=right:{$(\mathfrak{rr}_{3,1}\times \mathbb{R},J)$}] (a5) at (8,+0) {};
\node[circle,fill,inner sep=0pt,minimum size=3pt,label=left:{$(\mathfrak{r}_{2}',J_1(0,-1))$}] (a6) at (0,-2) {};
\node[circle,fill,inner sep=0pt,minimum size=3pt,label=right:{$(\mathfrak{rh}_{3} \times \mathbb{R},J)$}] (a7) at (4,-2) {};
\node[circle,fill,inner sep=0pt,minimum size=3pt,label=below:{$(\mathfrak{a}_{4}, J)$}] (a8) at (2,-4) {};

\draw [->] (a1) edge (a7) (a2) edge (a7) (a3) edge (a7) (a4) edge (a7) (a5) edge (a7);
\draw [->] (a1) edge (a6);
\draw [->] (a4) edge (a8) (a6) edge (a8);

\end{tikzpicture}
\end{center}

\begin{center}
\begin{tikzpicture}

\node[circle,fill,inner sep=0pt,minimum size=3pt,label=above:{\begin{tabular}{c} \scriptsize$(\mathfrak{d}'_{4,\delta},J_1)$ \\ \scriptsize$ a>0 $ \end{tabular}}] (a1) at (0,+0) {};
\node[circle,fill,inner sep=0pt,minimum size=3pt,label=above:{\begin{tabular}{c} \scriptsize$(\mathfrak{d}'_{4,\delta},J_4)$ \\ \scriptsize$ a>0 $ \end{tabular}}] (a2) at (2,+0) {};
\node[circle,fill,inner sep=0pt,minimum size=3pt,label=above:{\begin{tabular}{c} \scriptsize$(\mathfrak{d}'_{4,\delta},J_2)$ \\ \scriptsize$ a>0 $ \end{tabular}}] (a3) at (4,+0) {};
\node[circle,fill,inner sep=0pt,minimum size=3pt,label=above:{\begin{tabular}{c} \scriptsize$(\mathfrak{d}'_{4,\delta},J_3)$ \\ \scriptsize$ a>0 $ \end{tabular}}] (a4) at (6,+0) {};
\node[circle,fill,inner sep=0pt,minimum size=3pt,label=above:{\begin{tabular}{c} \scriptsize$(\mathfrak{d}_{4,a},J_k)$ \\ \scriptsize$a>\frac{1}{2}, a\neq1$\\ \scriptsize$k=1,2 $ \end{tabular}}] (a5) at (8,+0) {};
%%%%
\node[circle,fill,inner sep=0pt,minimum size=3pt,label=left:{\begin{tabular}{c} \scriptsize$(\mathfrak{r}'_{4,c,d},J_2)$ \\ \scriptsize$ c\neq -\frac{1}{2}$ \\ \scriptsize$  d>0$ \end{tabular}}] (a6) at (1,-2) {};
\node[circle,fill,inner sep=0pt,minimum size=3pt,label=left:{\begin{tabular}{c} \scriptsize$(\mathfrak{r}'_{4,c,d},J_1)$ \\ \scriptsize$  c\neq -\frac{1}{2}$ \\ \scriptsize$  d>0 $ \end{tabular}}] (a7) at (5,-2) {};
\node[circle,fill,inner sep=0pt,minimum size=3pt,label=above:{\begin{tabular}{c} \scriptsize$(\mathfrak{r}_{4,\alpha,\alpha},J)$ \\\scriptsize$ \alpha \in \mathbb{R} \setminus\{0\}$ \\ \scriptsize$ \alpha \neq \pm \frac{1}{2},1 $ \end{tabular}}] (a8) at (9,-2) {};
\node[circle,fill,inner sep=0pt,minimum size=3pt,label=right:{$(\mathfrak{rh}_{3}\times \mathbb{R},J)$}] (a9) at (5,-4) {};
\node[circle,fill,inner sep=0pt,minimum size=3pt,label=below:{$(\mathfrak{a}_{4}, J)$}] (a10) at (5,-6) {};
%%%%
\node[label={[rotate=90]above:{\scriptsize $(c,d)=(\frac{1}{2},\delta)$}}] (a11) at (1.25,-1) {};
\node[label={[rotate=90]above:{\scriptsize $(c,d)=(\frac{1}{2},\delta)$}}] (a11) at (5.25,-1) {};

%\node[label={[label distance=0.5cm,text depth=-1ex,rotate=-90]right:a long text}] at (2.5,1) {};

\draw [->] (a1) edge (a6) (a2) edge (a6);
\draw [->] (a3) edge (a7) (a4) edge (a7);
\draw [->] (a5) edge (a9);
\draw [->] (a6) edge (a9) (a7) edge (a9) (a8) edge (a9);
\draw [->] (a9) edge (a10);

\end{tikzpicture}
\end{center}

\begin{center}
\begin{tikzpicture}

\node[circle,fill,inner sep=0pt,minimum size=3pt,label=above:{$(\mathfrak{d}_{4,\sfrac{1}{2}},J_2)$}] (a1) at (0,+0) {};
\node[circle,fill,inner sep=0pt,minimum size=3pt,label=above:{$(\mathfrak{d}_{4,\sfrac{1}{2}},J_3)$}] (a2) at (2,+0) {};
\node[circle,fill,inner sep=0pt,minimum size=3pt,label=above:{$(\mathfrak{h}_{4},J)$}] (a3) at (4,+0) {};
%%%%%%
\node[circle,fill,inner sep=0pt,minimum size=3pt,label=left:{$(\mathfrak{r}_{4,\sfrac{1}{2},\sfrac{1}{2}},J)$}] (a4) at (1,-2) {};
\node[circle,fill,inner sep=0pt,minimum size=3pt,label=right:{$(\mathfrak{d}_{4,\sfrac{1}{2}},J_1)$}] (a5) at (3,-2) {};
%\node[label={[rotate=90]above:{$(\mathfrak{r}_{4,\sfrac{1}{2},\sfrac{1}{2}},J)$}}] (a4a) at (1.3,-1) {};
%\node[label={[rotate=90]above:{$(\mathfrak{d}_{4,\sfrac{1}{2}},J_1)$}}] (a5a) at (3.3,-1) {};
%%%%%%
\node[circle,fill,inner sep=0pt,minimum size=3pt,label=right:{$(\mathfrak{r}_{4,1},J)$}] (a6) at (5,-2) {};
\node[circle,fill,inner sep=0pt,minimum size=3pt,label=left:{$(\mathfrak{rh}_{3}\times \mathbb{R},J)$}] (a7) at (3,-4) {};
\node[circle,fill,inner sep=0pt,minimum size=3pt,label=right:{$(\mathfrak{r}_{4,1,1}, J )$}] (a8) at (5,-4) {};
\node[circle,fill,inner sep=0pt,minimum size=3pt,label=below:{$(\mathfrak{a}_4 , J )$}] (a9) at (3,-6) {};

\draw [->]  (a1) edge (a4) (a2) edge (a4) (a2) edge (a5) (a3) edge (a5) ;
\draw [->]  (a4) edge (a7) (a5) edge (a7) (a6) edge (a8);
\draw [->]  (a7) edge (a9) (a8) edge (a9);
\draw [->]  (a6) edge (a7);

\end{tikzpicture}
\end{center}